\documentclass[11pt]{amsart}
\usepackage[letterpaper,margin=1in]{geometry}
\usepackage[T1]{fontenc}
\usepackage{lmodern}
\usepackage{microtype}
\usepackage{amsmath,amssymb,amsthm,mathtools,bm,mathrsfs}
\usepackage{booktabs,tabularx,array}
\usepackage[dvipsnames]{xcolor}
\usepackage{hyperref}

\definecolor{navy}{HTML}{17365D}
\definecolor{teal}{HTML}{287D7D}
\definecolor{rust}{HTML}{A34A28}
\definecolor{softgray}{HTML}{F2F4F6}
\hypersetup{colorlinks=true,linkcolor=navy,citecolor=teal,urlcolor=rust,pdfauthor={},pdftitle={Free-boundary IPM equations with surface tension}}
\relax {}

\newcolumntype{Y}{>{\raggedright\arraybackslash}X}
\allowdisplaybreaks
\numberwithin{equation}{section}

\newtheorem{theorem}{Theorem}[section]
\newtheorem{proposition}[theorem]{Proposition}
\newtheorem{lemma}[theorem]{Lemma}
\newtheorem{corollary}[theorem]{Corollary}
\theoremstyle{definition}
\newtheorem{definition}[theorem]{Definition}
\theoremstyle{remark}
\newtheorem{remark}[theorem]{Remark}

\newcommand{\R}{\mathbb R}
\newcommand{\T}{\mathbb T}
\newcommand{\Om}{\Omega}
\newcommand{\Sig}{\Sigma}
\newcommand{\cO}{\mathcal O}
\newcommand{\DN}{\mathcal G}
\newcommand{\bG}{\mathbf G}
\newcommand{\bS}{\mathbf S}
\newcommand{\abs}[1]{\left|#1\right|}
\newcommand{\dd}{\,\mathrm d}
\DeclareMathOperator{\diver}{div}

\title[Free-boundary IPM equations]{On the free-boundary incompressible porous media equations with surface tension}

\author{Lizhe Wan}
\address{Beijing International Center for Mathematical Research, Peking University}
\curraddr{}
\email{wanlizhe@pku.edu.cn}

\author{Jiaqi Yang}
\address{School of Mathematics and Statistics, Northwestern Polytechnical University}
\curraddr{}
\email{yjqmath@nwpu.edu.cn, yjqmath@163.com}


\begin{document}

\begin{abstract}
We study the free-boundary incompressible porous media equations with surface tension in dimension $1+d$.
For fixed positive surface tension, we establish local well-posedness for arbitrary initial interfaces in $H^{s+1}$ and density perturbations in $H^s$, with $s>(d+3)/2$.
When $s>(d+5)/2$, we justify the local vanishing-surface-tension limit for stable interfaces and sufficiently small density perturbations.
In a finite-depth strip over $\R^d$ with a flat bottom and fixed positive surface tension, we prove global well-posedness and uniform anisotropic energy bounds for small compatible perturbations of an equilibrium with strictly stable stratification, under a horizontal low-frequency condition.
\end{abstract}

\maketitle

\setcounter{tocdepth}{1}
\tableofcontents

\section{Introduction}

We consider an incompressible fluid in a porous medium whose free surface is a graph and whose nonconstant density is transported by the flow.
The pressure at the free surface is prescribed by its curvature, and the resulting parabolic interface evolution is coupled to density transport in the moving fluid domain.
For this free-boundary incompressible porous media (IPM) system, we establish local well-posedness at fixed positive surface tension and justify the vanishing-surface-tension limit for stable interfaces and small density perturbations.
In a finite-depth strip with a flat bottom, we also prove global well-posedness near a strictly stably stratified equilibrium for small compatible data satisfying a horizontal low-frequency condition.

\subsection{The free-boundary IPM equations}

Let $d\geq1$, $\cO$ be $\R^d$ or $\T^d$, and let the fluid occupy
\begin{equation*}
  \Om_{f(t)}=\{(x,y)\in\cO\times\R:b(x)<y<f(t,x)\}.
\end{equation*}
The free upper interface is
\begin{equation*}
\Sig_{f(t)}=\{y=f(t,x)\}.
\end{equation*}
In finite depth, $b$ is a fixed smooth bottom profile strictly below the surface $\Sig_{f(t)}$.
In infinite depth, we set $b=-\infty$ and omit the bottom boundary condition.
After normalizing permeability, viscosity, and gravity, the equations read
\begin{equation*}
\begin{cases}
u+\nabla_{x,y}p=-\rho e_{d+1},\qquad \diver_{x,y}u=0 & \text{in }\Om_{f(t)},\\
\partial_t\rho+u\cdot\nabla_{x,y}\rho=0 & \text{in }\Om_{f(t)},\\
u\cdot\nu_b=0 & \text{on }y=b(x),\\
p=\sigma H(f),\qquad \partial_t f=u\cdot N_f & \text{on }\Sig_{f(t)},
\end{cases}
\end{equation*}
where $\sigma\geq0$, $N_f=(-\nabla f,1)$, and the curvature convention is
\begin{equation*}
  H(f)=-\diver_x\left(\frac{\nabla f}{\sqrt{1+\abs{\nabla f}^2}}\right).
\end{equation*}
Here $u$, $p$, and $\rho$ denote velocity, pressure, and density of the fluid.
$e_{d+1}$ is the upward vertical unit vector and $\nu_b$ is the unit normal to the bottom.
The first two equations express Darcy's law, incompressibility, and transport of the density.
In finite depth, the bottom is impermeable.
The dynamic condition $p=\sigma H(f)$ prescribes the pressure relative to the constant exterior pressure, while the kinematic condition states that the surface moves with the fluid.
In graph coordinates, the latter reads $\partial_t f=u_{d+1}-u_x\cdot\nabla f$ on $\Sig_{f(t)}$.

\subsection{The reduced boundary--bulk system}\label{sec:reduced}

We use the boundary--bulk reduction of Latocca and Nguyen~\cite{LatoccaNguyen2025} for the gravity-driven problem, modifying the dynamic condition to include capillarity.
Fix a smooth background profile $\gamma(y)$, we write the density as $\rho=\gamma+g$, where $g$ is the time-dependent perturbation:
\begin{equation*}
  \rho(t,x,y)=\gamma(y)+g(t,x,y),  \qquad  \Gamma(y)=\int_0^y\gamma(\eta)\,\dd\eta.
\end{equation*}
The normalization $\Gamma(0)=0$ allows the usual Sobolev composition estimates for $\Gamma(f)$ on $\mathcal O=\R^d$.
With the shifted pressure $q=p+\Gamma(y)$, Darcy's law and the boundary conditions become
\begin{equation}\label{eq:qproblem}
  u=-\nabla q-ge_{d+1},\qquad
  \Delta q=-\partial_y g,\qquad
  q\big|_{\Sig_f}=\sigma H(f)+\Gamma(f),\qquad
  \partial_{\nu_b}q=-\nu_{b,d+1}g.
\end{equation}
For a fixed interface, one can define the harmonic and density-dependent pressure operators through their boundary-value problems.
For a surface datum $h$ on $\Sig_f$, we define the operator $\bG[f]h$ by
\begin{equation*}
\begin{cases}
 \Delta_{x,y}\phi_h=0 & \text{in }\Om_f,\\
 \phi_h=h & \text{on }\Sig_f,\\
 \partial_{\nu_b}\phi_h=0 & \text{on }y=b(x).
\end{cases}
\qquad
 \bG[f]h:=\nabla_{x,y}\phi_h.
\end{equation*}
For a bulk datum $k$ on $\Om_f$, we define the operator $\bS[f]k$ by
\begin{equation*}
\begin{cases}
 \Delta_{x,y}\psi_k=-\partial_y k & \text{in }\Om_f,\\
 \psi_k=0 & \text{on }\Sig_f,\\
 \partial_{\nu_b}\psi_k=-\nu_{b,d+1}k & \text{on }y=b(x).
\end{cases}
\qquad
 \bS[f]k:=\nabla_{x,y}\psi_k.
\end{equation*}
In infinite depth, we omit the bottom conditions and solve the elliptic problems in the finite-energy class, with $\nabla\phi_h,\nabla\psi_k\in L^2(\Omega_f)$.
The Dirichlet--Neumann operator is the normal trace of $\bG[f]h$:
\[
  \DN[f]h=N_f\cdot\bG[f]h\big|_{\Sig_f}.
\]
Decomposing the pressure into its harmonic and density-dependent parts gives the expression of the velocity
\begin{equation*}
  u=-\bG[f]\bigl(\sigma H(f)+\Gamma(f)\bigr)-\bS[f]g-ge_{d+1}.
\end{equation*}
The interface and density perturbation therefore satisfy the system
\begin{equation} \label{e:coupled}
\begin{cases}
\partial_t f=-\DN[f]\bigl(\sigma H(f)+\Gamma(f)\bigr)-\bigl(N_f\cdot\bS[f]g+g\bigr)\big|_{\Sig_f},\\
\partial_tg+u\cdot\nabla g+\gamma'(y)u_{d+1}=0.
\end{cases}
\end{equation}
We refer to \cite[Section~1.1]{LatoccaNguyen2025} for the derivation of this formulation.
For $\sigma>0$, the interface equation is third-order parabolic at leading order, while the density is transported in the moving $(d+1)$-dimensional fluid domain.

To compare solutions with different interfaces, we pull the bulk variables back to a common domain.
Set $\mathcal S=\cO\times(-1,0)$ in finite depth and $\mathcal S=\cO\times(-\infty,0)$ in infinite depth, and choose a regularized graph diffeomorphism
\begin{equation*}
  \Phi_f(t,x,z)=(x,\varrho_f(t,x,z)),\qquad
  \Phi_f(\mathcal S)=\Om_{f(t)}.
\end{equation*}
We use the regularized extension of \cite[Lemma~3.3]{LatoccaNguyen2025}, whose explicit formula is given in \eqref{eq:explicit-flattening}.
By construction, $\varrho_f$ agrees with $f$ on $\{z=0\}$ and (in finite depth) with $b$ on $\{z=-1\}$.
For a sufficiently small common flattening parameter $\delta>0$, the vertical derivative $\partial_z\varrho_f$ is bounded below by a positive constant and $\varrho_f$ obeys tame Sobolev estimates with a $1/2$-derivative gain from the boundary to the bulk (see Lemma~\ref{app:ln-flattening}).
Define
\[
  \widetilde g=g\circ\Phi_f,\qquad \widetilde u=u\circ\Phi_f.
\]
The density equation becomes
\begin{equation}\label{eq:ALE}
 \partial_t\widetilde g+V_f\cdot\nabla_{x,z}\widetilde g =Q_f,\quad Q_f=-\gamma'(\varrho_f)\widetilde u_{d+1},\quad
 V_f=\left(\widetilde u_x,
 \frac{\widetilde u_{d+1}-\nabla_x\varrho_f\cdot\widetilde u_x -\partial_t\varrho_f}{\partial_z\varrho_f}\right).
\end{equation}
The kinematic condition and, in finite depth, the bottom impermeability condition imply that the transformed velocity $V_f$ has vanishing normal component on $\partial\mathcal S$.
Thus the transport estimates for the pulled-back density on $\mathcal S$ have no boundary flux.
The interface equation in \eqref{e:coupled}, with $g=\widetilde g\circ\Phi_f^{-1}$, and the transport equation \eqref{eq:ALE} form the system
\begin{equation*} \tag{FIPM}\label{e:IPM}
\begin{cases}
\partial_t f=-\DN[f]\bigl(\sigma H(f)+\Gamma(f)\bigr)
  -\bigl(N_f\cdot\bS[f]g+g\bigr)\big|_{\Sig_f},\\
\partial_t\widetilde g+V_f\cdot\nabla_{x,z}\widetilde g =Q_f.
\end{cases}
\end{equation*}
We can recover the pressure from \eqref{eq:qproblem} and $p=q-\Gamma(y)$, the velocity from Darcy's law, and the density from $\rho=\gamma+g$.

\subsection{Previous results}
When $\rho\equiv\rho_c$ is constant, the density transport equation is identically satisfied, and \eqref{e:coupled} reduces to the one-phase Muskat problem with surface tension.
Choosing $\gamma\equiv\rho_c$ and $g=0$ in the reduced formulation gives the interface equation
\[
 \partial_t f+\DN[f]\bigl(\sigma H(f)+\rho_c f\bigr)=0.
\]
The classical Muskat problem describes an interface separating fluids of constant density in a porous medium.
Its well-posedness depends on the interface geometry, the Rayleigh--Taylor condition, and the smoothing available in the evolution.
Using contour-dynamics formulations, C\'ordoba--C\'ordoba--Gancedo \cite{CordobaCordobaGancedo2011,CordobaCordobaGancedo2013} established local well-posedness in the stable regime in two and three dimensions.
In two dimensions without surface tension, Cheng--Granero-Belinch\'on--Shkoller \cite{ChengGraneroShkoller2016} proved local well-posedness for the one-phase problem with arbitrary $H^2$ initial interfaces via a fixed-domain energy method.
Nguyen--Pausader \cite{NguyenPausader2020} introduced a Dirichlet--Neumann and paradifferential approach, proving local well-posedness in subcritical Sobolev spaces in arbitrary dimension, with general viscosity contrast and Lipschitz bottoms.
See also Alazard--Lazar \cite{AlazardLazar2020} for a paralinearization in homogeneous Sobolev spaces in 2D.
For the equal-viscosity two-phase problem, Alazard and Q.-H. Nguyen \cite{AlazardNguyen2022,AlazardNguyen2023} developed critical Cauchy theories at the endpoint $H^{3/2}(\R)$ in 2D and for critical $\dot H^2(\R^2)\cap W^{1,\infty}(\R^2)$ data in 3D.

Closely related free-boundary Darcy flows arise in the Hele--Shaw problem.
Ambrose \cite{Ambrose2004} proved local well-posedness for two-phase Hele--Shaw flow without surface tension in the stable regime.
Global existence for small data was established by Constantin--Pugh \cite{ConstantinPugh1993} and Siegel--Caflisch--Howison \cite{SiegelCaflischHowison2004}. The latter also constructed finite-time singularities in the unstable regime.
Using a Dirichlet--Neumann formulation, Alazard--Meunier--Smets \cite{AlazardMeunierSmets2020} derived new Lyapunov functionals and applied them to the Cauchy problem.

For $\sigma>0$, capillarity gives third-order parabolic smoothing at the free boundary.
Classical and abstract parabolic theories for Hele--Shaw and Muskat flows with surface tension were developed by Escher--Simonett \cite{EscherSimonett1997} and Escher--Matioc \cite{EscherMatioc2011}.
Using paradifferential methods, Nguyen \cite{Nguyen2020} established large-data local well-posedness in subcritical Sobolev spaces for one- and two-phase Muskat flows with surface tension in arbitrary dimension.
The composition $\DN[f]H(f)$ is elliptic of order three and provides a $3/2$-derivative gain in $L^2_t$ for the interface.
This smoothing permits local well-posedness for arbitrary initial interfaces without a Rayleigh--Taylor stability condition.

More recently, Chen, Hu, and Q.-H. Nguyen \cite{ChenHuNguyen2026I} developed Schauder estimates for nonlocal quasilinear evolution equations, with applications including the Muskat problem with surface tension.
In a companion paper \cite{ChenHuNguyen2026II}, they proved local well-posedness for large graph interfaces in a log-critical class for the two-phase problem with unequal viscosities, arbitrary density contrast, and positive separation from rigid boundaries, without a Rayleigh--Taylor sign condition.
In these results, the density is constant in each phase; here it is a transported variable in the fluid domain.

In the stable regime without surface tension, global existence for small or critical data was studied by Constantin--C\'ordoba--Gancedo--Strain \cite{ConstantinCordobaGancedoStrain2013}, by Gancedo--Garc\'ia-Ju\'arez--Patel--Strain \cite{GancedoGarciaJuarezPatelStrain2019} in Wiener spaces for the viscosity-jump case, and by Nguyen \cite{NguyenBesov2021} in the critical Besov space $\dot B^1_{\infty,1}$.
For the 2D periodic one-phase problem, Dong--Gancedo--Nguyen \cite{DongGancedoNguyen2023} proved global well-posedness for arbitrary Lipschitz initial interfaces, showing that the strong solution coincides with the unique viscosity solution.

The vanishing-surface-tension limit ($\sigma\to 0^+$) links the capillary and gravity-driven regimes.
This singular limit was justified for smooth data by Ambrose \cite{Ambrose2014} in 2D and Ambrose--Liu \cite{AmbroseLiu2020} in 3D.
At subcritical Sobolev regularity, Flynn--Nguyen \cite{FlynnNguyen2020} proved quantitative convergence for stable initial interfaces, obtaining an $O(\sqrt{\sigma})$ rate in $C_TH^{s-1}$ (one derivative below the initial regularity $H^s$) and an $O(\sigma)$ rate in $C_TH^{s-2}$ when the initial curvature belongs to $L^2$.

On a fixed domain, the IPM equation with nonconstant density is an active-scalar equation: Darcy's law determines the velocity from the density.
C\'ordoba--Gancedo--Orive \cite{CordobaGancedoOrive2007} proved local well-posedness and continuation criteria in 2D.
For perturbations of a stably stratified background profile, long-time and global stability results were established by Elgindi \cite{Elgindi2017} on $\R^2$ or $\T^2$ and by Castro--C\'ordoba--Lear \cite{CastroCordobaLear2019} in a periodic strip with impermeable boundaries. Further results on stratified stability were obtained by
Bianchini--Crin-Barat--Paicu \cite{BianchiniCrinBaratPaicu2024}
on $\R^2$ using anisotropic estimates.
Park \cite{Park2025} developed an energy approach in a fixed
periodic channel, while Bianchini--Jo--Park--Wang
\cite{BianchiniJoParkWang2025} established asymptotic stability
in $H^k$, $k>2$, on $\T\times\R$.
For stability and instability of general stratified steady states
in a fixed periodic channel, see also
Banihashemi--Mohammadkhani--Nguyen
\cite{BanihashemiMohammadkhaniNguyen2026}.
For weak solutions and unstable configurations, nonuniqueness via convex integration was shown by C\'ordoba--Faraco--Gancedo \cite{CordobaFaracoGancedo2011} and Sz\'ekelyhidi \cite{Szekelyhidi2012}, while small-scale formation was studied by Kiselev--Yao \cite{KiselevYao2023}.

Latocca--Nguyen \cite{LatoccaNguyen2025} studied the free-boundary IPM system, where a moving interface is coupled to a nonconstant transported bulk density.
For $\sigma=0$, they proved local well-posedness under a Rayleigh--Taylor stability condition and a smallness assumption on the density perturbation.
Propagating the Sobolev regularity of the bulk density in this coupled setting requires $L^2_t$ estimates for the full bulk velocity in addition to interface energy bounds.
In the present work, we combine these fixed-domain bulk elliptic estimates with the third-order capillary smoothing of the free surface.

\subsection{Main results}

We write $D=d+1$ for the bulk dimension.
The interface $f=f(t,x)$ is defined on $[0,T]\times\cO$,
whereas the fixed-domain bulk variables
$\widetilde g,\widetilde u,V_f$ depend on $(t,x,z)\in[0,T]\times\mathcal S$.
Thus $H^r(\cO)$ measures regularity in $x$, while
$H^r(\mathcal S)$ is the isotropic Sobolev space in $(x,z)$.

For $\mathcal D=\cO$ or $\mathcal S$, we use
\[
 C_TH^r(\mathcal D)=C([0,T];H^r(\mathcal D)),
 \qquad
 L_T^pH^r(\mathcal D)=L^p([0,T];H^r(\mathcal D)).
\]
For interface functions, we define
\[
 Z_T^r=C_TH^r(\cO)\cap L_T^2H^{r+3/2}(\cO),
 \qquad
 \|f\|_{Z_T^r}
 =\|f\|_{C_TH^r(\cO)}
  +\|f\|_{L_T^2H^{r+3/2}(\cO)}.
\]
We omit the spatial domain from the notation when it is clear from context.

We take the initial interface and density perturbation in
\[
 (f_0,\widetilde g_0)
 \in H^{s+1}(\cO)\times H^s(\mathcal S),
 \qquad
 \widetilde g_0=g_0\circ\Phi_{f_0}.
\]
We compare data on different physical domains in this product space
after pulling them back to $\mathcal S$.
We use common flattening parameters on bounded sets of interfaces
with, in finite depth, a common positive depth bound.

\subsubsection{Local well-posedness for fixed surface tension}

For the next two theorems, let $d\ge1$ and
$\cO=\R^d$ or $\T^d$.
In finite depth, assume
$b=-\ell_0+b_0$, where $\ell_0>0$ and
$b_0\in\bigcap_{k\ge0}H^k(\cO)$.
In infinite depth, we use the finite-energy pressure normalization specified above.
At Sobolev index $s$, assume that $\gamma$ is smooth with bounded derivatives through order $\lceil s\rceil+3$ on the relevant vertical range.

\begin{theorem}[Large-data local well-posedness]
\label{thm:local-target}
Let $\sigma>0$ and $s>(d+3)/2$.
For any initial data
\[
 f_0\in H^{s+1}(\cO),\qquad g_0\in H^s(\Omega_{f_0}),
\]
satisfying $\inf_x(f_0-b)\ge2h_0>0$ in finite depth, there exist $T>0$ and a unique solution of \eqref{e:IPM} with
\begin{equation}\label{eq:main-regularity}
f\in Z_T^{s+1},\qquad
 \widetilde g\in C_TH^s(\mathcal S),\qquad 
 \widetilde u,V_f\in L_T^2H^s(\mathcal S),\qquad
 f_t\in L_T^2H^{s-1/2}(\cO).
\end{equation}
The solution depends continuously on
$(f_0,\widetilde g_0)$ in
$H^{s+1}(\cO)\times H^s(\mathcal S)$,
with values in the spaces in \eqref{eq:main-regularity}.
\end{theorem}

On bounded subsets of $H^{s+1}\times H^s$,
with a common positive depth bound in finite depth,
the solutions have a common existence interval.
On this interval, the map
$(f_0,\widetilde g_0)\mapsto(f,\widetilde g)$
is locally Lipschitz for the lower-order topologies
\[
 H^s\times H^{s-1}
 \longrightarrow Z_T^s\times C_TH^{s-1}.
\]
The total density also preserves its lower bound:
if $\rho_0\ge c_\rho>0$, then
$\rho(t)\ge c_\rho$ on $[0,T]$.

\subsubsection{Vanishing surface tension}

We introduce the stability functional
\[
 \mathcal T(f)=\gamma(f)-\DN[f]\Gamma(f)
\]
as in \cite[(1.29), Theorem~1.4]{LatoccaNguyen2025}.
When $\gamma$ is a positive constant,
$\mathcal T(f)=\gamma(1-\DN[f]f)>0$
by \cite[Proposition~4.3]{NguyenPausader2020}.

\begin{theorem}[Vanishing surface tension limit]
\label{thm:vanishing-surface-tension}
Let $s>(d+5)/2$ and fix $R_0,a_0,h_0>0$.
There exist $\varepsilon,T,C>0$, depending only on $d,s,R_0,a_0,h_0$ and the fixed geometry and background, but independent of $\sigma\in[0,1]$, with the following properties.

Suppose the common initial data satisfy
\[
 \|f_0\|_{H^{s+1}}\le R_0,\qquad
 \|\widetilde g_0\|_{H^s}\le\varepsilon,\qquad
 \mathcal T(f_0)\ge2a_0,
\]
and $f_0-b\ge2h_0$ in finite depth.
For every $\sigma\in[0,1]$, the solution of
\eqref{e:IPM} exists uniquely on $[0,T]$ and satisfies
$\mathcal T(f^\sigma)\ge a_0$ and, in finite depth,
$f^\sigma-b\ge h_0$.
For $\sigma>0$, it has the regularity
\eqref{eq:main-regularity}; for $\sigma=0$,
\[
 f^0\in C_TH^s\cap L_T^2H^{s+1/2},
 \qquad \widetilde g^0\in C_TH^s.
\]
Uniformly for $\sigma\in[0,1]$, we have the upper bound
\[
\begin{aligned}
 \|f^\sigma\|_{C_TH^s}^2
 +\|\widetilde g^\sigma\|_{C_TH^s}^2
 +\sigma\|f^\sigma\|_{C_TH^{s+1}}^2
 +\|f^\sigma\|_{L_T^2H^{s+1/2}}^2
 +\sigma\|f^\sigma\|_{L_T^2H^{s+3/2}}^2 +\sigma^2\|f^\sigma\|_{L_T^2H^{s+5/2}}^2 \le C,
\end{aligned}
\]
with the $\sigma$-weighted terms omitted at $\sigma=0$.
Moreover,
\begin{equation}\label{eq:main-tension-rate}
\begin{aligned}
 \|f^\sigma-f^0\|_{C_TH^{s-1}}
 +\|f^\sigma-f^0\|_{L_T^2H^{s-1/2}}
 +\|\widetilde g^\sigma-\widetilde g^0\|_{C_TH^{s-1}}
 \le C\sqrt{\sigma},
 \qquad 0<\sigma\le1.
\end{aligned}
\end{equation}
\end{theorem}

The associated velocity, transport field, and density
source also satisfy
\[
 \|\widetilde u^\sigma\|_{L_T^2H^s}
 +\|V_{f^\sigma}\|_{L_T^2H^s}
 +\|Q_{f^\sigma}\|_{L_T^2H^s}
 \le C,
 \qquad 0\le\sigma\le1.
\]

\subsubsection{Global well-posedness near a strictly stratified equilibrium}\label{sec:global-intro}

For the global result, we restrict to $\mathcal O=\R^d$, a flat bottom $b=-1$, and fixed positive surface tension $\sigma>0$.
The additional assumptions and notation in this subsubsection apply to Theorem~\ref{thm:global} and Sections~\ref{sec:global-linear-estimates}--\ref{sec:global-estimates}.

We consider perturbations of a hydrostatic equilibrium with strictly stable stratification.
Fix a smooth background density profile $\gamma$ on an open neighborhood $I_0$ of $[-1,0]$ such that
\begin{equation}\label{eq:global-profile}
 \gamma(0)=c>0,\qquad 0<c_*\le-\gamma'(y)\le C_*\qquad\text{for all }y\in I_0.
\end{equation}
In the global argument, $f$ remains sufficiently small in $L^\infty$ that $[-1,f(t,x)]$ is contained in $I_0$.
In this subsubsection and Sections~\ref{sec:global-linear-estimates}--\ref{sec:global-estimates}, constants may depend on $\gamma$ and on $d,m,c$, and $\sigma$.
We also require the total density to be constant on the free surface and on the bottom:
\begin{equation}\label{eq:global-boundary-density}
 \rho(t,x,f(t,x))=c,\qquad \rho(t,x,-1)=\gamma(-1).
\end{equation}
Both boundaries are material surfaces: fluid particles
initially on either boundary remain on that boundary.
Hence the constant boundary densities are preserved by transport.
At $\{y=-1\}$, the second condition in \eqref{eq:global-boundary-density} reads $g|_{y=-1}=0$.

Set $\mathcal S=\R^d\times I$, where $I=(-1,0)$.
Whenever $\rho(t,x,\cdot)$ is strictly decreasing and satisfies \eqref{eq:global-boundary-density}, we define the semi-Lagrangian density coordinate map $\Phi_\eta(t,\cdot):\mathcal S\to\Om_{f(t)}$ by
\begin{equation}\label{eq:global-density-map}
 \rho(t,x,z+\eta(t,x,z))=\gamma(z),\qquad
 \Phi_\eta(x,z)=(x,z+\eta(t,x,z)).
\end{equation}
By \eqref{eq:global-boundary-density}, the function $\eta$ and its Jacobian $J_\eta$ satisfy
\[
 \eta|_{-1}=0,\qquad \eta|_{z=0}=f,\qquad
 J_\eta:=1+\eta_z=\frac{\gamma'(z)}{\rho_y\circ\Phi_\eta}>0.
\]
We denote the velocity pulled back to the density coordinates by $U=(U_x,U_{d+1})=u\circ\Phi_\eta$.
We measure $\eta$ and $U$ in the anisotropic Sobolev space $X_m$, with norm
\begin{equation}\label{eq:global-Xm}
 \|h\|_{X_m}^2=\|h\|_{L^2(\mathcal S)}^2
       +\|h_{zz}\|_{L^2(\mathcal S)}^2+\||D_x|^m h\|_{L^2(\mathcal S)}^2.
\end{equation}
For vector-valued functions, the squared $X_m$ norm is the sum of the squared component norms.
To control horizontal transport and boundary traces, we estimate the time integral of the velocity in the Fourier--Wiener norm
\begin{equation}\label{eq:global-Wiener}
 \mathcal W_d(h)=\inf_{\widetilde h|_{\mathcal S}=h}
 \int_{\R^d}\int_\R(1+|\xi|)|\widehat{\widetilde h}(\xi,\zeta)|\dd\zeta\dd\xi,
\end{equation}
where the infimum is taken over all extensions $\widetilde h$ of $h$ to $\R^d\times\R$.
The norm $\mathcal W_d(h)$ controls the $L^\infty(\mathcal S)$ norms of $h$ and $\nabla_xh$, as well as their pointwise traces on $\{z=0\}$ and $\{z=-1\}$.
To obtain $L^1_t$ integrability of $\mathcal W_d(U(t))$, we impose the low-frequency condition
\begin{equation}\label{eq:global-low-norm}
 \mathcal B_{\rm low}(h)=\int_{0<|\xi|<1} |\xi|^{-1}\|\widehat h(\xi,\cdot)\|_{H^2(I)}\dd\xi.
\end{equation}

\begin{definition}[Compatible initial data]\label{def:global-data}
Let $m\ge2d+5$ be an integer.
We call $(f_0,\widetilde g_0)\in H^{m+1}(\R^d)\times H^m(\mathcal S)$ a \emph{compatible initial datum} if the following conditions hold.
The initial fluid depth $1+f_0(x)$ is strictly positive.
The initial density $\rho_0=\gamma+\widetilde g_0\circ\Phi_{f_0}^{-1}$ satisfies \eqref{eq:global-boundary-density} and the strict stratification bound $-\rho_{0,y}\ge c_*/2$.
The combined size
\begin{equation}\label{eq:global-epsilon}
 \varepsilon=\|f_0\|_{H^{m+1}}
       +\|\widetilde g_0\|_{H^m}+\mathcal B_{\rm low}(\eta_0)
\end{equation}
is finite, where $\eta_0$ is defined by \eqref{eq:global-density-map} at $t=0$.
\end{definition}

\begin{theorem}[Global well-posedness]\label{thm:global}
Let $d\ge1$, $\sigma>0$, and let $m\ge2d+5$ be an integer.
Assume that $\gamma$ satisfies \eqref{eq:global-profile}.
There exist $\varepsilon_*,C>0$ such that every
compatible initial datum in Definition~\ref{def:global-data}
with $\varepsilon\le\varepsilon_*$ generates a unique
global strong solution of \eqref{e:IPM}.
The fluid depth and density remain positive, and
$\rho_y<0$ for all time.
Moreover,
\begin{equation}\label{eq:global-estimate}
\sup_{t\ge0}\left(
   \|\eta(t)\|_{X_m}^2+\|f(t)\|_{H^{m+1}}^2
 \right)
 +\int_0^\infty\|U(t)\|_{X_m}^2\,\dd t + \varepsilon\int_0^\infty\mathcal W_d(U(t))\,\dd t
 \le C\varepsilon^2.
\end{equation}

On every finite interval $[0,T]$, the solution satisfies
\begin{equation}\label{eq:global-local-class}
f\in Z_T^{m+1},\qquad
 \widetilde g\in C_TH^m(\mathcal S),\qquad \widetilde u,V_f\in L_T^2H^m(\mathcal S),\qquad
 f_t\in L_T^2H^{m-1/2}.
\end{equation}
On the compatible small-data class, the solution depends
continuously on $(f_0,\widetilde g_0)$ in
$H^{m+1}\times H^m$, with values in the spaces
in \eqref{eq:global-local-class}.
\end{theorem}

The weight $|\xi|^{-1}$ in \eqref{eq:global-low-norm} comes from integrating the horizontal velocity in time.
The following estimates give sufficient conditions for $\mathcal B_{\rm low}(h)$ to be finite:
\[
 \mathcal B_{\rm low}(h)\le C
 \begin{cases}
  \|h\|_{L_x^2H_z^2}, & d\ge3,\\
  \|h\|_{L_x^1H_z^2}, & d=2,\\
  \displaystyle\int_{\R}|x|\|h(x,\cdot)\|_{H_z^2}\dd x,
  & d=1,\quad \displaystyle\int_{\R}h(x,\cdot)\dd x=0.
 \end{cases}
\]
The first two bounds follow from the Cauchy--Schwarz inequality and the Fourier $L^\infty$ bound.
In dimension one, the cancellation gives
\[
 \|\widehat h(\xi,\cdot)\|_{H_z^2} =\left\|\int_{\R}(e^{-ix\xi}-1)h(x,\cdot)\dd x\right\|_{H_z^2} \le |\xi|\int_{\R}|x|\|h(x,\cdot)\|_{H_z^2}\dd x.
\]

\subsection{Proof strategy and organization}

The local and vanishing-surface-tension arguments use the same formulation on a fixed domain, but rely on different estimates for the interface.
The regularized flattening map of \cite{AlazardBurqZuily2014,NguyenPausader2020,LatoccaNguyen2025} gains half a derivative from the boundary to the bulk; this gain allows us to estimate the pressure and compare the velocities associated with different interfaces.
The kinematic condition makes the transformed transport field tangent to the fixed boundary, hence the density estimates have no boundary flux.

For fixed $\sigma>0$, the capillary operator has a positive third-order principal symbol.
The reduction in \cite{Nguyen2020} and the paradifferential calculus of \cite{BCD2011} yield an $H^{s+1}$ interface energy with $L_T^2H^{s+5/2}$ dissipation.
Together with the density bound, the elliptic estimates control $\widetilde u$ in $L_T^2H^s$.
The interface equation and the regularized extension also control the time derivative of the coordinate map, giving the same bound for $V_f$.
These estimates close the density transport estimate on a sufficiently short interval, with no stability assumption on the initial interface.

We choose the approximate transport fields to remain tangent to the fixed boundary while the interface and density are updated at different stages of the iteration.
The iterates satisfy uniform bounds on a common lifespan.
Comparisons one derivative below the initial regularity then give a contraction and identify the limit as a solution of the coupled system.
To recover continuous dependence at the top indices, we smooth a reference transport solution, in the spirit of Bona--Smith \cite{BonaSmith1975}, and combine lower-order stability with estimates for the smoothing error.

For the vanishing-tension limit, the fixed-$\sigma$
estimates must be replaced by bounds that remain uniform
as $\sigma\to0$.
Stability supplies first-order surface dissipation
that remains coercive in this limit
\cite{LatoccaNguyen2025,FlynnNguyen2020},
while smallness of the density perturbation allows
its contribution to the interface energy to be absorbed.
To control the bulk velocity uniformly, we supplement
the unweighted $H^s$ energy with a higher surface energy
weighted by $\sigma$.
The additional dissipation controls the capillary
contribution to the velocity.
Together with the tame elliptic and transport estimates,
this gives a common lifespan independent of $\sigma$.

On this interval, we compare with the zero-surface-tension
solution given by the stable local theory of
Latocca--Nguyen~\cite[Theorem~1.4]{LatoccaNguyen2025}.
We treat the capillary terms as forcing in the surface
and density difference equations.
The weighted estimates bound their contributions
by $C\sqrt{\sigma}$.
Coupling the two difference estimates then yields
the convergence rate in \eqref{eq:main-tension-rate}.

In Section~\ref{sec:elliptic} and Section~\ref{sec:estimates}, we establish the fixed-domain elliptic estimates, and the capillary and transport estimates respectively.
Then in Section~\ref{sec:construction}, we prove local well-posedness and the continuation criterion, while Section~\ref{sec:vanishing-surface-tension} is devoted to proving the uniform estimates and the vanishing-tension limit.
Next, in Section~\ref{sec:global-linear-estimates} we begin with the strategy for global well-posedness and develop the density-coordinate formulation and the forced linear velocity estimate.
In Section~\ref{sec:global-estimates}, we prove the nonlinear estimates and complete the proof of Theorem~\ref{thm:global}.
Finally, in Appendix~\ref{app:latocca-nguyen} we collect auxiliary estimates and the local theory without surface tension.

Throughout the paper, interface norms are taken over $\mathcal O$ and pulled-back bulk norms over $\mathcal S$.
For vector fields, $L^2$ and Sobolev norms use the Euclidean norm on the components; norms of tuples are understood as the sum of the norms of their entries.
Constants denoted by $C_R$ may change from line to line, with their precise dependence specified at the beginning of each section.

\section{Elliptic estimates on the fixed domain}\label{sec:elliptic}

\subsection{Flattening and the variational pressure equation}

We pull the pressure problem back to the fixed domain $\mathcal S$ using the graph map $\Phi_f(x,z)=(x,\varrho_f(x,z))$.
In finite depth, we decompose $b=-\ell_0+b_0$, where $-\ell_0$ is the flat reference bottom.
Following \cite[Lemma~3.3]{LatoccaNguyen2025}, recalled
in Lemma~\ref{app:ln-flattening}, set
\begin{equation}\label{eq:explicit-flattening}
\begin{aligned}
 \varrho_f&=\varrho_*+\mathscr E f,\\
 \varrho_*&=
 \begin{cases}
 z,&\text{in infinite depth},\\
 \ell_0z-ze^{-\delta(z+1)\langle D_x\rangle}b_0,
     &\text{in finite depth},
 \end{cases}\\
 \mathscr E h&=
 \begin{cases}
 e^{\delta z\langle D_x\rangle}h,
     &\text{in infinite depth},\\
 (1+z)e^{\delta z\langle D_x\rangle}h,
     &\text{in finite depth}.
 \end{cases}
\end{aligned}
\end{equation}
Thus $\varrho_f|_{z=0}=f$ and, in finite depth, $\varrho_f|_{z=-1}=b$.
Assume $s>(d+3)/2$, $\|f\|_{H^{s+1}}\le R$, and,
in finite depth, $f-b\ge h_0>0$.
For these interfaces, a common sufficiently small $\delta>0$ ensures
$J_f:=\partial_z\varrho_f\ge c_R>0$.

The extension gains one half derivative from the interface to the bulk.
Let $I=(-1,0)$ in finite depth and $I=(-\infty,0)$ in infinite depth.
Taking the Fourier transform in $x$ gives, for every integer $j\ge0$,
\[
\begin{aligned}
 \int_I
 \bigl|\partial_z^j\widehat{\mathscr E h}(\xi,z)\bigr|^2
 \,\dd z
 \le C_{j,\delta}\langle\xi\rangle^{2j}
       |\widehat h(\xi)|^2
       \int_I e^{2\delta z\langle\xi\rangle}\,\dd z\le C_{j,\delta}\langle\xi\rangle^{2j-1}
       |\widehat h(\xi)|^2.
\end{aligned}
\]
For integers $k\ge0$, Plancherel's theorem and summation
over $\partial_x^\alpha\partial_z^j$ with
$|\alpha|+j\le k$ give the corresponding $H^k$ bound.
Interpolation then yields
\[
 \|\mathscr E h\|_{H^r(\mathcal S)}
 \le C_{r,\delta}\|h\|_{H^{r-1/2}(\cO)},
 \qquad r\ge0.
\]
Since the bottom is fixed and $\mathscr Eh$ is linear in $h$, 
\[
 \varrho_{f_1}-\varrho_{f_2}=\mathscr E(f_1-f_2), \qquad \partial_t\varrho_f=\mathscr E f_t.
\]

The transformed pressure coefficient is
\begin{equation}\label{eq:matrix-A}
 A_f=J_fD\Phi_f^{-1}D\Phi_f^{-T}
 =\begin{pmatrix}
 J_fI_d&-\nabla_x\varrho_f\\
 -(\nabla_x\varrho_f)^T&
 (1+|\nabla_x\varrho_f|^2)/J_f
 \end{pmatrix}.
\end{equation}
Note that the flat reference matrices for $A_f$ and $D\Phi_f^{-T}$ are $A_*=M_*=I$ in infinite depth and
\[
 A_*=\operatorname{diag}(\ell_0I_d,\ell_0^{-1}),
 \qquad M_*=\operatorname{diag}(I_d,\ell_0^{-1})
\]
in finite depth.

In the estimates below, $C_R$ depends only on $d,s,R,\delta$ and, in finite depth, on $\ell_0$, $\|b_0\|_{H^{s+1}}$, and $h_0^{-1}$.
The extension estimate, Sobolev product estimates, and the lower bound for $J_f$ give
\begin{align}
 &\|A_f-A_*\|_{H^{s+1/2}} +\|D\Phi_f^{-T}-M_*\|_{H^{s+1/2}} \le C_R, \label{eq:geometry-single}\\
 &\|A_{f_1}-A_{f_2}\|_{H^r} +\|D\Phi_{f_1}^{-T}-D\Phi_{f_2}^{-T}\|_{H^r}
 \le C_R\|f_1-f_2\|_{H^{r+1/2}}, \quad r=s,s-1,
 \label{eq:geometry-difference}\\
 &\|\partial_t\varrho_f\|_{H^r}
 \le C_R\|f_t\|_{H^{r-1/2}}, \quad r=s,s-1.
 \label{eq:time-extension}
\end{align}
The same bounds for the flattening map and its inverse show that the physical and fixed-domain Sobolev norms are equivalent:
\begin{equation}\label{eq:sobolev-pullback}
 C_R^{-1}\|a\|_{H^r(\Omega_f)} \le \|a\circ\Phi_f\|_{H^r(\mathcal S)} \le C_R\|a\|_{H^r(\Omega_f)}, \qquad r=s,s-1.
\end{equation}

Let $v=q\circ\Phi_f$ and $\widetilde g=g\circ\Phi_f$ be the pullbacks of $q$ and $g$ to the fixed domain $\mathcal S$.
Using $J_fD\Phi_f^{-1}e_{d+1}=e_z$, we transform the weak pressure equation in \eqref{eq:qproblem} into
\begin{equation}\label{eq:pressure-weak}
 \int_{\mathcal S}A_f\nabla v\cdot\nabla\phi\,\dd x\dd z
 =-\int_{\mathcal S}\widetilde g\,\partial_z\phi\, \dd x\dd z
\end{equation}
for every energy test function $\phi$ that vanishes at $z = 0$.
The fixed-domain problem therefore becomes
\begin{equation}\label{eq:fixed-pressure}
 \operatorname{div}(A_f\nabla v)
 =-\partial_z\widetilde g, \qquad
 v|_{z=0}=h_f:=\sigma H(f)+\Gamma(f),
\end{equation}
with $(A_f\nabla v+\widetilde g e_z)\cdot e_z=0$
at $z=-1$ in finite depth.
The chain rule gives the pulled-back velocity:
\begin{equation}\label{eq:velocity-fixed}
 \widetilde u =-D\Phi_f^{-T}\nabla v-\widetilde g e_{d+1}.
\end{equation}

\subsection{Gradient estimates and elliptic comparison}

To control the full pressure gradient in bulk Sobolev spaces, we use localization and the constant-coefficient boundary estimates of \cite[Theorems~9.3 and~9.19]{Grubb2009}.
\begin{lemma}[Gradient estimate on the fixed domain]\label{lem:strip-elliptic}
Let $s>D/2+1$, $r\in\{s-1,s\}$, and let $A$ be a symmetric matrix on $\mathcal S$ such that
\[
 A(X)\zeta\cdot\zeta\ge\kappa|\zeta|^2,
 \qquad \|A-A_*\|_{H^{s+1/2}(\mathcal S)}\le R,
\]
for some $\kappa>0$ and all $\zeta\in\R^D$.
Let $v$ be the variational solution (as in \eqref{eq:pressure-weak}) of
\[
 \operatorname{div}(A\nabla v)=\operatorname{div}W,
 \qquad v|_{z=0}=h,
\]
with $(A\nabla v-W)\cdot e_z=0$ at $z=-1$ in finite depth.
Then
\begin{equation}\label{eq:strip-elliptic}
 \|\nabla v\|_{H^r(\mathcal S)} \le C_R\bigl( \|h\|_{H^{r+1/2}(\cO)}+\|W\|_{H^r(\mathcal S)} \bigr).
\end{equation}
Here $C_R$ also depends on $\kappa^{-1}$ and the fixed reference matrix $A_*$.
\end{lemma}
\begin{proof}
Fix $r\in\{s-1,s\}$.
For smooth coefficients and data, the estimate can be reduced to the case where $v|_{z=0}= 0$.
Indeed, one can choose an extension $Eh$ satisfying
\[
 Eh|_{z=0}=h,\qquad \|\nabla Eh\|_{H^q}\leq C\|h\|_{H^{q+1/2}}, \qquad q\in\{0,r\},
\]
and set
\[
 u=v-Eh,\qquad \widetilde W=W-A\nabla Eh.
\]
The function $u$ vanishes at $z = 0$ and satisfies
\[
 \operatorname{div}(A\nabla u)=\operatorname{div}\widetilde W,
\qquad
 \|\widetilde W\|_{H^r} \leq C_R\bigl(\|W\|_{H^r}+\|h\|_{H^{r+1/2}}\bigr),
\]
with $(A\nabla u-\widetilde W)\cdot e_z=0$ at the bottom.
Taking $u$ as a test function gives
\[
 \int_{\mathcal S}A\nabla u\cdot\nabla u
 =\int_{\mathcal S}\widetilde W\cdot\nabla u.
\]
By uniform ellipticity and the Cauchy--Schwarz inequality,
we obtain
\begin{equation}\label{eq:elliptic-energy}
 \|\nabla u\|_{L^2}\leq C_R\|\widetilde W\|_{L^2}, \qquad \|\nabla v\|_{L^2} \leq C_R\bigl(\|h\|_{H^{1/2}}+\|W\|_{L^2}\bigr).
\end{equation}

For a constant symmetric positive definite matrix $A_0$, the full-space resolvent bound and the half-space Poisson estimates in \cite[Theorems~9.3 and~9.19]{Grubb2009}, combined with Sobolev interpolation and a linear change of variables preserving the boundary plane, yield the following result:
For $w$ that solves
\[
 (1-\operatorname{div}(A_0\nabla))w =k-\operatorname{div}F,
\]
on the full space or a half-space, with either $w=0$ or $(A_0\nabla w-F)\cdot\nu=0$ on the boundary,
it satisfies the estimate
\begin{equation}\label{eq:constant-elliptic}
 \|w\|_{H^{r+1}} \leq C\bigl(\|F\|_{H^r}+\|k\|_{H^{r-1}}\bigr).
\end{equation}

To freeze the coefficients, fix $0<\varepsilon<\min\{1,r-D/2\}$.
We choose a uniformly locally finite cover by interior and boundary patches $U_j$, with cutoffs satisfying
$\sum_j\chi_j^2=1$. 
Set
\[
 A_j=A(X_j),\qquad w_j=\chi_j(u-c_j),\qquad a_j=\theta_j(A_j-A),
\]
where $\theta_j$ is supported in $U_j$ and equals one near $\operatorname{supp}\chi_j$.
Here $c_j=0$ on top boundary patches and is a local average otherwise. 
By taking the patches sufficiently small, we make the variation of $A$ small in $L^\infty$ and obtain the following $H^r$ bound:
\[
 \|a_j\|_{L^\infty}\leq\eta,\qquad \|a_j\|_{H^r}\leq C_R
\]
uniformly in $j$.
Then
\[
 A_j\nabla w_j-F_j =\chi_j(A\nabla u-\widetilde W),\qquad
 (1-\operatorname{div}(A_j\nabla))w_j=(w_j-S_j)-\operatorname{div}F_j,
\]
with
\[
 F_j=\chi_j\widetilde W +A(u-c_j)\nabla\chi_j+a_j\nabla w_j,\qquad S_j=(A\nabla u-\widetilde W)\cdot\nabla\chi_j.
\]
In particular, $w_j=0$ on top boundary patches, while
$(A_j\nabla w_j-F_j)\cdot\nu=0$ on the bottom.

Local Poincar\'e inequalities and Sobolev product estimates give
\[
\|w_j\|_{H^{r-1}} +\|\chi_j\widetilde W+A(u-c_j)\nabla\chi_j\|_{H^r} +\|S_j\|_{H^{r-1}}\leq C_R\bigl( \|\widetilde W\|_{H^r(U_j)} +\|\nabla u\|_{H^{r-\varepsilon}(U_j)}\bigr).
\]
Since $r-\varepsilon>D/2$, Sobolev embedding and the product estimate bound the coefficient error by
\[
 \|a_j\nabla w_j\|_{H^r}\leq C\|a_j\|_{L^\infty}\|\nabla w_j\|_{H^r} +C\|a_j\|_{H^r}\|\nabla w_j\|_{L^\infty}\leq C\eta\|\nabla w_j\|_{H^r} +C_R\|\nabla w_j\|_{H^{r-\varepsilon}}.
\]
Applying \eqref{eq:constant-elliptic} with $w=w_j$, $F=F_j$, and $k=w_j-S_j$, and using the preceding bounds, we obtain
\[
\begin{aligned}
 \|w_j\|_{H^{r+1}}
 \le C\bigl(\|F_j\|_{H^r}+\|w_j-S_j\|_{H^{r-1}}\bigr)\le C\eta\|w_j\|_{H^{r+1}}
 +C_R\bigl(\|\widetilde W\|_{H^r(U_j)} +\|\nabla u\|_{H^{r-\varepsilon}(U_j)}\bigr).
\end{aligned}
\]
Here local Poincar\'e inequalities also give
$\|\nabla w_j\|_{H^{r-\varepsilon}}\le C_R\|\nabla u\|_{H^{r-\varepsilon}(U_j)}$.
We choose $\eta$ so that $C\eta\le1/2$ and absorb the first term.
Then summing all the terms, and using bounded overlap yield
\[
 \sum_j\|w_j\|_{H^{r+1}}^2\le C_R\sum_j\bigl(\|\widetilde W\|_{H^r(U_j)}^2+\|\nabla u\|_{H^{r-\varepsilon}(U_j)}^2\bigr)\le C_R\bigl(\|\widetilde W\|_{H^r}^2+\|\nabla u\|_{H^{r-\varepsilon}}^2\bigr).
\]
Since $\sum_j\chi_j^2=1$ and $w_j=\chi_j(u-c_j)$,
\[
 \nabla u=\sum_j\chi_j\nabla w_j -\sum_j\chi_j(u-c_j)\nabla\chi_j.
\]
Bounded overlap and local Poincar\'e inequalities therefore give
\[
 \|\nabla u\|_{H^r}
 \le C_R\Big(\Bigl(\sum_j\|w_j\|_{H^{r+1}}^2\Bigr)^{1/2} +\|\nabla u\|_{H^{r-\varepsilon}}\Big).
\]
Combining these estimates, we obtain
\[
 \|\nabla u\|_{H^r} \leq C_R\bigl( \|\widetilde W\|_{H^r} +\|\nabla u\|_{H^{r-\varepsilon}}\bigr).
\]
For every $\theta>0$, interpolation gives
\[
 \|\nabla u\|_{H^{r-\varepsilon}} \leq \theta\|\nabla u\|_{H^r} +C_\theta\|\nabla u\|_{L^2}.
\]
Choosing $\theta$ sufficiently small and using \eqref{eq:elliptic-energy}, we obtain
\[
 \|\nabla u\|_{H^r}\leq C_R\bigl(\|\widetilde W\|_{H^r}+\|\nabla u\|_{L^2}\bigr)
 \leq C_R\|\widetilde W\|_{H^r}.
\]
Adding the boundary extension gives
\[
 \|\nabla v\|_{H^r} \leq C_R\bigl(\|W\|_{H^r}+\|h\|_{H^{r+1/2}}\bigr).
\]
By using smooth approximation, energy convergence, and weak lower semicontinuity, the above estimate  can be extended from smooth solutions to the stated variational solutions.
\end{proof}

Subtracting the fixed-domain equations \eqref{eq:pressure-weak} introduces the source $-(A_{f_1}-A_{f_2})\nabla v_2$ in divergence form.
We then have the following elliptic comparison estimate.
\begin{lemma}[Elliptic comparison]\label{lem:elliptic-comparison}
Suppose $f_1,f_2\in H^{s+1}(\cO)$ satisfy $\|f_j\|_{H^{s+1}}\le R$ (and $\inf_x(f_j-b)\ge h_0>0$ in finite depth), and let $v_j$ solve \eqref{eq:fixed-pressure} with top Dirichlet data $h_j$ and bulk data $\widetilde g_j$.
For $r=s,s-1$, writing $\delta f=f_1-f_2$ and similarly for the other variables, one has
\begin{equation}\label{eq:elliptic-comparison}
 \|\nabla\delta v\|_{H^r}
 \leq C_R\left(\|\delta h\|_{H^{r+1/2}}+\|\delta\widetilde g\|_{H^r} +\|\delta f\|_{H^{r+1/2}} (\|h_2\|_{H^{r+1/2}}+\|\widetilde g_2\|_{H^r})\right).
\end{equation}
\end{lemma}

\begin{proof}
Fix $r\in\{s-1,s\}$ and set
\[
 W_\delta =-\delta\widetilde g\,e_z
  -(A_{f_1}-A_{f_2})\nabla v_2.
\]
Subtracting \eqref{eq:pressure-weak} gives, in the variational sense,
\[
 \operatorname{div}(A_{f_1}\nabla\delta v) =\operatorname{div}W_\delta,
 \qquad \delta v|_{z=0}=\delta h,
\]
with $(A_{f_1}\nabla\delta v-W_\delta)\cdot e_z=0$ at the bottom in finite depth.
Lemma~\ref{lem:strip-elliptic} then yields
\[
 \|\nabla\delta v\|_{H^r} \leq C_R\bigl( \|\delta h\|_{H^{r+1/2}}+\|W_\delta\|_{H^r}\bigr)\leq C_R\bigl(  \|\delta h\|_{H^{r+1/2}}  +\|\delta\widetilde g\|_{H^r}  +\|A_{f_1}-A_{f_2}\|_{H^r} \|\nabla v_2\|_{H^r}\bigr).
\]
By \eqref{eq:geometry-difference} and
Lemma~\ref{lem:strip-elliptic},
\[
 \|A_{f_1}-A_{f_2}\|_{H^r} \leq C_R\|\delta f\|_{H^{r+1/2}},\quad
 \|\nabla v_2\|_{H^r} \leq C_R\bigl( \|h_2\|_{H^{r+1/2}}+\|\widetilde g_2\|_{H^r}\bigr).
\]
Using these bounds together, we prove \eqref{eq:elliptic-comparison}.
\end{proof}

\subsection{Curvature, velocity, and boundary forcing}
In this subsection, we apply the elliptic estimates to the pressure with boundary value $\sigma H(f)+\Gamma(f)$ and to the density-dependent boundary forcing.
Let $r=s$ or $s-1$, and suppose that the interfaces $f_j$ are bounded in $H^{s+1}$.

Write $H(f)=-\operatorname{div}A(\nabla f)$, where $A(p)=p/\sqrt{1+|p|^2}$.
The standard Sobolev composition and difference estimates \cite[Theorems~A.7--A.8]{Nguyen2020}, based on \cite[Theorem~2.89 and Corollary~2.90]{BCD2011}, applied to $A$ and $A-\operatorname{Id}$ (whose derivative vanishes at the origin), give
\begin{align}
 \|H(f)\|_{H^{r+1/2}}&\leq C_R\|f\|_{H^{r+5/2}},\label{eq:curv-single}\\
 \|H(f_1)-H(f_2)\|_{H^{r+1/2}}&\leq C_R\bigl(\|\delta f\|_{H^{r+5/2}} +\|\delta f\|_{H^{r+1}}\sum_{j=1}^2\|f_j\|_{H^{r+5/2}}\bigr).
 \label{eq:curv-diff}
\end{align}

Since $\Gamma(0)=0$, the Sobolev composition estimates give,
for $r=s,s-1$,
\begin{align*}
 \|\Gamma(f)\|_{H^{r+1/2}} \leq C_R\|f\|_{H^{r+1/2}},\quad \|\Gamma(f_1)-\Gamma(f_2)\|_{H^{r+1/2}}
 \leq C_R\|f_1-f_2\|_{H^{r+1/2}}.
\end{align*}

We write the density-dependent boundary term as
\[
 \mathcal B(f,\widetilde g)=
 (N_f\cdot\bS[f](\widetilde g\circ\Phi_f^{-1})+\widetilde g\circ\Phi_f^{-1})|_{\Sigma_f}.
\]
Set $\delta\mathcal B=\mathcal B(f_1,\widetilde g_1)-\mathcal B(f_2,\widetilde g_2)$.
For the single-solution and difference estimates, we apply Lemma~\ref{lem:strip-elliptic} and Lemma~\ref{lem:elliptic-comparison} respectively, with zero top data.
The trace theorem and the geometric coefficient estimates then give
\begin{align}
 &\|\mathcal B(f,\widetilde g)\|_{H^{s-1/2}}\leq C_R\|\widetilde g\|_{H^s},\label{eq:B-single}\\
 &\|\delta\mathcal B\|_{H^{r-1/2}}\leq C_R\bigl(\|\delta\widetilde g\|_{H^r} +\|\delta f\|_{H^{r+1}}\bigr),\quad r=s,s-1.\label{eq:B-diff}
\end{align}
The constant $C_R$ in \eqref{eq:B-diff} also depends on $\max_{j=1,2}\|\widetilde g_j\|_{H^s}$.
Since $\nabla f\in H^s$, multiplication by $N_f=(-\nabla f,1)$ preserves these trace estimates.
The $H^{s+1}$ bound on $f$ controls the coefficients of the pullback to $\mathcal S$.
Thus \eqref{eq:velocity-fixed}, \eqref{eq:strip-elliptic}, and \eqref{eq:curv-single} give
\begin{equation}\label{eq:velocity-single}
 \|\widetilde u\|_{H^s}\leq C_R\bigl(\|f\|_{H^{s+5/2}}+\|\widetilde g\|_{H^s}\bigr).
\end{equation}

\section{Estimates for fixed positive surface tension}\label{sec:estimates}
For fixed $\sigma>0$, energy estimates for the forced interface equation and the elliptic bounds control the transport field $V_f$ and source $Q_f$ in $L_T^2H^s$.
We shall use these bounds to obtain the transport and comparison estimates needed for the construction of solutions in Section~\ref{sec:construction}.

\subsection{The forced scalar interface equation}
Consider the forced interface equation:
\begin{equation}\label{eq:scalar-forced}
 f_t+\mathcal C(f)=F,\qquad
 \mathcal C(f):=\sigma\DN[f]H(f)+\DN[f]\Gamma(f).
\end{equation}

In the coupled system, $F=-\mathcal B(f,\widetilde g)$.
 \eqref{eq:B-single} gives
\[
 \|F\|_{L_T^2H^{s-1/2}} \le C_R\sqrt T\,\|\widetilde g\|_{C_TH^s}.
\]
At the $H^{s+1}$ energy level, the forcing contributes
\[
\begin{aligned}
 &2\bigl|\operatorname{Re}\langle F,f\rangle_{H^{s+1}}\bigr| =2\bigl|\operatorname{Re}\langle \langle D\rangle^{s-1/2}F, \langle D\rangle^{s+5/2}f\rangle_{L^2}\bigr|\\
 \le&2\|F\|_{H^{s-1/2}}\|f\|_{H^{s+5/2}}\le\varepsilon\|f\|_{H^{s+5/2}}^2 +C_\varepsilon\|F\|_{H^{s-1/2}}^2.
\end{aligned}
\]
For sufficiently small $\varepsilon>0$, the first term is absorbed into the capillary dissipation, while the second is integrable in time.
The following proposition adapts the capillary analysis of \cite[Theorem~1.2]{Nguyen2020} to this forcing and gives the energy, stability, and existence results used in the iteration later.

\begin{proposition}
\label{prop:scalar}
Fix $\sigma>0$.
All depth conditions below apply only in finite depth. Every smooth solution $f$ of \eqref{eq:scalar-forced} with $\|f\|_{C_TH^{s+1}}\le R$ and $f-b\ge h_0>0$ satisfies
\begin{equation}\label{eq:scalar-energy}
 \frac{\dd}{\dd t}\|f\|_{H^{s+1}}^2 +c_R\|f\|_{H^{s+5/2}}^2 \le C_R\|f\|_{H^{s+1}}^2 +C_R\|F\|_{H^{s-1/2}}^2,
\end{equation}
where $c_R,C_R>0$ may depend on the fixed $\sigma$.

For two solutions $f_1,f_2\in Z_T^{s+1}$ with forcings
$F_1,F_2$, satisfying the same interface and depth bounds
and $\|f_1\|_{Z_T^{s+1}}+\|f_2\|_{Z_T^{s+1}}\le M$,
one has, for $T\le1$ and $p=s,s+1$,
\begin{equation}\label{eq:scalar-stability}
 \|f_1-f_2\|_{Z_T^p} \le C_{R,M}\left( \|f_1(0)-f_2(0)\|_{H^p} +\|F_1-F_2\|_{L_T^2H^{p-3/2}} \right).
\end{equation}

For every $f_0\in H^{s+1}$ with $f_0-b\ge2h_0$ and every $F\in L_T^2H^{s-1/2}$, there exists $T_0\in(0,\min\{T,1\}]$ such that \eqref{eq:scalar-forced} has a unique solution
\[
 f\in Z_{T_0}^{s+1},\qquad
 f_t\in L_{T_0}^2H^{s-1/2},\qquad f(0)=f_0,
\]
with $f-b\ge h_0$ on $[0,T_0]$.
The lifespan and bounds are uniform for bounded initial data with a common positive depth bound, provided the forcing norms are uniformly small on a common short interval.
\end{proposition}

\begin{proof}
\emph{1. Operator estimates.}
Fix $p\in\{s,s+1\}$ and
$0<\alpha<\min\{1/2,s-1-d/2\}$.
The principal paradifferential symbols of the operators $\DN[f]$ and $H(f)$ are $\lambda_f$ and $\ell_f$, respectively:
\begin{equation}\label{eq:surface-symbols}
 \lambda_f =\sqrt{(1+|\nabla f|^2)|\xi|^2-(\nabla f\cdot\xi)^2},\qquad
 \ell_f =\frac{\lambda_f^2}{(1+|\nabla f|^2)^{3/2}},
 \qquad \lambda_f\ell_f\geq c_R|\xi|^3.
\end{equation}
The product $\lambda_f\ell_f$ is the principal symbol of $\DN[f]H(f)$.
For each $p$, the capillary estimates below require only a common $H^p$ bound on the interfaces; their constants depend on that bound, $p$, $\alpha$, and the fixed bottom and depth bounds.
The single-interface estimate is \cite[(3.10)]{Nguyen2020}:
\begin{equation}\label{eq:cap-para-single}
\|\DN[f]H(f)-T_{\lambda_f\ell_f}f\|_{H^{p-3/2+\alpha}} \leq C_R\|f\|_{H^{p+3/2}}.
\end{equation}
For the difference estimate, set
\[
 \delta f=f_1-f_2,\qquad K_p(t)=\sum_{j=1}^2\|f_j(t)\|_{H^{p+3/2}},
\]
The difference calculation in \cite[(3.24)--(3.33)]{Nguyen2020}, with zero gravity, gives
\begin{equation}\label{eq:cap-para-difference}
 \begin{aligned}
 \|\DN[f_1]H(f_1)-\DN[f_2]H(f_2)-T_{\lambda_{f_1}\ell_{f_1}}\delta f \|_{H^{p-3/2}}\leq C_R\bigl( \|\delta f\|_{H^{p+3/2-\alpha}}+K_p(t)\|\delta f\|_{H^p}\bigr)
 \end{aligned}
\end{equation}
for both $p=s$ and $p=s+1$.
The fixed-bottom and periodic settings are covered by Theorems~2.2--2.4 and Remark~1.4 in \cite{Nguyen2020}.

Note that $\Gamma(0)=0$, the Dirichlet--Neumann mapping and difference estimates give
\[
 \|\DN[f]\Gamma(f)\|_{H^{p-3/2}}\leq C_R\|f\|_{H^p},\qquad
 \|\DN[f_1]\Gamma(f_1)-\DN[f_2]\Gamma(f_2)\|_{H^{p-3/2}}\leq C_R\|\delta f\|_{H^p}.
\]
Consequently,
\begin{equation}\label{eq:scalar-operator-bound}
 \|\mathcal C(f)\|_{H^{p-3/2}}
 \le C_R\|f\|_{H^{p+3/2}}.
\end{equation}

\emph{2. Energy and stability.}
The adjoint and composition formulas show that $T_{\lambda_f\ell_f}-(T_{\sqrt{\lambda_f\ell_f}})^*T_{\sqrt{\lambda_f\ell_f}}$ has order at most $3-\alpha$.
Using the ellipticity of $\sqrt{\lambda_f\ell_f}$ and estimating the commutator with $\langle D\rangle^p$, we obtain,
\[
 \operatorname{Re}\langle T_{\lambda_f\ell_f}w,w\rangle_{H^p} \geq c_R\|w\|_{H^{p+3/2}}^2 -C_R\|w\|_{H^{p+3/2}} \|w\|_{H^{p+3/2-\alpha}}.
\]
For every $\varepsilon>0$, Sobolev interpolation and Young's inequality give
\[
 \|w\|_{H^{p+3/2}}\|w\|_{H^{p+3/2-\alpha}} \leq \varepsilon\|w\|_{H^{p+3/2}}^2 +C_\varepsilon\|w\|_{H^p}^2.
\]

We test \eqref{eq:scalar-forced} against $f$ at the $H^p$ level and estimate the remainder and forcing by Sobolev duality.
The remainder in \eqref{eq:cap-para-single} contributes at most $C_R\|f\|_{H^{p+3/2}} \|f\|_{H^{p+3/2-\alpha}}$.
The bound for $\DN[f]\Gamma(f)$ and the forcing contribute at most
\[
 \bigl(C_R\|f\|_{H^p}
       +\|F\|_{H^{p-3/2}}\bigr)\|f\|_{H^{p+3/2}}.
\]
Choosing $\varepsilon$ sufficiently small, we apply Young's inequality to absorb these terms into the capillary dissipation and obtain
\[
 \frac{\dd}{\dd t}\|f\|_{H^p}^2
 +c_R\|f\|_{H^{p+3/2}}^2
 \leq C_R\|f\|_{H^p}^2
      +C_R\|F\|_{H^{p-3/2}}^2.
\]
Taking $p=s+1$, we get the inequality \eqref{eq:scalar-energy}.

For the difference equation, \eqref{eq:cap-para-difference}
and the same argument give
\[
 \begin{aligned}
 \frac{\dd}{\dd t}\|\delta f\|_{H^p}^2
+c_R\|\delta f\|_{H^{p+3/2}}^2
 \leq C_R(1+K_p(t)^2)\|\delta f\|_{H^p}^2
+C_R\|F_1-F_2\|_{H^{p-3/2}}^2.
 \end{aligned}
\]
Since $\|K_p\|_{L_T^2}\leq M$, Gronwall's inequality yields
\[
 \|\delta f\|_{Z_T^p}^2 \leq C_R \exp\bigl(C_R(T+\|K_p\|_{L_T^2}^2)\bigr)\bigl( \|\delta f(0)\|_{H^p}^2 +\|F_1-F_2\|_{L_T^2H^{p-3/2}}^2\bigr).
\]
This proves \eqref{eq:scalar-stability}.
Spatial regularization extends these estimates to the stated energy class.

\emph{3. Construction and passage to the limit.}
Choose $R>2\|f_0\|_{H^{s+1}}$ and,  let the initial depth be at least $2h_0>0$ in the case of finite depth.
Let $P_N$ be the orthogonal Fourier projection to
$|\xi|\le 2^N$. Then
\begin{equation}\label{eq:fourier-tail}
 \|(I-P_N)h\|_{H^{s-3/2}}
 \le C2^{-N}\|h\|_{H^{s-1/2}}.
\end{equation}
Since
\[
 \|(I-P_N)f_0\|_{L^\infty} \leq C_\beta2^{-\beta N}\|f_0\|_{H^{s+1}}, \qquad 0<\beta<s+1-d/2,
\]
the projected data have the required depth for all sufficiently large $N$.
On the range of $P_N$, we solve the frequency-truncated equation
\[
 \partial_tf_N+P_N\mathcal C(f_N)=P_NF, \qquad f_N(0)=P_Nf_0.
\]
For each fixed $N$, the truncated operator $f_N\mapsto P_N\mathcal C(f_N)$ is locally Lipschitz on the frequency-localized space $P_NH^{s+1}(\cO)$ within the open subset of interfaces satisfying the positive-depth condition.
Since $P_NF\in L_T^2H^{s+1}(\cO)\subset L_T^1H^{s+1}(\cO)$, the Cauchy--Lipschitz theorem yields a unique local solution $f_N\in C([0,T_N];P_NH^{s+1}(\cO))$.

On an interval where $\|f_N\|_{H^{s+1}}\leq R$ and the depth is at least $h_0$, the identities $P_Nf_N=f_N$ and $P_N^*=P_N$ give
\[
 \begin{aligned}
\|f_N(t)\|_{H^{s+1}}^2 +c_R\int_0^t\|f_N(\tau)\|_{H^{s+5/2}}^2\,\dd\tau\leq\|f_0\|_{H^{s+1}}^2+C_RtR^2 +C_R\|F\|_{L_t^2H^{s-1/2}}^2.
 \end{aligned}
\]
The equation also yields
\[
 \|f_{N,t}\|_{L_t^2H^{s-1/2}} \leq C_R\|f_N\|_{L_t^2H^{s+5/2}} +\|F\|_{L_t^2H^{s-1/2}}.
\]
Since $H^{s-1/2}(\cO)\hookrightarrow L^\infty(\cO)$,
\[
 \sup_{\tau\leq t}\|f_N(\tau)-f_0\|_{L^\infty} \leq \|(P_N-I)f_0\|_{L^\infty} +Ct^{1/2}\|f_{N,t}\|_{L_t^2H^{s-1/2}}.
\]
Taking $t$ and the forcing norm sufficiently small, both bootstrap bounds are preserved on a common interval $[0,T_0]$, and we have the uniform estimate
\begin{equation}\label{eq:scalar-uniform}
 \sup_N (\|f_N\|_{Z_{T_0}^{s+1}} +\|f_{N,t}\|_{L_{T_0}^2H^{s-1/2}})\leq M_0.
\end{equation}
In infinite depth, the same argument applies without the bootstrap assumption on the depth.

We rewrite the approximate equation as
\[
 \partial_tf_N+\mathcal C(f_N)=F+E_N,
 \qquad E_N=(I-P_N)\bigl(\mathcal C(f_N)-F\bigr).
\]
By \eqref{eq:scalar-operator-bound} with $p=s+1$, \eqref{eq:fourier-tail}, and \eqref{eq:scalar-uniform},
we obtain
\[
 \begin{aligned}
\|E_N\|_{L_{T_0}^2H^{s-3/2}} \leq C2^{-N}\bigl( \|\mathcal C(f_N)\|_{L_{T_0}^2H^{s-1/2}} +\|F\|_{L_{T_0}^2H^{s-1/2}}\bigr)\leq C2^{-N}.
 \end{aligned}
\]
Applying \eqref{eq:scalar-stability} at $p=s$ gives
\[
 \begin{aligned}
 \|f_N-f_L\|_{Z_{T_0}^s} \leq C\bigl(\|(P_N-P_L)f_0\|_{H^s}+2^{-N}+2^{-L}\bigr)\leq C2^{-\min\{N,L\}}.
 \end{aligned}
\]
Thus $f_N\to f$ strongly in $Z_{T_0}^s$.
By weak-$*$ compactness, the limit also satisfies
\[
 f\in L_{T_0}^\infty H^{s+1} \cap L_{T_0}^2H^{s+5/2}, \qquad f_t\in L_{T_0}^2H^{s-1/2}.
\]
The operator difference bounds at $p=s$ give
\[
 \|\mathcal C(f_N)-\mathcal C(f)\|_{L_{T_0}^2H^{s-3/2}} \leq C_{R,M_0}\|f_N-f\|_{Z_{T_0}^s}\longrightarrow0.
\]
Thus $f$ solves \eqref{eq:scalar-forced}.

Since $f\in L^2([0,T_0];H^{s+5/2})$ and $f_t\in L^2([0,T_0];H^{s-1/2})$, the Lions--Magenes time-continuity lemma gives $f\in C([0,T_0];H^{s+1})$; see \cite[Theorem~6.41]{HunterPDE}.
The initial condition follows from the convergence of the approximations.
Uniqueness of \eqref{eq:scalar-forced} follows from the difference estimate \eqref{eq:scalar-stability}.
The existence time and estimates are uniform under the stated bounds on the initial data and depth and the smallness condition on the forcing.
\end{proof}

The scalar stability estimate \eqref{eq:scalar-stability} and the boundary forcing comparison \eqref{eq:B-diff} give the following estimate for the coupled interface equation.

\begin{corollary}[Coupled interface comparison]\label{cor:surface-comparison}
For two solutions with bounded $Z_T^{s+1}\times C_TH^s$ norms and common positive depth, we define
\[
 A_r=\|f_1-f_2\|_{Z_T^{r+1}},\qquad
 G_r=\|\widetilde g_1-\widetilde g_2\|_{C_TH^r}, \qquad r=s,s-1
\]
After decreasing $T$ in terms of those bounds,
\begin{equation}\label{eq:surface-comparison}
 A_r\leq C_R\bigl(\|f_1(0)-f_2(0)\|_{H^{r+1}}+\sqrt T\,G_r\bigr).
\end{equation}
\end{corollary}

\begin{proof}
We apply scalar stability with $F_j=-\mathcal B(f_j,\widetilde g_j)$ in \eqref{eq:scalar-stability} at $p=r+1$.
The trace comparison \eqref{eq:B-diff} bounds the forcing difference in $L_T^2H^{r-1/2}$ by $C_R\sqrt T(G_r+A_r)$.
For small $T$, we absorb the term containing $A_r$ on the left to obtain \eqref{eq:surface-comparison}.
\end{proof}

\subsection{Transport field and time-integrated comparisons}

We estimate the velocity and source in the fixed-domain transport equation \eqref{eq:ALE}.
At $z=0$, the kinematic condition $f_t=\widetilde u|_{z=0}\cdot N_f$ gives $V_z=0$.
In finite depth, the normal component also vanishes at the bottom by impermeability and $\varrho_{f,t}=0$.
Thus $V_f$ is tangent to the fixed boundaries.

On a bounded $Z_T^{s+1}\times C_TH^s$ set,
\eqref{eq:velocity-single} and the kinematic identity
$f_t=N_f\cdot\widetilde u|_{z=0}$ give
\begin{equation*}
 \|\widetilde u\|_{L_T^2H^s}+\|V_f\|_{L_T^2H^s}
 +\|Q_f\|_{L_T^2H^s}+\|f_t\|_{L_T^2H^{s-1/2}}
 \leq C_R\bigl(\|f\|_{L_T^2H^{s+5/2}}+\sqrt T\,\|\widetilde g\|_{C_TH^s}\bigr).
\end{equation*}
Indeed, the trace theorem, multiplication by $N_f=(-\nabla f,1)$, and \eqref{eq:time-extension} give
\[
\begin{aligned}
 \|f_t\|_{H^{s-1/2}}
 &=\|\widetilde u|_{z=0}\cdot N_f\|_{H^{s-1/2}}
 \le C_R\|\widetilde u|_{z=0}\|_{H^{s-1/2}}
 \le C_R\|\widetilde u\|_{H^s},\\
 \|\partial_t\varrho_f\|_{H^s}
 &\le C_R\|f_t\|_{H^{s-1/2}}
 \le C_R\|\widetilde u\|_{H^s}.
\end{aligned}
\]
The geometric coefficient bounds \eqref{eq:geometry-single} and \eqref{eq:ALE} then imply
\[
 \|V_f\|_{H^s}\le C_R\bigl(\|\widetilde u\|_{H^s}+\|\partial_t\varrho_f\|_{H^s}\bigr) \le C_R\|\widetilde u\|_{H^s}.
\]
For $Q_f$, we take $a=\ell_0$ in finite depth and $a=1$ in infinite depth.
The flattening estimates give
$\|\varrho_f-az\|_{H^s\cap L^\infty}\le C_R$.
Lemma~\ref{app:ln-multiplier}, which recalls \cite[Proposition~E.1]{LatoccaNguyen2025}, applies with $F=\gamma'$ and yields
\[
 \|Q_f\|_{H^s}
 =\|\gamma'(\varrho_f)\widetilde u_{d+1}\|_{H^s}
 \le C_R\|\widetilde u_{d+1}\|_{H^s\cap L^\infty}
 \le C_R\|\widetilde u\|_{H^s}.
\]
Integrating these inequalities in time and using \eqref{eq:velocity-single} gives the displayed estimate.

For two solutions, Lemma~\ref{lem:elliptic-comparison} and \eqref{eq:curv-diff} imply, for $r=s,s-1$,
\begin{equation*}
 \|\delta\widetilde u(t)\|_{H^r} \leq C_R\left(\|\delta f(t)\|_{H^{r+5/2}} +\|\delta\widetilde g(t)\|_{H^r} +(1+K_{r+1}(t))\|\delta f(t)\|_{H^{r+1}}\right).
\end{equation*}
The difference of the kinematic traces gives the corresponding bound for $\delta f_t$ in $H^{r-1/2}$.
Applying the time-extension bound \eqref{eq:time-extension} and subtracting the field formulas \eqref{eq:ALE}, we obtain
\begin{equation}\label{eq:ALE-difference}
 \|\delta\widetilde u\|_{L_T^2H^r}+\|\delta V\|_{L_T^2H^r}+\|\delta Q\|_{L_T^2H^r}+\|\delta f_t\|_{L_T^2H^{r-1/2}} \leq C_R(A_r+\sqrt T\,G_r).
\end{equation}
The constant in \eqref{eq:ALE-difference} depends only on the common solution bounds and fixed parameters, since
\[
 \|K_{r+1}\|_{L^2(0,T)} \leq \sum_{j=1}^2\|f_j\|_{L_T^2H^{r+5/2}} \leq C_R.
\]
By the Cauchy--Schwarz inequality in time,
\begin{equation}\label{eq:ALE-L1-difference}
 \begin{aligned}
 \|\delta V\|_{L_T^1H^r}+\|\delta Q\|_{L_T^1H^r}
 \leq \sqrt T\bigl( \|\delta V\|_{L_T^2H^r} +\|\delta Q\|_{L_T^2H^r}\bigr)\leq C_R(\sqrt T\,A_r+T G_r).
 \end{aligned}
\end{equation}
\subsection{Transport estimates and lower-order uniqueness}

We apply the transport estimate for tangent vector fields in Lemma~\ref{app:ln-transport} at index $s$.
The hypotheses follow from $s>D/2+1$ and $V\in L_T^2H^s$.
In particular, $V\in L_T^2L^2$ and $\nabla V\in L_T^1H^{s-1}$.
For the transport equation
\[
(\partial_t+V\cdot\nabla) h=Q,
\]
the lemma gives a unique solution $h\in C_TH^s$ satisfying
\begin{equation}\label{eq:transport-high}
 \|h\|_{C_TH^s} \leq C_s\exp(C\|V\|_{L_T^1H^s}) \bigl(\|h(0)\|_{H^s}+C\|Q\|_{L_T^1H^s}\bigr).
\end{equation}
We also need a transport estimate at $s-1$ with coefficient $\|V\|_{L_T^1H^s}$.
For $r=s-1>D/2$, the Sobolev commutator bound gives
\[
 \|[\langle D\rangle^r,V]\nabla h\|_{L^2} \leq C\|V\|_{H^{r+1}}\|h\|_{H^r}.
\]

On $\cO\times\mathbb R$, the energy estimate and the commutator bound give the transport estimate, including $r=1+D/2$;
see \cite{BCD2011} for the commutator argument.
To obtain the strip estimate, we extend $V$, $Q$, and $h(0)$ to $\cO\times\mathbb R$ and solve the extended equation as in Lemma~\ref{lem:reference-smoothing}.
Since the extended velocity is Lipschitz and tangent to $\partial\mathcal S$, its flow preserves $\mathcal S$.
By uniqueness, the extended solution restricts to $h$, which therefore satisfies the same estimate.

For two solutions, the difference equation is
\begin{equation*}
 \partial_t\delta\widetilde g+V_1\cdot\nabla\delta\widetilde g
 =\delta Q-\delta V\cdot\nabla\widetilde g_2.
\end{equation*}
Since $H^{s-1}(\mathcal S)$ is an algebra,
$\|\delta V\cdot\nabla\widetilde g_2\|_{H^{s-1}}
\leq C_R\|\delta V\|_{H^{s-1}}$.
Applying the lower-order transport estimate at index $s-1$ together with \eqref{eq:ALE-L1-difference} yields
\[
 G_{s-1}\leq C_R\bigl(\|\delta\widetilde g(0)\|_{H^{s-1}}+\sqrt T\,A_{s-1}+T G_{s-1}\bigr).
\]
Substituting the bound for $A_{s-1}$ from \eqref{eq:surface-comparison} leaves only the initial differences and $G_{s-1}$ on the right.
For a sufficiently short common interval, we absorb the remaining unknown term and obtain
\begin{equation}\label{eq:lower-stability}
 \|\delta f\|_{Z_T^s}+\|\delta\widetilde g\|_{C_TH^{s-1}} \leq C_R\bigl(\|\delta f(0)\|_{H^s} +\|\delta\widetilde g(0)\|_{H^{s-1}}\bigr).
\end{equation}
Once existence is established, this estimate gives uniqueness in the class \eqref{eq:main-regularity} and Lipschitz dependence in the stated lower-order topology.

\section{The proof of local well-posedness}\label{sec:construction}\label{sec:continuity}

In this section, we construct solutions by an iteration whose transport fields remain tangent to the fixed boundary, then prove continuous dependence, positivity of the density, and a continuation criterion.

\subsection{An iteration preserving boundary tangency}

Following \cite[(6.10)--(6.15)]{LatoccaNguyen2025}, we construct each iterate by solving the scalar surface equation of Proposition~\ref{prop:scalar} and a linear transport equation for the density.

We begin with the unforced surface equation
\[
 \partial_t f_1+\mathcal C(f_1)=0,\qquad
 f_1(0)=f_0,\qquad \widetilde g_1=0.
\]
Let $\widetilde u_1$ be the corresponding harmonic Darcy velocity.
We define $V_1$ by \eqref{eq:ALE}, using $\partial_t f_1=-\mathcal C(f_1)$, and set
\[
 Q_1=-\gamma'(\varrho_{f_1})\widetilde u_{1,d+1}.
\]
The identity $\partial_t f_1=-\mathcal C(f_1)$ gives
$V_{1,z}|_{z=0}=0$.
In finite depth, the bottom condition gives tangency there as well.

For $n\ge2$, we solve
\begin{equation}\label{eq:iterate-surface}
 \begin{cases}
 \partial_t f_n+\mathcal C(f_n)
 =-\mathcal B(f_{n-1},\widetilde g_{n-1}),
 \qquad f_n(0)=f_0,\\
 \partial_t\widetilde g_n
 +V_{n-1}\cdot\nabla\widetilde g_n=Q_{n-1},
 \qquad \widetilde g_n(0)=\widetilde g_0.
 \end{cases}
\end{equation}
Set $g_n=\widetilde g_n\circ\Phi_{f_n}^{-1}$ and define the associated harmonic and density contributions by
\[
 C_n=\bigl(\bG[f_n](\sigma H(f_n)+\Gamma(f_n))\bigr) \circ\Phi_{f_n},
 \qquad B_n=(\bS[f_n]g_n)\circ\Phi_{f_n}
       +\widetilde g_n e_{d+1}.
\]
The pulled-back physical velocity and density source are
\begin{equation}\label{eq:iterate-physical}
 \widetilde u_n=-C_n-B_n, \qquad
 Q_n=-\gamma'(\varrho_{f_n})\widetilde u_{n,d+1}.
\end{equation}
The first equation in \eqref{eq:iterate-surface}
uses the density forcing from the preceding iterate.
To preserve boundary tangency, we use the same density flux in the vertical component of the transport field:
\begin{equation}\label{eq:iterate-field}
 V_{n,x}=\widetilde u_{n,x},\qquad 
 V_{n,z}=J_{f_n}^{-1} \bigl(-n_{f_n}\cdot C_n -n_{f_{n-1}}\cdot B_{n-1} -\partial_t\varrho_{f_n}\bigr),
\end{equation}
where $n_f=(-\nabla_x\varrho_f,1)$, and then
$n_f|_{z=0}=N_f$.

At $z=0$, the first equation in \eqref{eq:iterate-surface} gives
\[
 \left.(J_{f_n}V_{n,z})\right|_{z=0} =-\mathcal C(f_n) -\mathcal B(f_{n-1},\widetilde g_{n-1}) -\partial_t f_n =0.
\]
In finite depth, the bottom pressure conditions gives
\[
 n_{f_n}\cdot C_n =n_{f_{n-1}}\cdot B_{n-1} =\partial_t\varrho_{f_n}=0
 \qquad\text{at }z=-1.
\]
Thus $V_{n,z}=0$ on each fixed boundary.

\subsection{Uniform bounds and a common lifespan}\label{sec:local-uniform-lifespan}

Fix constants
\[
 R_f>\|f_0\|_{H^{s+1}}+1,\qquad R_g>2C_s(\|\widetilde g_0\|_{H^s}+1),
\]
and choose $L$ satisfying
\[
 c_{R_f}L^2>\|f_0\|_{H^{s+1}}^2+1,
\]
where $c_{R_f}>0$ is the coercivity constant in \eqref{eq:scalar-energy} under the interface and depth bounds.
We prove by induction that
\begin{equation}\label{eq:uniform-iterates}
 \sup_n\|f_n\|_{C_TH^{s+1}}\leq R_f,\qquad
 \sup_n\|f_n\|_{L_T^2H^{s+5/2}}\leq L,\qquad
 \sup_n\|\widetilde g_n\|_{C_TH^s}\leq R_g,
\end{equation}
together with $f_n-b\geq h_0$ in finite depth and
\begin{equation}\label{eq:iterate-transport-bound}
 \|V_n\|_{L_T^2H^s}+\|Q_n\|_{L_T^2H^s}
 \leq C_R(L+\sqrt T\,R_g).
\end{equation}
For $n=1$, applying \eqref{eq:scalar-energy} with $F=0$
gives the first two bounds in \eqref{eq:uniform-iterates}
for sufficiently small $T$.
The density bound follows from $\widetilde g_1=0$.
By \eqref{eq:velocity-single}, the kinematic identity
$f_{1,t}=N_{f_1}\cdot\widetilde u_1|_{z=0}$, and
\eqref{eq:time-extension},
\[
 \|\widetilde u_1\|_{L_T^2H^s}
 +\|f_{1,t}\|_{L_T^2H^{s-1/2}}
 +\|\partial_t\varrho_{f_1}\|_{L_T^2H^s}
 \le C_R\|f_1\|_{L_T^2H^{s+5/2}}
 \le C_RL.
\]
The formulas in \eqref{eq:ALE} then give
\[
 \|V_1\|_{L_T^2H^s}+\|Q_1\|_{L_T^2H^s}
 \le C_RL,
\]
which proves \eqref{eq:iterate-transport-bound} for $n=1$.
In finite depth,
$\sup_{t\le T}\|f_1(t)-f_0\|_{L^\infty}
\le C_R\sqrt T\,L$
preserves the depth bound for sufficiently small $T$.

Suppose that all these bounds hold through step $n-1$.
By \eqref{eq:B-single}, the forcing in the first equation
of \eqref{eq:iterate-surface} satisfies
\[
 \|\mathcal B(f_{n-1},\widetilde g_{n-1})\|_{L_T^2H^{s-1/2}} \leq C_{R_f}\sqrt T\,R_g.
\]
On any interval where $\|f_n(t)\|_{H^{s+1}}\le R_f$
and, in finite depth, $f_n(t)-b\ge h_0$,
\eqref{eq:scalar-energy} and Gronwall's inequality give
\[
 \|f_n(t)\|_{H^{s+1}}^2 +c_{R_f}\int_0^t\|f_n(\tau)\|_{H^{s+5/2}}^2\,\dd\tau \leq e^{C_{R_f}t} \bigl(\|f_0\|_{H^{s+1}}^2+C_{R_f}tR_g^2\bigr).
\]
For $T$ sufficiently small, the right-hand side is less than $\|f_0\|_{H^{s+1}}^2+1$.
By the choices of $R_f$ and $L$, this gives the first two bounds in \eqref{eq:uniform-iterates}.
The first equation in \eqref{eq:iterate-surface} gives
\begin{equation}\label{eq:iterate-time-bound}
 \|f_{n,t}\|_{L_T^2H^{s-1/2}}
 \leq C_R(L+\sqrt T\,R_g).
\end{equation}
Since $H^{s-1/2}(\cO)\hookrightarrow L^\infty(\cO)$,
\[
 \sup_{t\leq T}\|f_n(t)-f_0\|_{L^\infty} \leq C_R\sqrt T(L+\sqrt T\,R_g).
\]
After decreasing $T$, this bound preserves the depth and permits the same flattening parameters at each step.

Applying \eqref{eq:transport-high} and the Cauchy--Schwarz inequality in time gives
\[
 \|\widetilde g_n\|_{C_TH^s} \leq C_s e^{C_R\sqrt T(L+\sqrt T\,R_g)} \bigl(\|\widetilde g_0\|_{H^s} +C_R\sqrt T(L+\sqrt T\,R_g)\bigr).
\]
For sufficiently small $T$, the right-hand side is less than $R_g$.
Applying Lemma~\ref{lem:strip-elliptic} with
$(h,W)=(\sigma H(f_n)+\Gamma(f_n),0)$ for $C_n$
and $(h,W)=(0,-\widetilde g_j e_z)$ for $B_j$,
and using \eqref{eq:curv-single} and Sobolev composition
for $\Gamma$, we obtain
\[
\begin{aligned}
 \|C_n\|_{L_T^2H^s}
 &\le C_R\|f_n\|_{L_T^2H^{s+5/2}}\le C_RL,\\
 \|B_j\|_{L_T^2H^s}
 &\le C_R\sqrt T\,\|\widetilde g_j\|_{C_TH^s}
 \le C_R\sqrt T\,R_g,\qquad j=n-1,n.
\end{aligned}
\]
By \eqref{eq:time-extension} and \eqref{eq:iterate-time-bound},
\[
 \|\partial_t\varrho_{f_n}\|_{L_T^2H^s}
 \le C_R\|f_{n,t}\|_{L_T^2H^{s-1/2}}
 \le C_R(L+\sqrt T\,R_g).
\]
Using these bounds in \eqref{eq:iterate-physical}--\eqref{eq:iterate-field}, we obtain
\[
\begin{aligned}
 &\|\widetilde u_n\|_{L_T^2H^s}+\|Q_n\|_{L_T^2H^s}
 \le C_R\bigl(\|C_n\|_{L_T^2H^s}+\|B_n\|_{L_T^2H^s}\bigr)
 \le C_R(L+\sqrt T\,R_g),\\
 &\|V_n\|_{L_T^2H^s}\le C_R\bigl(\|\widetilde u_n\|_{L_T^2H^s} +\|C_n\|_{L_T^2H^s}+\|B_{n-1}\|_{L_T^2H^s}+\|\partial_t\varrho_{f_n}\|_{L_T^2H^s}\bigr)\le C_R(L+\sqrt T\,R_g).
\end{aligned}
\]
This proves \eqref{eq:iterate-transport-bound} at index $n$.
This closes the induction on a common interval $[0,T]$.

\subsection{Convergence of the iteration}

For $n\ge2$, set
\[
 X_n=\|f_n-f_{n-1}\|_{Z_T^s},\qquad
 Y_n=\|\widetilde g_n-\widetilde g_{n-1}\|_{C_TH^{s-1}},
 \qquad D_n=X_n+Y_n.
\]
The uniform bounds \eqref{eq:uniform-iterates} imply that $D_2+D_3<\infty$.
Since successive surface iterates have the same initial data, \eqref{eq:scalar-stability} and \eqref{eq:B-diff} give
\[
 X_n\leq C_R\sqrt T\,(X_{n-1}+Y_{n-1}),
 \qquad n\ge3.
\]

We next compare the variables in \eqref{eq:iterate-physical}--\eqref{eq:iterate-field}.
To estimate $C_m-C_{m-1}$, we apply Lemma~\ref{lem:elliptic-comparison} at $r=s-1$ and bound the top data using \eqref{eq:curv-diff} and the composition estimates for $\Gamma$.
For $B_m-B_{m-1}$, the same lemma applies with zero top data.
The geometric coefficient bounds then give
 \begin{align*}
 \|C_m-C_{m-1}\|_{L_T^2H^{s-1}} &\leq C_R X_m,\\
 \|B_m-B_{m-1}\|_{L_T^2H^{s-1}} &\leq C_R\sqrt T\,(X_m+Y_m), \qquad m\ge2.
 \end{align*}
The same estimate controls $B_{m-1}-B_{m-2}$, the difference of density fluxes appearing in $V_{m,z}-V_{m-1,z}$.
Subtracting the first equations in \eqref{eq:iterate-surface} gives
\[
 \|\partial_t(f_m-f_{m-1})\|_{L_T^2H^{s-3/2}} \leq C_R\bigl( X_m+\sqrt T\,(X_{m-1}+Y_{m-1})\bigr), \qquad m\ge3.
\]
Since
\[
 \partial_t(\varrho_{f_m}-\varrho_{f_{m-1}}) =\mathscr E\,\partial_t(f_m-f_{m-1}),
\]
the time-extension estimate \eqref{eq:time-extension} controls the time derivative of the difference of the coordinate maps.
Using these bounds in
\eqref{eq:iterate-physical}--\eqref{eq:iterate-field},
together with \eqref{eq:geometry-single}
and \eqref{eq:geometry-difference} at $r=s-1$, we obtain 
\[
\|V_m-V_{m-1}\|_{L_T^2H^{s-1}} +\|Q_m-Q_{m-1}\|_{L_T^2H^{s-1}}\leq C_R\bigl( X_m+X_{m-1}+\sqrt T\,(Y_m+Y_{m-1})\bigr), \quad m\ge3.
\]

For $n\ge4$, subtracting the density equations gives
\[
 \begin{aligned}
\partial_t(\widetilde g_n-\widetilde g_{n-1})
 +V_{n-1}\cdot\nabla(\widetilde g_n-\widetilde g_{n-1})=Q_{n-1}-Q_{n-2}
 -(V_{n-1}-V_{n-2})\cdot\nabla\widetilde g_{n-1},
 \end{aligned}
\]
with zero initial datum.
The lower-order transport estimate, the algebra property of $H^{s-1}$, and the Cauchy--Schwarz inequality in time imply
\[
 Y_n\leq C_R\bigl( \sqrt T\,(X_{n-1}+X_{n-2}) +T(Y_{n-1}+Y_{n-2})\bigr).
\]
Combining this estimate with the bound for $X_n$ and choosing $T\le1$ sufficiently small gives
\[
 D_n\leq q_T(D_{n-1}+D_{n-2}), \qquad n\ge4,\qquad
 q_T=C_R\sqrt T<\frac12.
\]
We sum the recurrence for $n=4,\ldots,N$ and add $D_2+D_3$ to both sides.
Since $D_j\ge0$, both $\sum_{j=3}^{N-1}D_j$
and $\sum_{j=2}^{N-2}D_j$
are bounded by $\sum_{j=2}^{N}D_j$.
Hence
\[
 \sum_{n=2}^N D_n \leq D_2+D_3+2q_T\sum_{n=2}^N D_n.
\]
Absorbing the last term and letting $N\to\infty$, we obtain
\[
 \sum_{n=2}^{\infty}D_n \leq \frac{D_2+D_3}{1-2q_T}<\infty.
\]
Thus $(f_n,\widetilde g_n)$ converges strongly in $Z_T^s\times C_TH^{s-1}$.

\subsection{Identification of the limit and time continuity}

Let $(f,\widetilde g)$ be the strong limit.
The uniform bounds and weak lower semicontinuity yield
\[
 f\in L_T^\infty H^{s+1}\cap L_T^2H^{s+5/2},
 \qquad f_t\in L_T^2H^{s-1/2},
 \qquad \widetilde g\in L_T^\infty H^s.
\]

\emph{The surface equation.}
The capillary comparison \eqref{eq:cap-para-difference}, the order-three bound for the principal operator, and \eqref{eq:B-diff} imply
\[
 \begin{aligned}
 \mathcal C(f_n)&\longrightarrow\mathcal C(f),\\
 \mathcal B(f_{n-1},\widetilde g_{n-1}) &\longrightarrow\mathcal B(f,\widetilde g)
 \end{aligned}
 \qquad\text{in }L_T^2H^{s-3/2}.
\]
Passing to the limit in the first equation of \eqref{eq:iterate-surface}, we get the limiting equation
\[
 f_t+\mathcal C(f)=-\mathcal B(f,\widetilde g).
\]
The equations for the iterates and the convergence above also give $f_{n,t}\to f_t$ in $L_T^2H^{s-3/2}$.

\emph{The transport equation.}
By \eqref{eq:elliptic-comparison} and
\eqref{eq:time-extension}, both with $r=s-1$,
we can pass to the limit in \eqref{eq:iterate-field}:
\[
 (\widetilde u_n,V_n,Q_n) \longrightarrow(\widetilde u,V_f,Q_f) \quad\text{in }L_T^2H^{s-1}(\mathcal S).
\]
Here $V_f$ is the field defined by \eqref{eq:ALE}.
Weak compactness also gives $V_f,Q_f\in L_T^2H^s(\mathcal S)$, while the boundary tangency of $V_n$ is preserved in the limit.

Since $s-1>D/2$, the Sobolev product estimate gives
\[
 V_{n-1}\cdot\nabla\widetilde g_n \longrightarrow V_f\cdot\nabla\widetilde g \quad\text{in }L_T^2H^{s-2}.
\]
Thus the limiting density satisfies
\[
 \partial_t\widetilde g+V_f\cdot\nabla\widetilde g=Q_f
\]
in distributions.
The strong convergence at the lower indices also gives $f(0)=f_0$ and $\widetilde g(0)=\widetilde g_0$.

\emph{Time continuity and conclusion.}
The bounds on $f$ and $f_t$ imply
$f\in C_TH^{s+1}$ by the standard time-continuity argument.
The tangent transport theorem gives a solution in $C_TH^s$ with field $V_f$, source $Q_f$, and initial datum $\widetilde g_0$.
Uniqueness in $L^2$ identifies this solution with the distributional limit $\widetilde g$.

The physical pressure and velocity are recovered from \eqref{eq:fixed-pressure} and \eqref{eq:velocity-fixed}.
Changing variables by $\Phi_f$ gives a solution of \eqref{e:IPM} with regularity \eqref{eq:main-regularity}.
Uniqueness and lower-order Lipschitz dependence follow from the difference estimate \eqref{eq:lower-stability}.

\subsection{Smoothing the reference transport solution}

A direct $H^s$ estimate for the difference of two densities
would require an $H^{s+1}$ bound on one of them.
We therefore fix one density solution and approximate it
by smooth functions.
These approximations converge in $C_TH^s$,
and their errors in the transport equation tend to zero
in $L_T^1H^s$.

\begin{lemma}[Reference transport smoothing]
\label{lem:reference-smoothing}
Let $s>D/2+1$ and let $V\in L_T^2H^s(\mathcal S)$
satisfy $V_{d+1}=0$ on $\partial\mathcal S$, and  $Q\in L_T^1H^s(\mathcal S)$.
Suppose that $h\in C_TH^s(\mathcal S)$ solves
\[
 \partial_t h+V\cdot\nabla h=Q.
\]
There exist spatially smooth functions $h^K$ such that
\begin{equation}\label{eq:reference-smoothing}
 \begin{gathered}
 h^K\longrightarrow h\quad\text{in }C_TH^s,\\
 \partial_t h^K+V\cdot\nabla h^K=Q+r_K, \qquad
 \|r_K\|_{L_T^1H^s}\longrightarrow0.
 \end{gathered}
\end{equation}
Moreover,
\[
 \sup_K\|h^K\|_{C_TH^s}\leq C_R, \qquad
 \|h^K\|_{C_TH^{s+1}}\leq C_R2^K,
\]
where $R$ bounds
\[
 \|h(0)\|_{H^s} +\|Q\|_{L_T^1H^s} +T^{1/2}\|V\|_{L_T^2H^s}.
\]
The constants also depend on $s$ and the fixed strip.
\end{lemma}

\begin{proof}
Fix a bounded Sobolev extension operator $E$ from $\mathcal S$ to $\cO\times\mathbb R$, we solve the extended transport equation
\begin{equation}\label{eq:extended-transport}
 \partial_t h^e+(EV)\cdot\nabla h^e=EQ, \qquad h^e(0)=Eh(0).
\end{equation}
The standard $H^s$ transport estimate gives
\[
 \|h^e\|_{C_TH^s(\cO\times\mathbb R)} \leq C\exp\bigl(C\|V\|_{L_T^1H^s}\bigr) \bigl(\|h(0)\|_{H^s}+\|Q\|_{L_T^1H^s}\bigr) \leq C_R.
\]
The field $EV\in L_T^1W^{1,\infty}$ has zero normal trace on $\partial\mathcal S$, hence its flow preserves $\mathcal S$.
Uniqueness for the transport equation then gives $h^e|_{\mathcal S}=h$.

Let $S_K$ be smooth Fourier cutoffs on $\cO\times\mathbb R$ at frequencies $\lesssim2^K$.
We first prove a commutator estimate uniform in $K$:
\begin{equation}\label{eq:Friedrichs-Hs}
 \begin{aligned}
 \|[U\cdot\nabla,S_K]w\|_{H^s}
 &\leq C( \|\nabla U\|_{L^\infty}\|w\|_{H^s} +\|U\|_{H^s}\|\nabla w\|_{L^\infty})\\ &\leq C\|U\|_{H^s}\|w\|_{H^s}, \qquad U,w\in H^s(\cO\times\mathbb R).
 \end{aligned}
\end{equation}
Indeed,
\[
 \langle D\rangle^s[U\cdot\nabla,S_K]w ={}[\langle D\rangle^s,U\cdot\nabla]S_Kw -S_K[\langle D\rangle^s,U\cdot\nabla]w+[U\cdot\nabla,S_K]\langle D\rangle^sw.
\]
The Kato--Ponce commutator estimate bounds the first two terms in $L^2$ by the first right-hand side of \eqref{eq:Friedrichs-Hs}.
For the last term, integration by parts in the convolution formula gives
\[
 \|[U\cdot\nabla,S_K]z\|_{L^2}\leq C\|\nabla U\|_{L^\infty}\|z\|_{L^2}.
\]
Here the cutoff kernels satisfy, uniformly in $K$,
\[
 \|k_K\|_{L^1(\cO\times\mathbb R)} +\int_{\cO\times\mathbb R}  \operatorname{dist}_{\cO\times\mathbb R}(Y,0) |\nabla k_K(Y)|\,\dd Y \leq C.
\]
In the periodic case, these kernels are obtained by
periodizing in the horizontal variables.
Taking $z=\langle D\rangle^sw$ proves \eqref{eq:Friedrichs-Hs} for every real $s>D/2+1$.

For fixed $U\in H^s$ and $w\in H^{s+1}$,
\[
 [U\cdot\nabla,S_K]w =U\cdot\nabla(S_K-I)w-(S_K-I)(U\cdot\nabla w) \longrightarrow0 \quad\text{in }H^s.
\]
By density and the uniform bound \eqref{eq:Friedrichs-Hs}, the same convergence holds for $U,w\in H^s$.
At almost every time, we apply this result with $U=EV$ and $w=h^e$.
The commutator estimate gives
\[
 \|[(EV)\cdot\nabla,S_K]h^e\|_{H^s}
 \leq C\|EV\|_{H^s}\|h^e\|_{H^s},
\]
whose right-hand side belongs to $L_T^1$.
The dominated convergence theorem gives
\[
 [(EV)\cdot\nabla,S_K]h^e \longrightarrow0 \quad\text{in }L_T^1H^s(\cO\times\mathbb R).
\]

Since $S_K$ localizes to frequencies $\lesssim2^K$,
$h^K$ is smooth in space and
\[
 \|(S_Kh^e)|_{\mathcal S}\|_{C_TH^s} +2^{-K}\|(S_Kh^e)|_{\mathcal S}\|_{C_TH^{s+1}} \leq C\|h^e\|_{C_TH^s(\cO\times\mathbb R)} \leq C_R.
\]
The set $h^e([0,T])$ is compact in $H^s(\cO\times\mathbb R)$.
Since $S_K\to I$ strongly and the operators are uniformly bounded, $(S_Kh^e)|_{\mathcal S}\to h$ in $C_TH^s(\mathcal S)$.

Finally, applying $S_K$ to \eqref{eq:extended-transport}
and restricting to $\mathcal S$ yields
\[
 (\partial_t +V\cdot\nabla) (S_Kh^e)|_{\mathcal S}=Q+r_K,
 \qquad r_K=\bigl((S_K-I)(EQ) +[(EV)\cdot\nabla,S_K]h^e\bigr)|_{\mathcal S}.
\]
Both terms in $r_K$ tend to zero in $L_T^1H^s$, so that $\|r_K\|_{L_T^1H^s(\mathcal S)}\to0$.
Setting $h^K=(S_Kh^e)|_{\mathcal S}$, we obtain the estimate \eqref{eq:reference-smoothing}.
\end{proof}

\subsection{Continuous dependence}

Suppose
$(f_{n,0},\widetilde g_{n,0})
 \longrightarrow (f_0,\widetilde g_0)$
 in $H^{s+1}\times H^s$.
Let $(f_n,\widetilde g_n)$ and $(f,\widetilde g)$ be the corresponding solutions, with associated vector fields $(\widetilde u_n,V_n,Q_n)$ and $(\widetilde u,V,Q)$.
For all sufficiently large $n$, the uniform construction provides a common interval $[0,T]$, uniform bounds in $Z_T^{s+1}\times C_TH^s$, and a common positive depth bound in finite depth.
Constants $C_R$ below depend only on these bounds and the fixed parameters.

Set
\[
 A_n=\|f_n-f\|_{Z_T^{s+1}},\qquad
 G_n=\|\widetilde g_n-\widetilde g\|_{C_TH^s}.
\]
The lower-order stability estimate \eqref{eq:lower-stability} and the field comparison \eqref{eq:ALE-difference} give
\begin{equation}\label{eq:lower-field-convergence}
 \|V_n-V\|_{L_T^1H^{s-1}}\longrightarrow0.
\end{equation}
The surface comparison \eqref{eq:surface-comparison}
at $r=s$ gives
\begin{equation}\label{eq:top-surface-continuity}
 A_n\leq C_R\bigl( \|f_{n,0}-f_0\|_{H^{s+1}}+\sqrt T\,G_n\bigr).
\end{equation}

To estimate $G_n$, we apply Lemma~\ref{lem:reference-smoothing} to the fixed reference density $\widetilde g$.
Denote the approximation by $g^K$ and its residual by $r_K$. Then
\[
 \partial_t g^K+V\cdot\nabla g^K=Q+r_K, \qquad \varepsilon_K:= \|g^K-\widetilde g\|_{C_TH^s} +\|r_K\|_{L_T^1H^s}\longrightarrow0.
\]
For $w_{n,K}=\widetilde g_n-g^K$, subtraction gives
\[
 \partial_t w_{n,K}+V_n\cdot\nabla w_{n,K} =Q_n-Q-r_K-(V_n-V)\cdot\nabla g^K,
\]
with initial value
$w_{n,K}(0)=\widetilde g_{n,0}-g^K(0)$.
Since $s>D/2+1$, the Sobolev product estimate yields
\begin{equation}\label{eq:smooth-reference-product}
 \begin{aligned}
 \|(V_n-V)\cdot\nabla g^K\|_{H^s} &\leq C\bigl( \|V_n-V\|_{H^s}\|\nabla g^K\|_{L^\infty} +\|V_n-V\|_{L^\infty}\|g^K\|_{H^{s+1}}\bigr)\\
 &\leq C_R\|V_n-V\|_{H^s} +C_K\|V_n-V\|_{H^{s-1}}.
 \end{aligned}
\end{equation}
The uniform $C_TH^s$ bound for $g^K$ makes $C_R$ independent of $K$.

Applying \eqref{eq:transport-high} to $w_{n,K}$ and using
\eqref{eq:smooth-reference-product},
\eqref{eq:ALE-L1-difference}, and the definition of
$\varepsilon_K$, we obtain
\begin{equation}\label{eq:top-density-continuity}
 \begin{aligned}
 G_n\leq{}& C_R\bigl(\|\widetilde g_{n,0}-\widetilde g_0\|_{H^s} +\varepsilon_K+\sqrt T\,A_n+T G_n\bigr)+C_K\|V_n-V\|_{L_T^1H^{s-1}}.
 \end{aligned}
\end{equation}
We substitute \eqref{eq:top-surface-continuity} into \eqref{eq:top-density-continuity} and choose $T\leq1$ sufficiently small which is independent of $n$ and $K$ to absorb the terms containing $G_n$.
Using \eqref{eq:top-surface-continuity} again, we obtain
\[
 \begin{aligned}
 A_n+G_n\leq C_R\bigl( \|f_{n,0}-f_0\|_{H^{s+1}} +\|\widetilde g_{n,0}-\widetilde g_0\|_{H^s} +\varepsilon_K\bigr)+C_K\|V_n-V\|_{L_T^1H^{s-1}}.
 \end{aligned}
\]
For each fixed $K$, convergence of the initial data and
\eqref{eq:lower-field-convergence} imply
\[
 \limsup_{n\to\infty}(A_n+G_n)\leq C_R\varepsilon_K.
\]
We then let $K\to\infty$ to obtain $A_n+G_n\to0$.

Finally, \eqref{eq:ALE-difference} at $r=s$ gives
\[
 \begin{aligned}
\|\widetilde u_n-\widetilde u\|_{L_T^2H^s}
 +\|V_n-V\|_{L_T^2H^s}
 +\|Q_n-Q\|_{L_T^2H^s}
 +\|\partial_t(f_n-f)\|_{L_T^2H^{s-1/2}}
 \longrightarrow0.
 \end{aligned}
\]
This proves continuous dependence in the topology of Theorem~\ref{thm:local-target}.
On a longer common interval with uniform solution bounds and, in finite depth, a uniform positive lower bound on the depth, we obtain continuous dependence by iterating the short-time argument over finitely many subintervals.
\subsection{Positivity and continuation}\label{sec:continuation}

By the equation of the density perturbation \eqref{eq:ALE},
\[
\partial_t\varrho_f+V_f\cdot\nabla\varrho_f =\widetilde u_{d+1}.
\]
Hence the pulled-back total density $\widetilde\rho =\gamma(\varrho_f)+\widetilde g$ satisfies
\[
(\partial_t+V_f\cdot\nabla)\widetilde\rho=\gamma'(\varrho_f)\widetilde u_{d+1}+Q_f=0.
\]
The field $V_f\in L_T^1W^{1,\infty}$ is tangent to $\partial\mathcal S$, hence its flow preserves $\mathcal S$.
The total density is constant along characteristics.
Hence $\rho_0\geq c_\rho>0$ implies $\rho(t)\geq c_\rho$.
Together with existence, uniqueness, and continuous dependence, this proves Theorem~\ref{thm:local-target}.

The same construction gives the following continuation criterion.

\begin{corollary}[Continuation criterion]\label{cor:local-continuation}
Under the assumptions of Theorem~\ref{thm:local-target}, let a solution exist on $[0,T_*)$, with $T_*<\infty$, in the class \eqref{eq:main-regularity} on every compact subinterval.
Suppose that the background bounds in that theorem remain uniform and that
\[
 \sup_{t<T_*}\bigl(\|f(t)\|_{H^{s+1}} +\|g(t)\|_{H^s(\Omega_{f(t)})}\bigr)<\infty.
\]
In finite depth, assume also
\[
 \inf_{t<T_*,x}(f(t,x)-b(x))>0.
\]
Then the solution extends beyond $T_*$ in the same regularity class.
\end{corollary}

\begin{proof}
The Sobolev equivalence \eqref{eq:sobolev-pullback} gives a uniform bound for $(f(t),\widetilde g(t))$ in $H^{s+1}\times H^s$.
The construction of local solution in Section~\ref{sec:local-uniform-lifespan} then provides a time $\tau>0$, independent of $t_0<T_*$, on which the solution can be continued at $t_0$.
We choose $t_0$ so that $t_0<T_*<t_0+\tau$.
By uniqueness, the restarted solution agrees with the original solution on $[t_0,T_*)$ and extends it beyond $T_*$.
\end{proof}

\section{Vanishing surface tension limit}
\label{sec:vanishing-surface-tension}

Under the hypotheses of Theorem~\ref{thm:vanishing-surface-tension}, we establish estimates on a common time interval for $0<\sigma\leq1$.
Throughout this section, we assume $s>D/2+2$.

When $\|f\|_{H^s}\leq R$, the constants $C_R$ depend only on $R$, the fixed background and bottom, and the lower bounds for depth and stability.
In particular, they are independent of $\sigma$ and of the higher interface norms, which we keep explicit.
These bounds allow us to choose the same sufficiently small flattening parameter in \eqref{eq:explicit-flattening} for all solutions.

\subsection{Uniform tame elliptic and transport estimates}

We use the notations
\begin{equation} \label{FjGDef}
F_j(t)=\|f(t)\|_{H^{s+j}},\qquad G(t)=\|\widetilde g(t)\|_{H^s}, \qquad F_0(t)\leq R.
\end{equation}
We omit the time argument when it is clear from context.
The coefficients in \eqref{eq:matrix-A} satisfy
\begin{equation}\label{eq:vs-geometry}
 \begin{aligned}
 &\|A_f-A_*\|_{H^{s-1/2}}+\|D\Phi_f^{-T}-M_*\|_{H^{s-1/2}}\leq C_R,\\
 &\|A_f-A_*\|_{H^s}+\|D\Phi_f^{-T}-M_*\|_{H^s}\leq C_R(1+F_{1/2}).
 \end{aligned}
\end{equation}
Since $s-1/2>D/2+1$, the coefficient Lipschitz norms and ellipticity constants are uniformly controlled.
The next lemma separates the low-order bounds, whose constants depend only on $F_0\le R$, from the $H^s$ gradient estimate, which depends linearly on $F_{1/2}$.

\begin{lemma}[Tame elliptic estimates with $H^s$ interface control]
\label{lem:vs-tame-elliptic}
Let $s>D/2+1$ and assume \eqref{eq:vs-geometry}.
\begin{enumerate}
\item For the variational problem in Lemma~\ref{lem:strip-elliptic},
\begin{align}
 &\|\nabla v\|_{H^{s-1}} \leq C_R\bigl( \|h\|_{H^{s-1/2}}+\|W\|_{H^{s-1}}\bigr), \label{eq:vs-elliptic-low}\\
 &\|\nabla v\|_{H^s}\leq C_R\Bigl( \|h\|_{H^{s+1/2}}+\|W\|_{H^s} +(1+F_{1/2}) \bigl(\|h\|_{H^{s-1/2}}+\|W\|_{H^{s-1}}\bigr) \Bigr). \label{eq:vs-elliptic-high}
\end{align}
Moreover, the low-order estimate \eqref{eq:vs-elliptic-low} holds for every regularity index $0\leq r\leq s-1/2$ (with $h\in H^{r+1/2}(\cO)$ and $W\in H^r(\mathcal S)$), and in particular at $r=s-3/2$ and $r=s-1/2$, with constants $C_R$ depending only on $F_0\le R$.

\item For an energy solution of
$\operatorname{div}(A\nabla v)=k$ with zero top Dirichlet
data and homogeneous bottom conormal data, 
\begin{align}
 &\|\nabla v\|_{H^{s-1}} \leq C_R\bigl(\|k\|_{H^{s-2}}+\|\nabla v\|_{L^2}\bigr),
 \label{eq:vs-scalar-low}\\
 &\|\nabla v\|_{H^s}\leq C_R\Bigl( \|k\|_{H^{s-1}} +(1+F_{1/2}) \bigl(\|k\|_{H^{s-2}}+\|\nabla v\|_{L^2}\bigr)\Bigr).
 \label{eq:vs-scalar-high}
\end{align}
\end{enumerate}
\end{lemma}

\begin{proof}
Set $r_*=s-1/2>D/2$.
The first bound in \eqref{eq:vs-geometry} controls both $\|A-A_*\|_{H^{r_*}}$ and a positive H\"older norm of $A$ in terms of $R$.
Consequently, the partition of unity and patch size in the coefficient-freezing argument of Lemma~\ref{lem:strip-elliptic} depend only on $R$.
After absorbing the small coefficient oscillation,
that argument gives
\[
 \|\nabla v\|_{H^{r_*}} \leq C_R\bigl(\|h\|_{H^{r_*+1/2}}+\|W\|_{H^{r_*}} +\|\nabla v\|_{H^{r_*-\beta}}\bigr), \qquad  0<\beta<\min\{1,r_*-D/2\}
\]
We absorb the last term by interpolation and the energy estimate
\[
 \|\nabla v\|_{L^2}
 \leq C_R\bigl(\|h\|_{H^{1/2}}+\|W\|_{L^2}\bigr).
\]
For fixed $A$, interpolation of the linear solution
operator between the energy level and $r_*$ then yields
\[
 \|\nabla v\|_{H^r}
 \leq C_R\bigl(
   \|h\|_{H^{r+1/2}}+\|W\|_{H^r}\bigr),
 \qquad 0\leq r\leq r_*.
\]
This proves \eqref{eq:vs-elliptic-low} and the asserted
fractional-index estimates.

For the top-order bound \eqref{eq:vs-elliptic-high}, we use the second estimate in \eqref{eq:vs-geometry} on the patches chosen above.
Each localized coefficient difference $a=a_j$ and localized solution $w=w_j$ then satisfy
\[
 \|a\|_{L^\infty}\leq\eta,
 \qquad
 \|a\|_{H^s}\leq C_R(1+F_{1/2}),
\]
and hence
\[
 \|a\nabla w\|_{H^s}
 \leq C\eta\|\nabla w\|_{H^s}
      +C_R(1+F_{1/2})\|\nabla w\|_{H^{s-1}},
\]
since $s-1>D/2$.
We absorb the first term and estimate the extension and cutoff terms by the same product bound.
It follows that
\[
 \begin{aligned}
 \|\nabla v\|_{H^s}
 \leq C_R\Bigl(
\|h\|_{H^{s+1/2}}+\|W\|_{H^s}+(1+F_{1/2})
      \bigl(\|\nabla v\|_{H^{s-1}}
            +\|h\|_{H^{s-1/2}}\bigr)
 \Bigr).
 \end{aligned}
\]
Substituting \eqref{eq:vs-elliptic-low} proves
\eqref{eq:vs-elliptic-high}.

For a scalar source, the constant-coefficient estimate
\eqref{eq:constant-elliptic} requires $k$ in $H^{r-1}$.
The same localization at $r=s-1$, followed by interpolation
against $\|\nabla v\|_{L^2}$, yields
\eqref{eq:vs-scalar-low}.
At $r=s$, it gives
\[
 \|\nabla v\|_{H^s}
 \leq C_R\bigl(
   \|k\|_{H^{s-1}}
   +(1+F_{1/2})\|\nabla v\|_{H^{s-1}}\bigr).
\]
Combining these two estimates proves
\eqref{eq:vs-scalar-high}.
The localization preserves the conormal condition as in
Lemma~\ref{lem:strip-elliptic}, and approximation extends
the estimates to variational solutions.
\end{proof}

We apply \eqref{eq:vs-elliptic-low}--\eqref{eq:vs-elliptic-high} with $h=0$ and $W=-\widetilde g e_z$.
The trace and product estimates then give
\begin{equation}\label{eq:vs-bulk-trace}
 \|\mathcal B(f,\widetilde g)\|_{H^{s-1/2}}
 \leq C_R G(1+F_{1/2}).
\end{equation}
To control the normal in this estimate, we use
$\|\nabla f\|_{H^{s-1/2}}\leq F_{1/2}$,
while \eqref{eq:vs-elliptic-low} bounds the boundary pressure gradient in $L^\infty$ by $C_R G$.

To estimate the capillary velocity contribution $-\sigma(\bG[f]H(f))\circ\Phi_f$, we apply \eqref{eq:vs-elliptic-high} with boundary datum $h=\sigma H(f)$ and bulk flux $W=0$.
Sobolev composition applied to
$H(f)=-\operatorname{div}_x(\nabla f/\sqrt{1+|\nabla f|^2})$
gives
\[
 \|H(f)\|_{H^{s-1/2}}\le C_R F_{3/2},
 \qquad
 \|H(f)\|_{H^{s+1/2}}\le C_R F_{5/2}.
\]
Together with \eqref{eq:vs-elliptic-high} and
\eqref{eq:vs-geometry}, these estimates yield
\[
\|\sigma(\bG[f]H(f))\circ\Phi_f\|_{H^s}
 \leq C_R\sigma\bigl(F_{5/2}+(1+F_{1/2})F_{3/2}\bigr)\leq C_R\sigma F_{5/2},
\]
where the last inequality follows from
\[
 F_{1/2}F_{3/2}\leq F_0F_2\leq R F_{5/2}.
\]
Thus the $L_T^2H^s$ estimate for this contribution requires only the weighted norm $\sigma\|f\|_{L_T^2H^{s+5/2}}$.

Applying \eqref{eq:vs-elliptic-low}--\eqref{eq:vs-elliptic-high}
with $h=\sigma H(f)+\Gamma(f)$ and $W=-\widetilde g e_z$,
and using Sobolev composition for $\Gamma(f)$
and \eqref{eq:velocity-fixed}, we obtain
\begin{align}
 &\|\widetilde u\|_{H^{s-1}}
 \leq C_R(1+G+\sigma F_{3/2}),
 \label{eq:vs-u-low}\\
 &\|\widetilde u\|_{H^s}
 \leq C_R\bigl(
 (1+G)(1+F_{1/2})+\sigma F_{5/2}\bigr).
 \label{eq:vs-u-high}
\end{align}

The kinematic condition $f_t=N_f\cdot\widetilde u|_{z=0}$ and the trace estimates give
\[
 \begin{aligned}
 \|f_t\|_{H^{s-3/2}}
 \leq C_R\|\widetilde u\|_{H^{s-1}},\quad
 \|f_t\|_{H^{s-1/2}}
 \leq C_R\bigl(
 \|\widetilde u\|_{H^s}
 +F_{1/2}\|\widetilde u\|_{H^{s-1}}\bigr).
 \end{aligned}
\]
Combining these bounds with the time-extension estimate \eqref{eq:time-extension} and the product estimates in \eqref{eq:ALE}, we get
\begin{equation}\label{eq:vs-ALE-high}
 \begin{aligned}
 \|V_f\|_{H^s}+\|Q_f\|_{H^s}
       +\|f_t\|_{H^{s-1/2}}\leq C_R\bigl(
 (1+G)(1+F_{1/2})+\sigma F_{5/2}\bigr).
 \end{aligned}
\end{equation}
Here the additional capillary product is bounded by $\sigma F_{1/2}F_{3/2}\leq C_R\sigma F_{5/2}$ as above.

\subsection{Linearization of the surface operators}

The weighted $H^{s+1}$ interface estimate requires the equation for $\nabla_x f$.
Set the first order and third order operators
\[
 L_0(f)=\DN[f]\Gamma(f),\qquad
 L_3(f)=\DN[f]H(f).
\]
Let $\phi_{\Gamma(f)}$ be the harmonic extension of $\Gamma(f)$, and define
\[
B=\partial_y\phi_{\Gamma(f)}|_{\Sigma_f},\qquad c=\gamma(f)-B=\frac{\mathcal T(f)}{1+|\nabla f|^2}, \qquad W=\nabla_x\phi_{\Gamma(f)}|_{\Sigma_f}=c\nabla f.
\]
We separate the principal parts of the linearized operators and keep the higher interface norms explicit in the remainder estimates.

\begin{lemma}[Tame linearized surface operators]
\label{lem:vs-linearized}
Assume $s>D/2+2$, $F_0\le R$, $\mathcal T(f)\ge a_0>0$, and a common positive depth bound in finite depth.
For smooth directions $w$, with the bottom held fixed,
\begin{align}
 &D_fL_0(f)w =T_{\lambda_fc}w-T_W\cdot\nabla_xw+R_0(f)w,
 \\
 &\|R_0(f)w\|_{H^{s-1/2}} \le C_R(1+F_{1/2})\|w\|_{H^s}, \label{eq:vs-lin-gravity-bound}\\
 &D_fL_3(f)w= T_{\lambda_f\ell_f}w+R_3(f)w, \label{eq:vs-lin-cap}\\
 &\|R_3(f)w\|_{H^{s-3/2}} \le C_R\Bigl(  \|w\|_{H^{s+3/2-\alpha}} +(1+F_{3/2})\|w\|_{H^s}\Bigr), \label{eq:vs-lin-cap-bound}
\end{align}
where $0<\alpha<1/2$ is fixed.
Moreover, for $i=1,\ldots,d$ and $j=0,3$,
\[
 \partial_iL_j(f) =D_fL_j(f)\partial_i f+K_{j,i}(f),
\]
where $K_{j,i}=0$ in infinite depth, while in finite depth
\begin{equation}\label{eq:vs-bottom-bound}
 \|K_{0,i}(f)\|_{H^{s-1/2}} \le C_R(1+F_{1/2}), \qquad \|K_{3,i}(f)\|_{H^{s-3/2}}\le C_R.
\end{equation}
\end{lemma}

\begin{proof}
Let $\phi_a$ denote the harmonic extension of $a$, with homogeneous bottom Neumann data when present.
The proof of the lemma is divided into the following steps.

\emph{1. Shape derivative.}
The fixed-bottom shape derivative is
\begin{equation}\label{eq:vs-shape-gravity}
 D_fL_0(f)w =\DN[f](cw)-\operatorname{div}_x(Ww).
\end{equation}
Define
\[
 \mathcal R_f(a,w) =\DN[f](aw)-a\DN[f]w-w\DN[f]a, \qquad E_b=\DN[f]B+\operatorname{div}_xW.
\]
Since $c=\gamma(f)-B$, we obtain
\[
 D_fL_0(f)w =c\DN[f]w-W\cdot\nabla_xw +w\DN[f]\gamma(f)+\mathcal R_f(\gamma(f),w) -\mathcal R_f(B,w)-wE_b.
\]
\emph{2. Product remainder.}
For $a=B$ or $a=\gamma(f)-\gamma(0)$,
Lemma~\ref{lem:vs-tame-elliptic} and Sobolev composition give
\[
 \|a\|_{H^{s-1}}\le C_R,\qquad
 \|a\|_{H^{s-1/2}}\le C_R(1+F_{1/2}).
\]
Adding a constant to $a$ leaves $\mathcal R_f(a,w)$ unchanged.
The function $Z=\phi_{aw}-\phi_a\phi_w$ has zero top trace, homogeneous bottom Neumann data, and
\[
 \Delta Z=-2\nabla \phi_a\cdot\nabla\phi_w,\qquad
 N_f\cdot\nabla Z|_{\Sigma_f}=\mathcal R_f(a,w).
\]
For $Z\circ\Phi_f$, the flattened equation is
\[
 \operatorname{div}(A_f\nabla(Z\circ\Phi_f))=k,\qquad
k=-2J_f(\nabla\phi_a\cdot\nabla\phi_w)\circ\Phi_f.
\]
The lower-index elliptic estimates and tame multiplication give
\[
 \|k\|_{H^{s-2}}\le C_R\|w\|_{H^s},\qquad \|k\|_{H^{s-1}} \le C_R(1+F_{1/2})\|w\|_{H^s}.
\]
The energy estimate \eqref{eq:elliptic-energy} with $W=0$ and maximum principle also give, with the norms below taken over $\Omega_f$,
\begin{align*}
 \|\nabla Z\|_{L^2}\le \|\nabla\phi_{aw}\|_{L^2} +\|\phi_a\|_{L^\infty}\|\nabla\phi_w\|_{L^2}+\|\phi_w\|_{L^\infty}\|\nabla\phi_a\|_{L^2}\le C_R\|w\|_{H^s}.
\end{align*}
Applying \eqref{eq:vs-scalar-low}--\eqref{eq:vs-scalar-high} therefore yields
\[
 \|\nabla (Z\circ\Phi_f)\|_{H^{s-1}}\le C_R\|w\|_{H^s},\qquad 
 \|\nabla(Z\circ\Phi_f)\|_{H^s}\le C_R(1+F_{1/2})\|w\|_{H^s}.
\]
The trace theorem and tame multiplication imply
\begin{align}
 \|\mathcal R_f(a,w)\|_{H^{s-3/2}}
 &\le C_R\|w\|_{H^s}, \label{eq:vs-product-defect-low}\\
 \|\mathcal R_f(a,w)\|_{H^{s-1/2}}
 &\le C_R(1+F_{1/2})\|w\|_{H^s}. \label{eq:vs-product-defect}
\end{align}
Here the high-order conormal coefficient is multiplied by
$\nabla(Z\circ\Phi_f)|_{z=0}$ in $L^\infty$, which is controlled by the preceding $H^{s-1}$ bound.

\emph{3. Bottom contributions.}
The function $(\phi_B-\partial_y\phi_{\Gamma(f)})$ is harmonic and vanishes at the top.
$\phi_{\Gamma(f)}$ is a harmonic function, so that
\[
 N_f\cdot\nabla (\phi_B-\partial_y\phi_{\Gamma(f)})|_{\Sigma_f}=E_b.
\]
In infinite depth, uniqueness in the energy class gives $\phi_B-\partial_y\phi_{\Gamma(f)}=0$, hence $E_b=0$.
Horizontal translation invariance also gives $K_{j,i}=0$.

In finite depth, we write $N_b=(-\nabla b,1)$.
The bottom flux of $(\phi_B-\partial_y\phi_{\Gamma(f)})$ is
\[
 N_b\cdot\nabla (\phi_B-\partial_y\phi_{\Gamma(f)}) =-N_b\cdot\nabla\partial_y\phi_{\Gamma(f)}
 \qquad\text{on }y=b(x).
\]
For the spatial derivative corrections, we fix $i$.
Translation covariance gives
\[
 \partial_iL_j(f) =D_fL_j(f)\partial_i f+K_{j,i},\qquad K_{j,i}=N_f\cdot\nabla\psi_j|_{\Sigma_f},
\]
where $\psi_j$ is harmonic with zero top trace and has bottom flux
\[
 N_b\cdot\nabla\psi_j=
 \begin{cases}
 \nabla_x(\partial_i b)\cdot\nabla_x\phi_{\Gamma(f)} -(\partial_i b)N_b\cdot\partial_y\nabla\phi_{\Gamma(f)},&j=0,\\
 \nabla_x(\partial_i b)\cdot\nabla_x\phi_{H(f)} -(\partial_i b)N_b\cdot\partial_y\nabla\phi_{H(f)},&j=3,
 \end{cases}
 \qquad\text{on }y=b(x).
\]

The positive depth bound separates a fixed neighborhood of the bottom from $\Sigma_f$.
Elliptic regularity there controls the bottom flux norms by
\[
\|\nabla\phi_{\Gamma(f)}\|_{L^2(\Omega_f)}  +\|\nabla\phi_{H(f)}\|_{L^2(\Omega_f)} \le C_R.
\]
Extending these fluxes to the fluid domain and applying
Lemma~\ref{lem:vs-tame-elliptic}, we prove \eqref{eq:vs-bottom-bound} using the trace theorem and
\[
 \|E_b\|_{H^{s-3/2}}\le C_R,\qquad \|E_b\|_{H^{s-1/2}}\le C_R(1+F_{1/2}).
\]

\emph{4. Estimate of the principal parts.}
The coefficient bounds in Step~2 also hold for
$c-\gamma(0)$ and $W$; their Lipschitz norms are bounded by $C_R$.
The paralinearization estimates give
\[
 \|\DN[f]w-T_{\lambda_f}w\|_{H^{s-3/2}}\le C_R\|w\|_{H^s},\qquad
 \|\DN[f]w-T_{\lambda_f}w\|_{H^{s-1/2}}\le C_R(1+F_{1/2})\|w\|_{H^s}.
\]
For any surface remainder $\mathfrak r$, tame multiplication gives
\begin{equation}\label{eq:vs-remainder-multiplication}
 \|c\mathfrak r\|_{H^{s-1/2}}
 \le C_R\Bigl(\|\mathfrak r\|_{H^{s-1/2}}
   +(1+F_{1/2})\|\mathfrak r\|_{H^{s-3/2}} \Bigr).
\end{equation}
Applying this to $\mathfrak r=\DN[f]w-T_{\lambda_f}w$, and using Bony's decomposition and symbol composition, we obtain
\[
\begin{aligned}
 \|c\DN[f]w-W\cdot\nabla_xw
-T_{\lambda_fc}w+T_W\cdot\nabla_xw\|_{H^{s-1/2}}\le C_R(1+F_{1/2})\|w\|_{H^s}.
\end{aligned}
\]
Moreover, the tame elliptic estimate with datum $\gamma(f)-\gamma(0)$ yields
\[
 \|\DN[f]\gamma(f)\|_{H^{s-1/2}} \le C_R(1+F_{1/2}).
\]
Together with Steps~2--3, these estimates control the right-hand side of \eqref{eq:vs-shape-gravity} and prove \eqref{eq:vs-lin-gravity-bound}.

For the capillary operator, we apply
\eqref{eq:cap-para-difference} at $p=s$ to $f+\tau w$ and $f$.
For sufficiently small $\tau$, the perturbed interfaces satisfy common base bounds.
Dividing by $\tau$ and letting $\tau\to0$ gives
\[
 \|D_fL_3(f)w-T_{\lambda_f\ell_f}w\|_{H^{s-3/2}}
 \le C_R\Bigl( \|w\|_{H^{s+3/2-\alpha}} +(1+F_{3/2})\|w\|_{H^s}\Bigr).
\]
This proves \eqref{eq:vs-lin-cap}--\eqref{eq:vs-lin-cap-bound}.
\end{proof}

\subsection{The unweighted and weighted energy estimates}

We combine the unweighted interface and density energies with a weighted higher-order interface energy:
\begin{equation}\label{eq:vs-energy-def}
 E_s^\sigma=F_0^2+G^2+\sigma F_1^2,
 \qquad D_s^\sigma= F_{1/2}^2+\sigma F_{3/2}^2+\sigma^2F_{5/2}^2,
\end{equation}
where $F_j$ and $G$ are defined in \eqref{FjGDef}.
Stability and smallness of the density perturbation give the unweighted interface estimate.
The weighted energy gives the additional regularity needed to control the bulk velocity uniformly in surface tension.

\begin{proposition}[Uniform interface energies]
\label{prop:vs-energies}
Assume $s>D/2+2$.
On a set with $F_0\leq R$, $\mathcal T(f)\geq a_0>0$,
and a common positive depth bound, there exists $\varepsilon_b>0$, independent of $0<\sigma\leq1$,
such that smooth solutions with $G\leq\varepsilon_b$
satisfy
\begin{equation}\label{eq:vs-low-energy}
 \frac{\dd}{\dd t}F_0^2
 +c_R\bigl(F_{1/2}^2+\sigma F_{3/2}^2\bigr) \leq C_R(F_0^2+G^2).
\end{equation}
Moreover, for $E_1=\sigma F_1^2$,
\begin{equation}\label{eq:vs-high-energy}
 \frac{\dd}{\dd t}E_1 +c_R\bigl(\sigma F_{3/2}^2+\sigma^2F_{5/2}^2\bigr)\leq C_R\bigl(1+F_{1/2}^2+\sigma F_{3/2}^2\bigr)E_1+C_R(G^2+\sigma)(1+F_{1/2}^2).
\end{equation}
All constants have the base-geometry dependence specified at the start of this section and are independent of $\sigma$.
\end{proposition}

\begin{proof}
We write the first equation in \eqref{e:IPM} as
\[
 \partial_t f+L_0(f)+\sigma L_3(f)=-\mathcal B.
\]

\emph{1. The base energy.}
We use the stable reduction of $L_0$ in
\cite[Section~4.2]{LatoccaNguyen2025} and \eqref{eq:cap-para-single} at $p=s$, whose constant depends only on $F_0\le R$ and the fixed geometry.
The principal symbols satisfy
\[
 \lambda_fc\geq c_R|\xi|, \qquad \lambda_f\ell_f\geq c_R|\xi|^3.
\]
The symbolic energy argument of \cite[Proposition~3.1]{FlynnNguyen2020}, including the transport term, then gives
\[
 \frac{\dd}{\dd t}F_0^2 +c_R\bigl(F_{1/2}^2+\sigma F_{3/2}^2\bigr)
 \leq C_R(1+\sigma)F_0^2 +C_RF_{1/2}F_{1/2-\alpha}+C_R\sigma F_{3/2}F_{3/2-\alpha} +C_RG(1+F_{1/2})F_{1/2},
\]
where $0<\alpha<1/2$ is fixed, and the last term follows from \eqref{eq:vs-bulk-trace}.
The $F_0^2$ term includes the low-frequency errors arising from the use of inhomogeneous Sobolev norms.

For every $\eta>0$, interpolation gives
\[
 C_RF_{1/2}F_{1/2-\alpha}+C_R\sigma F_{3/2}F_{3/2-\alpha}\leq \eta \bigl(F_{1/2}^2+\sigma F_{3/2}^2\bigr) +C_{R,\eta}(1+\sigma)F_0^2.
\]
The density term satisfies
\[
 C_RG(1+F_{1/2})F_{1/2}\leq (C_RG+\eta)F_{1/2}^2+C_{R,\eta}G^2.
\]
Choosing $\eta$ and then $\varepsilon_b$ sufficiently small, independently of $\sigma$, allows us to absorb the dissipative error terms and obtain \eqref{eq:vs-low-energy}.

\emph{2. The weighted higher energy.}
Lemma~\ref{lem:vs-linearized} gives
\[
 \partial_i f_t +(T_{\lambda_fc}-T_W\cdot\nabla+\sigma T_{\lambda_f\ell_f})\partial_i f = -R_0(f)\partial_i f-\sigma R_3(f)\partial_i f-\partial_i\mathcal B-K_{0,i}-\sigma K_{3,i}.
\]
We test this equation against $\sigma\partial_i f$ at the $H^s$ level and take the real part.
The time derivative contributes
$\frac{\sigma}{2}\frac{\dd}{\dd t}\|\partial_i f\|_{H^s}^2$.
As in Step~1, the stable and capillary principal operators contribute the dissipation
\[
 c_R\sigma\|\partial_i f\|_{H^{s+1/2}}^2
 +c_R\sigma^2\|\partial_i f\|_{H^{s+3/2}}^2.
\]
For every $\eta>0$, the transport contribution and the lower-order symbolic errors are bounded in absolute value by
\[
 \eta\sigma\|\partial_i f\|_{H^{s+1/2}}^2 +\eta\sigma^2\|\partial_i f\|_{H^{s+3/2}}^2 +C_{R,\eta}\sigma\|\partial_i f\|_{H^s}^2.
\]
This follows by interpolation as in Step~1, using $0<\sigma\leq1$.
For sufficiently small $\eta$, the first two terms are absorbed into the principal dissipation.

For the remainders $R_0(f)\partial_i f$ and $\sigma R_3(f)\partial_i f$, Lemma~\ref{lem:vs-linearized}, Sobolev duality, interpolation, and Young's inequality give
\[
 \begin{aligned}
 &\sigma\bigl|\langle R_0(f)\partial_i f,\partial_i f\rangle_{H^s}\bigr| +\sigma^2\bigl|\langle R_3(f)\partial_i f,\partial_i f\rangle_{H^s}\bigr|\\
 \leq & \eta\sigma\|\partial_i f\|_{H^{s+1/2}}^2
+\eta\sigma^2\|\partial_i f\|_{H^{s+3/2}}^2+C_{R,\eta}\bigl(1+F_{1/2}^2+\sigma F_{3/2}^2\bigr) \sigma\|\partial_i f\|_{H^s}^2.
 \end{aligned}
\]
Here we interpolate $\|\partial_i f\|_{H^{s+3/2-\alpha}}$ between $\|\partial_i f\|_{H^s}$ and $\|\partial_i f\|_{H^{s+3/2}}$ while keeping the weight $\sigma^2$, and use the inequality
\[
 \sigma^2(1+F_{3/2}^2) \leq \sigma(1+\sigma F_{3/2}^2).
\]

By \eqref{eq:vs-bulk-trace}, Sobolev duality and Young's
inequality give
\begin{equation}\label{eq:vs-high-bulk-pairing}
 \begin{aligned}
 \sigma|\langle\partial_i\mathcal B,\partial_i f\rangle_{H^s}|
 \leq C_R\sigma G(1+F_{1/2}) \|\partial_i f\|_{H^{s+3/2}}\leq \eta\sigma^2\|\partial_i f\|_{H^{s+3/2}}^2+C_{R,\eta}G^2(1+F_{1/2}^2).
 \end{aligned}
\end{equation}
This estimate requires only $\widetilde g\in H^s$.
Similarly, \eqref{eq:vs-bottom-bound} and $\sigma\leq1$ yield
\[
 \begin{aligned}
 \sigma|\langle K_{0,i},\partial_i f\rangle_{H^s}| +\sigma^2|\langle K_{3,i},\partial_i f\rangle_{H^s}|\leq \eta\sigma\|\partial_i f\|_{H^{s+1/2}}^2
+\eta\sigma^2\|\partial_i f\|_{H^{s+3/2}}^2 +C_{R,\eta}\sigma(1+F_{1/2}^2).
 \end{aligned}
\]

We choose $\eta>0$ sufficiently small to absorb the dissipative errors, sum over $i=1,\ldots,d$, and add $\sigma$ times the base energy inequality \eqref{eq:vs-low-energy}.
The inhomogeneous Sobolev norm identity
\[
 F_{q+1}^2 =F_q^2+\sum_{i=1}^d\|\partial_i f\|_{H^{s+q}}^2,
 \qquad q=0,\frac12,\frac32,
\]
then recovers the full time derivative of $E_1=\sigma F_1^2$ and the dissipative terms $c_R(\sigma F_{3/2}^2+\sigma^2F_{5/2}^2)$.
The additional right-hand terms satisfy $\sigma F_0^2\leq E_1$ and $\sigma G^2\leq G^2$.
This proves \eqref{eq:vs-high-energy}.

Spatial regularization and the tame estimates extend these inequalities, in integrated form, to the local energy solutions.
\end{proof}

\subsection{The common lifespan}

Assume $\|f_0\|_{H^{s+1}}\leq R_0$, $R>2(R_0+1)$, and fix initial lower bounds for stability and depth:
\[
 \mathcal T(f_0)\geq 2a_0,\qquad f_0-b\geq 2h_0.
\]
Again all depth conditions are omitted in infinite depth.
We first choose $\varepsilon_b>0$ small enough for Proposition~\ref{prop:vs-energies} and the zero-surface-tension difference estimates in Section~\ref{sec:limit-comparison}.
With $\varepsilon_b$ fixed, we choose $\varepsilon>0$ so that the zero-surface-tension local theory applies and $C_s\varepsilon\leq\varepsilon_b/4$.
Here $C_s$ includes the norm of the extension operator used in the transport estimate \eqref{eq:transport-high}.
Both smallness parameters depend only on the fixed data bounds and the lower bounds for stability and depth.

\emph{Step 1: Uniform estimates.}
For each fixed $0<\sigma\leq1$, we start with the solution from
Theorem~\ref{thm:local-target} and bootstrap
\[
 F_0\leq R,\qquad G\leq\varepsilon_b,\qquad
 \mathcal T(f)\geq a_0,\qquad f-b\geq h_0.
\]
As long as these bounds hold, \eqref{eq:vs-low-energy} gives
\[
 F_0(t)^2+c_R\int_0^t
       (F_{1/2}^2+\sigma F_{3/2}^2)\dd\tau
 \leq R_0^2+C_Rt(R^2+\varepsilon_b^2).
\]
Thus there are $T_*\leq1$ and $M_0$, independent of $\sigma$, such that, for $T\leq T_*$,
\[
 \sup_{t\leq T}F_0(t)\leq R/2,\qquad
 \int_0^T(F_{1/2}^2+\sigma F_{3/2}^2)\dd t\leq M_0.
\]
In \eqref{eq:vs-high-energy}, the coefficient and forcing satisfy
\[
 \int_0^T(1+F_{1/2}^2+\sigma F_{3/2}^2)\dd t\leq1+M_0,
 \qquad \int_0^T(G^2+\sigma)(1+F_{1/2}^2)\dd t
 \leq(\varepsilon_b^2+1)(1+M_0).
\]
Since $\sigma F_1(0)^2\leq R_0^2$, Gronwall's inequality yields
\begin{equation}\label{eq:vs-high-integrated}
 \sup_{t\leq T}\sigma F_1(t)^2 +\int_0^T\sigma^2F_{5/2}^2\dd t\leq M_1,
\end{equation}
where $M_1$ depends only on the fixed bounds.

Equations \eqref{eq:vs-u-high}--\eqref{eq:vs-ALE-high} consequently give
\begin{equation}\label{eq:vs-L2-fields}
 \|\widetilde u^\sigma\|_{L_T^2H^s}
 +\|V^\sigma\|_{L_T^2H^s}
 +\|Q^\sigma\|_{L_T^2H^s}
 +\|f_t^\sigma\|_{L_T^2H^{s-1/2}}\leq M_2,
\end{equation}
where
\[
 M_2=C_R\bigl((1+\varepsilon_b)(1+\sqrt{M_0}) +\sqrt{M_1}\bigr).
\]
We keep $M_0,M_1,M_2$ fixed in the subsequent reductions of $T$.

\emph{Step 2: Closing the bootstrap.}
By \eqref{eq:transport-high} and \eqref{eq:vs-L2-fields},
\[
 \sup_{t\leq T}G(t)\leq C_s e^{C\sqrt T M_2}(\varepsilon+C\sqrt T M_2),
 \qquad
 \sup_{t\leq T}\|f^\sigma(t)-f_0\|_{H^{s-1}}\leq\sqrt T M_2.
\]
The Dirichlet--Neumann mapping and difference estimates give
\[
 \|\mathcal T(f)-\mathcal T(k)\|_{H^{s-2}}\leq C_R\|f-k\|_{H^{s-1}}.
\]
To see this, we split
\[
 \DN[f]\Gamma(f)-\DN[k]\Gamma(k)=\DN[f](\Gamma(f)-\Gamma(k))
   +(\DN[f]-\DN[k])\Gamma(k).
\]
The mapping bound controls the first term, and \cite[Theorem~2.7]{FlynnNguyen2020} controls the second with datum $\Gamma(k)\in H^s$ and output in $H^{s-2}$.
Since $s-2>d/2$, it follows that
\[
 \sup_{t\leq T}\Bigl(
 \|f^\sigma(t)-f_0\|_{L^\infty}
 +\|\mathcal T(f^\sigma(t))-\mathcal T(f_0)\|_{L^\infty}
 \Bigr)\leq C_R\sqrt T M_2.
\]
Taking $T$ smaller if necessary, we obtain
\[
 F_0\leq R/2,\qquad G\leq\varepsilon_b/2,\qquad
 \mathcal T(f^\sigma)\geq3a_0/2,\qquad
 f^\sigma-b\geq3h_0/2.
\]
These strictly improve all bootstrap bounds.

\emph{Step 3: Continuation and the zero-surface-tension solution.}
For each fixed $\sigma>0$, \eqref{eq:vs-high-integrated} gives $F_1\leq\sqrt{M_1/\sigma}$ throughout the bootstrap interval.
This estimate and the improved density and depth bounds allow us to apply the continuation criterion in Section~\ref{sec:continuation}.
Hence, every capillary solution exists on $[0,T]$.

For $\sigma=0$, we use the fixed-index construction in Lemma~\ref{app:ln-local}.
At index $s$, the required background multiplier and difference bounds follow from Remark~1.3 and Proposition~E.1 of that paper using derivatives of $\gamma$ only through order $\lceil s\rceil+2$.
The resulting solution satisfies
\[
 f^0\in C_TH^s\cap L_T^2H^{s+1/2},
 \qquad \widetilde g^0\in C_TH^s.
\]
We restrict $T$ to the lifespan of $(f^0,\widetilde g^0)$.
The $\sigma=0$ versions of the base energy and transport estimates above show that this solution obeys the same bounds on the resulting common interval $[0,T]$.

Thus, on a common time interval $[0,T]$ depending only on $R_0,a_0,h_0$ and the fixed background and bottom geometry, we have
\[
 \sup_{0\leq\sigma\leq1}\Big( \sup_{t\leq T}E_s^\sigma(t) +\int_0^T D_s^\sigma(t)\dd t\Big)\leq C,
\]
where the weighted high-order terms in \eqref{eq:vs-energy-def} are omitted when $\sigma=0$.

\subsection{The capillary residuals}

On the common interval $[0,T]$, we compare each capillary solution with the zero-surface-tension solution.
The uniform weighted bounds control the capillary contributions to the surface equation, velocity, transport field, and density source at the lower comparison indices.

Let $(U_0,V_0,Q_0)(f,\widetilde g)$ denote the vector fields on $\mathcal S$ without surface-tension, with $U_0$ the pulled-back physical velocity.
In the definition of $V_0$, the coordinate time derivative is computed using
\[
 f_t=-L_0(f)-\mathcal B(f,\widetilde g).
\]
The operators $(U_0,V_0,Q_0)$ are defined for any sufficiently regular pair $(f,\widetilde g)$, in particular for the capillary solutions $(f^\sigma,\widetilde g^\sigma)$.

For a capillary solution, linearity of the pressure problem and of $\mathscr E$ gives
\[
 \widetilde u^\sigma
 =U_0(f^\sigma,\widetilde g^\sigma)+U_{\rm cap}^\sigma,
 \qquad U_{\rm cap}^\sigma =-\sigma(\bG[f^\sigma]H(f^\sigma))\circ\Phi_{f^\sigma},
\]
and
\begin{equation}\label{eq:vs-ALE-split}
 \begin{aligned}
 V^\sigma
 &=V_0(f^\sigma,\widetilde g^\sigma)+V_{\rm cap}^\sigma,\quad
 V_{{\rm cap},x}^\sigma=U_{{\rm cap},x}^\sigma,\quad
 V_{{\rm cap},z}^\sigma
 =J_f^{-1}\bigl((-\nabla_x\varrho_f,1)\cdot U_{\rm cap}^\sigma +\sigma\mathscr E L_3(f^\sigma)\bigr),\\
 Q^\sigma
 &=Q_0(f^\sigma,\widetilde g^\sigma)+Q_{\rm cap}^\sigma,\qquad
 Q_{\rm cap}^\sigma =-\gamma'(\varrho_f)U_{{\rm cap},d+1}^\sigma.
 \end{aligned}
\end{equation}
Here $J_f$ and $\varrho_f$ are evaluated at $f=f^\sigma$.
The term $+\sigma\mathscr EL_3(f^\sigma)$ in $V_{{\rm cap},z}^\sigma$ comes from substituting $f_t^\sigma=-L_0(f^\sigma)-\mathcal B(f^\sigma,\widetilde g^\sigma)-\sigma L_3(f^\sigma)$ into $-\partial_t\varrho_{f^\sigma}=-\mathscr E f_t^\sigma$.

The low-order elliptic estimate \eqref{eq:vs-elliptic-low} and the trace theorem yield
\[
 \|U_{\rm cap}^\sigma\|_{H^{s-1}}
 +\|\sigma L_3(f^\sigma)\|_{H^{s-3/2}}
 \leq C_R\sigma\|H(f^\sigma)\|_{H^{s-1/2}} \leq C_R\sigma F_{3/2}.
\]
The geometric and background coefficients are bounded on $H^{s-1}$, and $\mathscr E:H^{s-3/2}\to H^{s-1}$ is also bounded.
The common energy bound therefore gives
\begin{equation}\label{eq:vs-small-residual}
 \begin{aligned}
 &\|U_{\rm cap}^\sigma\|_{L_T^2H^{s-1}} +\|V_{\rm cap}^\sigma\|_{L_T^2H^{s-1}} +\|Q_{\rm cap}^\sigma\|_{L_T^2H^{s-1}} +\|\sigma L_3(f^\sigma)\|_{L_T^2H^{s-3/2}}\\
 &\qquad\leq C_R\sigma\|f^\sigma\|_{L_T^2H^{s+3/2}} \leq C\sqrt{\sigma}.
 \end{aligned}
\end{equation}

\subsection{The coupled difference estimate}\label{sec:limit-comparison}

Set the difference
\[
 \delta f=f^\sigma-f^0,\qquad
 \delta\widetilde g=\widetilde g^\sigma-\widetilde g^0,
\]
and denote their norms by
\[
 X=\|\delta f\|_{C_TH^{s-1}} +\|\delta f\|_{L_T^2H^{s-1/2}}, \qquad Y=\|\delta\widetilde g\|_{C_TH^{s-1}}.
\]
All constants below are independent of $\sigma$.

Subtracting the first equation in \eqref{e:IPM}
for the solutions with surface tensions $\sigma$ and $0$ gives
\[
 \begin{aligned} \partial_t\delta f +L_0(f^\sigma)-L_0(f^0) +\mathcal B(f^\sigma,\widetilde g^\sigma) -\mathcal B(f^0,\widetilde g^0)=-\sigma L_3(f^\sigma).
 \end{aligned}
\]
We apply the energy argument of
\cite[Proposition~5.4]{LatoccaNguyen2025}, recalled in
Lemma~\ref{app:ln-surface-difference}, at index $s-1$.
The additional forcing is estimated by
\[
 |\langle\sigma L_3(f^\sigma),\delta f\rangle_{H^{s-1}}|\leq
 \|\sigma L_3(f^\sigma)\|_{H^{s-3/2}}
 \|\delta f\|_{H^{s-1/2}}\leq \eta\|\delta f\|_{H^{s-1/2}}^2 +C_\eta\|\sigma L_3(f^\sigma)\|_{H^{s-3/2}}^2.
\]
Our choice of $\varepsilon_b$ allows us to absorb the density-dependent dissipation term.
The common interface bounds control the remaining Gronwall coefficient.
The bound \eqref{eq:vs-small-residual} then gives
\begin{equation}\label{eq:vs-surface-difference}
 X\leq C\bigl( \|\delta f(0)\|_{H^{s-1}}+\sqrt T\,Y+\sqrt\sigma \bigr).
\end{equation}

For $A_0=U_0,V_0,Q_0$, we again set the difference
\[
 \delta A_0=A_0(f^\sigma,\widetilde g^\sigma) -A_0(f^0,\widetilde g^0).
\]
Let $v^0=q^0\circ\Phi_{f^0}$ be the pulled-back pressure of the zero-surface-tension solution.
To bound $\delta U_0$, we apply the low-order elliptic estimate \eqref{eq:vs-elliptic-low} to the difference of the zero-surface-tension pressure equations.
As in the proof of Lemma~\ref{lem:elliptic-comparison}, the top Dirichlet datum is $\Gamma(f^\sigma)-\Gamma(f^0)$ and the divergence source includes the coefficient difference $-(A_{f^\sigma}-A_{f^0})\nabla v^0$.
Equation \eqref{eq:velocity-fixed} then controls $\|\delta U_0\|_{H^{s-1}}$.
Since
\[
 L_0(f)+\mathcal B(f,\widetilde g) =-N_f\cdot U_0(f,\widetilde g)|_{z=0},
\]
the trace theorem and Sobolev product estimates control
$\delta(L_0+\mathcal B)$ in $H^{s-3/2}$.
The extension estimate \eqref{eq:time-extension}
at $r=s-1$ controls
$\mathscr E\delta(L_0+\mathcal B)$ in $H^{s-1}$.
Using the zero-surface-tension field formulas in \eqref{eq:ALE},
we obtain
\begin{equation}\label{eq:vs-zero-field-difference}
 \begin{aligned}
 \|\delta U_0\|_{H^{s-1}}
 +\|\delta V_0\|_{H^{s-1}}
 +\|\delta Q_0\|_{H^{s-1}}\leq C_R\bigl(
 \|\delta f\|_{H^{s-1/2}}
 +\|\delta\widetilde g\|_{H^{s-1}}\bigr).
 \end{aligned}
\end{equation}

Subtracting the transport equations gives
\[
 \partial_t\delta\widetilde g
 +V^\sigma\cdot\nabla\delta\widetilde g
 =
 \delta Q_0+Q_{\rm cap}^\sigma
 -(\delta V_0+V_{\rm cap}^\sigma)\cdot\nabla\widetilde g^0.
\]
Since $H^{s-1}(\mathcal S)$ is an algebra,
\[
 \|(\delta V_0+V_{\rm cap}^\sigma)
                \cdot\nabla\widetilde g^0\|_{H^{s-1}}
 \leq C\|\widetilde g^0\|_{H^s}
       \bigl(\|\delta V_0\|_{H^{s-1}}
             +\|V_{\rm cap}^\sigma\|_{H^{s-1}}\bigr).
\]
We apply the $H^{s-1}$ transport estimate to this equation and use \eqref{eq:vs-zero-field-difference} and \eqref{eq:vs-small-residual} to bound its right-hand side.
Since $\|V^\sigma\|_{L_T^1H^s}\leq C\sqrt T$ by \eqref{eq:vs-L2-fields}, we obtain
\begin{equation}\label{eq:vs-density-difference}
 Y\leq C\bigl(
 \|\delta\widetilde g(0)\|_{H^{s-1}}
 +TY+\sqrt T\,X+\sqrt{T\sigma}
 \bigr).
\end{equation}

For identical initial data, combining
\eqref{eq:vs-surface-difference} and
\eqref{eq:vs-density-difference} yields
\[
 Y\leq CTY+C\sqrt{T\sigma}.
\]
Choosing $T>0$ small enough that $CT\leq 1/2$ allows us to absorb $CTY$ into the left-hand side and substitute the resulting bound into \eqref{eq:vs-surface-difference} to find
\[
 Y\leq C\sqrt{T\sigma},
 \qquad X+Y\leq C\sqrt\sigma.
\]
The field decompositions and
\eqref{eq:vs-zero-field-difference} also yield
\[
 \begin{aligned}
 \|\widetilde u^\sigma-\widetilde u^0\|_{L_T^2H^{s-1}}
 +\|V^\sigma-V^0\|_{L_T^2H^{s-1}}\leq C\bigl(X+\sqrt T\,Y+\sqrt\sigma\bigr)
 \leq C\sqrt\sigma.
 \end{aligned}
\]

Finally, interpolation with the uniform base bounds gives, for
$s-1\leq r<s$,
\[
 \begin{aligned}
 \|\delta f\|_{C_TH^r}
 +\|\delta f\|_{L_T^2H^{r+1/2}}
 +\|\delta\widetilde g\|_{C_TH^r}\leq C_r(X+Y)^{s-r}
 \leq C_r\sigma^{(s-r)/2}.
 \end{aligned}
\]

\begin{remark}[Comparison with the scalar limit]
\label{rem:vs-scalar-check}
If $\gamma>0$ is constant and $g(0)=0$, then $g=0$ and the system
reduces to the one-phase Muskat problem considered in
\cite[Theorem~1.1]{FlynnNguyen2020}.
At the flat infinite-depth equilibrium, the linearized equation is
\[
 \partial_tf+\gamma|D|f+\sigma|D|^3f=0.
\]
Its homogeneous energy identity is
\[
 \begin{aligned}
\frac12\frac{\dd}{\dd t}
 \bigl(\|f\|_{\dot H^s}^2
       +\sigma\|f\|_{\dot H^{s+1}}^2\bigr)
 +\gamma\|f\|_{\dot H^{s+1/2}}^2
 +(1+\gamma)\sigma\|f\|_{\dot H^{s+3/2}}^2
 +\sigma^2\|f\|_{\dot H^{s+5/2}}^2=0.
 \end{aligned}
\]
This identity motivates the three $\sigma$-weighted dissipative norms in \eqref{eq:vs-energy-def}.
For a flat bottom $b=-\ell_0$, the linearized equation becomes
\[
 \partial_tf+
 |D|\tanh(\ell_0|D|)\bigl(\gamma+\sigma|D|^2\bigr)f=0.
\]
It has the same high-frequency dissipation; low frequencies are
controlled by the energy on a finite time interval.
\end{remark}

\section{Stratified solutions: density coordinates and linear estimates}\label{sec:global-linear-estimates}

Throughout this section and Section~\ref{sec:global-estimates}, we work in the flat strip setting of Theorem~\ref{thm:global}, with the profile \eqref{eq:global-profile} and boundary densities \eqref{eq:global-boundary-density}.
We use the density coordinate map $\Phi_\eta$ and the graph-flattening map $\Phi_f$ introduced in Subsubsection~\ref{sec:global-intro}.
We first explain how the estimates in these two sections give global well-posedness.

\textbf{Density coordinates and weak bulk damping.}
The density coordinate map fixes each density level: $\rho\circ\Phi_\eta=\gamma(z)$, with $\eta|_{z=0}=f$ and $\eta|_{-1}=0$.
In these coordinates, $\eta$ describes the location of each density level, and the modified pressure $\Pi$ in \eqref{eq:global-Pi-definition} separates the hydrostatic contribution from the velocity.
Since we consider the strict stratification, the bulk part of the relative energy \eqref{eq:global-relative-energy} is coercive, and the associated dissipation is the squared $L^2$ norm of the velocity with Jacobian weight $J_\eta$.
However, the linearized bulk evolution has no uniform damping at high vertical frequencies.
An energy bound alone therefore does not give the time-integrable control of the velocity and its horizontal derivatives needed to bound $\eta$ uniformly in time.
We address this difficulty by combining an anisotropic energy estimate with a separate bound for $\int_0^T\mathcal W_d(U(t))\,\dd t$ that is independent of $T$.

\textbf{A velocity estimate integrable in time.}
The linear pressure problem \eqref{eq:global-flat-pressure} couples the free surface to the bulk unknown $\eta$.
We treat the geometric and capillary remainders as prescribed sources.
At each nonzero horizontal frequency, two decreasing energies for the homogeneous evolution control $\eta$ and its second normal derivative.
These energies yield a velocity kernel bound in $B^{1/2}_{2,1}((-1,0))$, which controls boundary traces.
The projection of the bulk source and the harmonic extension of the boundary pressure give the forced estimate of Proposition~\ref{prop:global-linear-CL-proposition}.
At low horizontal frequencies, integration in time costs a factor $|\xi|^{-1}$, which is the reason for the assumption on $\mathcal B_{\rm low}(\eta_0)$.

To apply this estimate to the nonlinear solution, we use the divergence-free Piola velocity $\mathsf v$ from \eqref{eq:global-Piola-source}.
The bulk source is a product of a derivative of $\eta$ with $U$.
When $\eta$ is small in $X_m$, its $\mathcal Y$ norm is a small multiple of $\mathcal Y(U)$.
For the capillary source, the tangential pressure relation bounds the velocity generated by the nonlinear curvature remainder by a small multiple of $\mathcal Y(U)$.
The difference between $\mathsf v$ and $U$ satisfies the same type of product bound.
Absorbing these contributions gives Proposition~\ref{prop:global-source-closure}, and hence an $L_t^1\mathcal W_d$ bound for $U$ independent of the interval length.

\textbf{Anisotropic energy and the smallness argument.}
The anisotropic energy controls two vertical derivatives and $m$ horizontal derivatives of $\eta$, together with $m+1$ derivatives of $f$.
For the normal estimate, we work in physical coordinates with $h=y-\gamma^{-1}(\rho)$, whose pullback is $\eta$.
Differentiating its material equation gives an equation for $\Delta h$ in which the higher profile derivatives cancel.
By adding a surface quadratic correction, we cancel the leading boundary flux, and obtain the normal energy \eqref{eq:global-normal-energy}.
For the highest horizontal derivatives, a corrected pressure and flux preserve the divergence structure and control the differentiated transport terms.
Trace estimates and the capillary pressure condition give the surface regularity needed for the remaining boundary terms.

Combining these estimates with the relative energy gives an energy $\mathcal E$ equivalent to $e$ and satisfying the following inequality, where $e$ and $\mathcal D$ are defined in \eqref{eq:global-ed}:
\[
 \mathcal E'+\kappa\mathcal D
 \le Ce^{1/2}\mathcal D+C\mathcal W_d(U)e.
\]
The first term on the right can be absorbed into the dissipation when $e$ is small.
Then, the velocity integral controls the second term by Gronwall's inequality, yielding a uniform small bound for $e$ and a bound for the total dissipation.
By the anisotropic embedding, this also keeps $J_\eta$ close to $1$.
In addition, it preserves the invertibility of the density coordinate map and strict vertical monotonicity, and keeps the fluid depth positive.

\textbf{Continuation in the local solution class.}
The anisotropic bound does not control the full $H^m$ norm of the density required by the local continuation criterion.
We estimate this norm separately on each finite time interval, using the physical density transport equation.
The tame velocity estimate \eqref{eq:global-local-tame} is linear in the density norm, with a constant independent of that norm.
The required higher surface regularity follows from the energy and dissipation through Lemma~\ref{lem:global-surface-estimates}.
Gronwall's inequality then gives the finite-time density bound of Proposition~\ref{prop:global-strong-energy}.
Together with the small interface bound and positive depth, it permits continuation by Corollary~\ref{prop:global-local}.
The local theory also gives uniqueness and continuous dependence on every finite interval.

In this section, we develop the coordinate formulation and the forced linear estimate.
The nonlinear estimates and continuation argument are postponed to Section~\ref{sec:global-estimates}.

We write $\langle\cdot,\cdot\rangle$ for the $L^2$ inner product on the relevant domain (taken to be linear in the first slot and conjugate-linear in the second).
When no domain subscript is written, $\langle\cdot,\cdot\rangle$ denotes the $L^2(\R^d)$ inner product with respect to $\dd x$ (or its extension by Sobolev duality), while inner products in the vertical variable $z$ are taken over $I=(-1,0)$ with respect to $\dd z$.

\subsection{Density coordinates and the relative energy}
\label{sec:global-coordinates}

In the density coordinates defined by
\eqref{eq:global-density-map}, we choose the modified pressure
\begin{equation}\label{eq:global-Pi-definition}
\Pi=p\circ\Phi_\eta+\Gamma(z)+\gamma(z)\eta.
\end{equation}
The density transport equation, Darcy's law, and incompressibility give
\begin{equation}\label{eq:global-coordinate-identities}
 \begin{gathered}
 g\circ\Phi_\eta=\gamma(z)-\gamma(z+\eta),
 \qquad
 \eta_t=U_{d+1}-\nabla_x\eta\cdot U_x,\\
 \nabla_x\Pi=-U_x-\nabla_x\eta\,U_{d+1},
 \qquad
 \Pi_z=-J_\eta U_{d+1}+\gamma'(z)\eta,\\
 \Pi|_{z=0}=cf+\sigma H(f),
 \qquad \Pi_z|_{-1}=0,\qquad  \diver_x(J_\eta U_x) +\partial_z(U_{d+1}-\nabla_x\eta\cdot U_x)=0.
 \end{gathered}
\end{equation}
Define the Piola velocity by
\[
 \mathsf v =\bigl(J_\eta U_x,\, U_{d+1}-\nabla_x\eta\cdot U_x\bigr).
\]
Then $\diver\mathsf v=0$ and $\mathsf v_{d+1}=\eta_t$, with boundary conditions $\mathsf v_{d+1}|_{-1}=0$ and $\mathsf v_{d+1}|_{z=0}=f_t$.

For solutions of \eqref{e:IPM} with finite energy,
the pressure identities imply
\[
 \nabla\Pi\cdot\mathsf v =-J_\eta|U|^2+\gamma'(z)\eta\eta_t.
\]
Integrating over $\mathcal S$ and using the boundary conditions yields the relative energy identity
\begin{equation}\label{eq:global-relative-energy}
 \begin{aligned}
 &\frac{\dd}{\dd t}\Bigl\{
 \frac12\int_{\mathcal S}(-\gamma'(z))\eta^2\dd x\dd z +\frac c2\|f\|_{L_x^2}^2 +\sigma\int_{\R^d}
\bigl(\sqrt{1+|\nabla f|^2}-1\bigr)\dd x
 \Bigr\} +\int_{\mathcal S}J_\eta|U|^2\dd x\dd z=0.
 \end{aligned}
\end{equation}

To separate the linear evolution from the nonlinear sources, we write $\Pi=\Pi_L+\Pi_N$, where the linear
pressure $\Pi_L$ solves
\begin{equation}\label{eq:global-flat-pressure}
 -\Delta\Pi_L=-\partial_z(\gamma'(z)\eta), \qquad
 \Pi_L|_{z=0}=(c-\sigma\Delta_x)f, \qquad
 \partial_z\Pi_L|_{-1}=0.
\end{equation}
The corresponding linear velocity associated with
$Z=(f,\eta)$ is
\[
 U_Z=-\nabla\Pi_L+\gamma'(z)\eta e_{d+1}.
\]
Set
\begin{equation}\label{eq:global-nonlinear-source}
 \mathscr G_x=\nabla_x\eta\,U_{d+1}-\eta_zU_x,\qquad  \mathscr G_{d+1}=\nabla_x\eta\cdot U_x+\eta_zU_{d+1},\qquad
 \mathcal N(f)=H(f)+\Delta_xf.
\end{equation}
The coordinate identities then give
\begin{equation}\label{eq:global-Piola-source}
 \mathsf v=U_Z-\nabla\Pi_N-\mathscr G.
\end{equation}
Taking the divergence and using the pressure boundary conditions, we obtain the equation for $\Pi_N$:
\begin{equation}\label{eq:global-nonlinear-pressure}
 -\Delta\Pi_N=\diver\mathscr G, \qquad \Pi_N|_{z=0}=\sigma\mathcal N(f), \qquad \partial_z\Pi_N|_{-1}=0.
\end{equation}
At the bottom, $\nabla_x\eta=U_{d+1}=0$, so that $\mathscr G_{d+1}|_{-1}=0$.

Thus the equations on the fixed strip $\mathcal S$ read
\begin{equation}\label{eq:global-density-system}
 \begin{cases}
 \eta_t=\mathsf v_{d+1},\qquad \diver\mathsf v=0,\\
 \mathsf v+\nabla\Pi
   =\gamma'(z)\eta e_{d+1}-\mathscr G,\\
 \Pi|_{z=0}=(c-\sigma\Delta_x)f+\sigma\mathcal N(f),
 \qquad f=\eta|_{z=0},\\
 \eta|_{-1}=\mathsf v_{d+1}|_{-1}=0.
 \end{cases}
\end{equation}
Replacing $\mathscr G$ and $\sigma\mathcal N(f)$ by prescribed sources $G$ and $\vartheta$ gives the linear system \eqref{eq:global-forced-linear-system}.
We establish its velocity estimate below and apply it to the nonlinear sources in Section~\ref{sec:global-estimates}.

Recall from \eqref{eq:global-Xm} that
\[
 \|h\|_{X_m}^2
 =\|h\|_{L^2(\mathcal S)}^2
 +\|\partial_z^2h\|_{L^2(\mathcal S)}^2
 +\||D_x|^m h\|_{L^2(\mathcal S)}^2.
\]
In the rest of this paper, we use the anisotropic energy and dissipation
\begin{equation}\label{eq:global-ed}
e(t)=\|\eta(t)\|_{X_m}^2+\|f(t)\|_{H^{m+1}}^2,
 \qquad \mathcal D(t)=\|U(t)\|_{X_m}^2.
\end{equation}

The choice of $X_m$ comes from the bulk equation.
For the homogeneous linearized system, eliminating the pressure gives
$\Delta\eta_t=\gamma'(z)\Delta_x\eta$.
Testing against $\Delta\eta/(-\gamma')$ gives an energy for
$\Delta\eta$ with a correction at the free surface.
Since $\eta_{zz}=\Delta\eta-\Delta_x\eta$, this energy controls
$\eta_{zz}$ once the horizontal derivatives are bounded.
The nonlinear counterpart is \eqref{eq:global-normal-coercivity}.
We take $m\ge 2d+5$ so that $\|\eta\|_{X_m}$ also controls
multiplication by $\nabla_{x,z}\eta$ in the velocity estimate;
see Lemma~\ref{lem:global-source-product}.

The following embedding controls the Lipschitz and Fourier--Wiener norms by $X_m$.

\begin{lemma}[Anisotropic embeddings]
\label{lem:global-embeddings}
Let $E_z$ be a fixed extension operator from $I$ to $\R$, bounded on both $L_z^2$ and $H_z^2$ and commuting with horizontal derivatives.
For every integer $m>2d$ and every $h\in X_m$,
\begin{equation}\label{eq:global-Fourier-embeddings}
 \mathcal W_d(h) +\int_{\R^d}\int_\R   (|\xi|+|\zeta|) |\widehat{E_zh}(\xi,\zeta)|\dd\zeta\dd\xi \le C\|h\|_{X_m}.
\end{equation}
In particular,
\[
 \|h\|_{L^\infty(\mathcal S)} +\|\nabla_{x,z}h\|_{L^\infty(\mathcal S)} \le C\|h\|_{X_m}.
\]
\end{lemma}

\begin{proof}
Interpolation on $I$ controls $\|\partial_z h\|_{L^2}$ by $\|h\|_{L^2}+\|\partial_z^2h\|_{L^2}$.
The extension therefore satisfies
\[
 \int_{\R^d}\int_\R (1+|\xi|^{2m}+\zeta^4)
 |\widehat{E_zh}(\xi,\zeta)|^2\dd\zeta\dd\xi \le C\|h\|_{X_m}^2.
\]
By the definition of $\mathcal W_d$ (see \eqref{eq:global-Wiener}) and the Cauchy--Schwarz inequality, it suffices to check the integrability of the reciprocal weights.
Integration in $\zeta$ gives
\[
 \int_\R\frac{(1+|\xi|)^2}{1+|\xi|^{2m}+\zeta^4}\dd\zeta
\le C(1+|\xi|)^2(1+|\xi|^{2m})^{-3/4},\qquad 
 \int_\R\frac{\zeta^2}{1+|\xi|^{2m}+\zeta^4}\dd\zeta
\le C(1+|\xi|^{2m})^{-1/4}.
\]
Both right-hand sides are integrable in $\xi\in\R^d$
when $m>2d$ and $d\ge1$.
This proves the estimate \eqref{eq:global-Fourier-embeddings}.
The $L^\infty$ bounds follow by Fourier inversion.
\end{proof}

\subsection{Initial coordinates and local well-posedness}

For local Sobolev estimates, we use the graph-flattening map
\eqref{eq:explicit-flattening} with $\ell_0=1$ and $b_0=0$:
\begin{equation}\label{eq:global-flattening}
 \Phi_f(t,x,z)=(x,\varrho_f(t,x,z)),\qquad
 \varrho_f=z+\mathscr E f,\qquad
 \mathscr E h=(1+z)e^{\delta z\langle D_x\rangle}h.
\end{equation}
We choose a fixed $\delta>0$ as in Lemma~\ref{app:ln-flattening}, so that $\Phi_f$ is a diffeomorphism on the bounded set of interfaces under consideration.
In these coordinates, the density satisfies \eqref{eq:ALE}, and the kinematic and bottom conditions give $V_{f,d+1}=0$ on both fixed boundaries.

\begin{lemma}[Regularity of the initial density coordinate map]\label{lem:global-initial-chart}
For compatible initial data in Definition~\ref{def:global-data} with $\|f_0\|_{H^{m+1}}+\|\widetilde g_0\|_{H^m}$ sufficiently small, the map
\[
 (x,y)\longmapsto\bigl(x,\gamma^{-1}(\gamma(y)+g_0(x,y))\bigr)
\]
is a $C^1$ diffeomorphism from $\Om_{f_0}$ onto $\mathcal S$.
Its inverse $(x,z)\mapsto(x,z+\eta_0(x,z))$ satisfies
\begin{equation}\label{eq:global-initial-coordinate}
 \|\eta_0\|_{H^m(\mathcal S)} \le C\bigl(\|f_0\|_{H^{m+1}}+\|\widetilde g_0\|_{H^m}\bigr),
 \qquad \|\eta_{0,z}\|_{L^\infty}<\frac12.
\end{equation}
The map $(f_0,\widetilde g_0)\mapsto\eta_0$ is continuous into $H^m(\mathcal S)$.
\end{lemma}

\begin{proof}
Since $\rho_0(x,\cdot)$ is strictly decreasing, its boundary values show that $y\mapsto\gamma^{-1}(\rho_0(x,y))$ maps $[-1,f_0(x)]$ bijectively onto $[-1,0]$, with positive derivative.
In the regular graph coordinates this map becomes
\[
 (x,z)\longmapsto \left(x,\gamma^{-1} \bigl(\gamma(\varrho_{f_0}(x,z)) +\widetilde g_0(x,z)\bigr)\right).
\]
Since $y-\gamma^{-1}(\gamma(y)+q)$ vanishes at $q=0$, Sobolev composition estimates show that this map differs from the identity in $H^m$ by at most $C(\|f_0\|_{H^{m+1}}+\|\widetilde g_0\|_{H^m})$.
The Sobolev estimates for inversion of a diffeomorphism, followed by composition with $\Phi_{f_0}$, give \eqref{eq:global-initial-coordinate} and continuous dependence.
\end{proof}

Theorem~\ref{thm:local-target} and Corollary~\ref{cor:local-continuation} give the following result in the present setting.
\begin{corollary}[Local well-posedness and continuation criterion]\label{prop:global-local}
Let $m\ge2d+5$ be an integer.
For compatible initial data as in Definition~\ref{def:global-data}, \eqref{e:IPM} has a unique local strong solution in \eqref{eq:global-local-class}, depending continuously on the data in those spaces.
The lifespan is uniform on bounded subsets of $H^{m+1}(\R^d)\times H^m(\mathcal S)$ with a uniform positive depth bound and a fixed admissible vertical range for $\gamma$.
In particular, a solution on $[0,T_*)$, with $T_*<\infty$, can be extended beyond $T_*$ if
\[
 \sup_{t<T_*}\|f(t)\|_{H^{m+1}}\quad\text{remains sufficiently small},
 \qquad \sup_{t<T_*}\|\widetilde g(t)\|_{H^m(\mathcal S)}<\infty.
\]
\end{corollary}

\subsection{The linear velocity estimate with prescribed sources}\label{sec:global-linear-section}

Next, we estimate the time integral of the velocity in terms of the initial datum $\eta_0$ and prescribed bulk and boundary sources.
For each nonzero horizontal frequency $\xi$, we write $r=|\xi|$.
The coefficient $-\gamma'(z)$ is time independent and satisfies the bounds in \eqref{eq:global-profile}.

To apply the linear estimate in Section~\ref{sec:global-estimates},
we need a bound for the time integral of the Fourier--Wiener
norm of the velocity. At each horizontal frequency,
$B^{1/2}_{2,1}(I)$ controls the Fourier $L^1$ norm in the
vertical variable after extension, as well as the boundary traces.
Moreover, incompressibility gives
$\|\partial_z\widehat v_{d+1}\|_2\le r\|\widehat v_x\|_2$.
Interpolation therefore bounds
$\|\widehat v_{d+1}\|_{B^{1/2}_{2,1}(I)}$ by
$C\langle r\rangle^{1/2}\|\widehat v\|_2$.
This leads us to set
\begin{equation}\label{eq:global-linear-space-Y}
 Y_r(G)=\|G\|_{B^{1/2}_{2,1}(I)}+\langle \xi \rangle^{1/2}\|G\|_{L^2},
 \qquad \mathcal Y(F)=\int_{\mathbb R^d}(1+|\xi|)Y_r(\widehat F(\xi,\cdot))\,\dd\xi.
\end{equation}
After a fixed vertical extension, the Fourier $L^1_\zeta$ norm
of each vertical dyadic block at frequency $2^j$ is bounded
by $C2^{j/2}$ times its $L^2_z$ norm.
Summing over $j$ and integrating in $\xi$ with weight $1+|\xi|$
gives
\begin{equation}\label{eq:global-linear-Y-Wiener}
 \mathcal W_d(F)\le C\mathcal Y(F).
\end{equation}
Moreover, the Cauchy--Schwarz inequality and the embedding $H^1(I)\hookrightarrow B^{1/2}_{2,1}(I)$ give $\mathcal Y(F)\le C\|F\|_{H^{d+3}}$.

For a boundary pressure $\vartheta$, we define $q_\vartheta$ and $U_\vartheta$ by
\begin{equation} \label{Defqvartheta}
 \Delta q_\vartheta=0,\qquad q_\vartheta|_{z=0}=\vartheta, \qquad \partial_zq_\vartheta|_{-1}=0, \qquad U_\vartheta=-\nabla q_\vartheta.
\end{equation}

\begin{proposition}[Forced linear velocity estimate]
\label{prop:global-linear-CL-proposition}
Let $\gamma=\gamma(z)$ be a time-independent profile with $0<c_*\le-\gamma'(z)\le C_*$, let $c>0$, and fix $\sigma\ge0$.
Suppose that $(f,\eta,v,\Pi)$ solves
\begin{equation}\label{eq:global-forced-linear-system}
 \begin{gathered}
 \eta_t=v_{d+1},\qquad \operatorname{div}v=0,\qquad
 v+\nabla\Pi=\gamma'(z)\eta e_{d+1}-G
          \quad\text{in }\mathcal S,\\
 f=\eta|_{z=0},\qquad \Pi|_{z=0}=(c-\sigma\Delta_x)f+\vartheta, \qquad \eta|_{-1}=v_{d+1}|_{-1}=0.
 \end{gathered}
\end{equation}
Then for $m\ge2d+5$ and every $T>0$,
\begin{equation}\label{eq:global-linear-CL-estimate}
 \int_0^T\mathcal Y(v(t))\,\dd t
 \le C\Big( \|\eta_0\|_{H^m(\mathcal S)}+\mathcal B_{\rm low}(\eta_0) +\int_0^T\bigl(\mathcal Y(G(t)) +\mathcal Y(U_\vartheta(t))\bigr)\,\dd t\Big),
\end{equation}
where $\eta_0=\eta(0)$ and $\mathcal B_{\rm low}(\eta_0)$ is defined as in \eqref{eq:global-low-norm}.
The constant $C$ is independent of $T$.
\end{proposition}

In Lemma \ref{lem:direct-refined-energy} below, we first derive a homogeneous velocity estimate from two nonincreasing energies.
Applying this estimate to cosine blocks and summing gives
the kernel bound for initial fields with zero top normal trace in Lemma \ref{lem:direct-velocity-kernel}.
For a general initial field, we separate a harmonic field with the same top normal trace and estimate its contribution by Duhamel's formula.
Finally, we apply the kernel bound to the forced system after projecting the bulk source onto the divergence-free subspace and extending the boundary pressure harmonically.
This gives the proof of Proposition \ref{prop:global-linear-CL-proposition}.

Fix $\xi\ne0$, set $r=|\xi|$, and suppress the Fourier hats.
Let $\mathcal H_\xi$ be the closed subspace of
$L^2(I;\mathbb C^{d+1})$ consisting of fields $G=(G_x,G_{d+1})$
such that
\[
 i\xi\cdot G_x+\partial_zG_{d+1}=0,
 \qquad G_{d+1}|_{z=-1}=0.
\]
Integrating the divergence constraint gives $G_{d+1}|_{z=0}=-i\xi\cdot\int_I G_x\,\dd z$, and the Cauchy--Schwarz inequality yields
\[
 |G_{d+1}(0)| \le |\xi|\int_I|G_x(z)|\,\dd z \le r\|G_x\|_{L^2(I)}.
\]

In the homogeneous Darcy equation, the bulk forcing is $\gamma'(z)\eta e_{d+1}$ and the top pressure is $(c+\sigma r^2)f$. 
To obtain the bulk and surface energy terms, we use a divergence-free test field with vertical component $\eta$, so that its top normal trace is $f$.
Choosing the horizontal component parallel to $\xi$ and imposing the divergence constraint gives
\[
 D_\xi\eta=\left(\frac{i\xi}{r^2}\partial_z\eta,\eta\right)
 \in\mathcal H_\xi.
\]

Testing Darcy's law against an arbitrary $\psi\in\mathcal H_\xi$
and integrating by parts now gives $U=-L_rD_\xi\eta$, where
\begin{equation}\label{eq:direct-velocity-operator}
 \langle L_rG,\psi\rangle_{L^2(I)}
 =\int_I(-\gamma'(z))G_{d+1}\overline{\psi_{d+1}}\,\dd z
  +(c+\sigma r^2)
    \left.(G_{d+1}\overline{\psi_{d+1}})\right|_{z=0}.
\end{equation}
By the trace bound above, this form defines a bounded,
nonnegative self-adjoint operator on $\mathcal H_\xi$.
In particular, taking $\psi=D_\xi\eta$ gives
$-\operatorname{Re}\langle U,D_\xi\eta\rangle_{L^2(I)}=2E(t)$,
where $E(t)$ is the basic energy in the following lemma.
\begin{lemma}[A direct homogeneous estimate]\label{lem:direct-refined-energy}
At fixed $\xi\ne0$, consider
\[
 \eta_t=U_{d+1},\qquad U=-\nabla\Pi+\gamma'(z)\eta e_{d+1},\qquad
 \operatorname{div}U=0,
\]
with $\eta|_{z=-1}=U_{d+1}|_{z=-1}=0$, $f=\eta|_{z=0}$, and $\Pi|_{z=0}=(c+\sigma r^2)f$.
Let
\[
 E(t)=\frac12\int_I(-\gamma'(z))|\eta|^2\dd z+\frac{c+\sigma r^2}{2}|f|^2,
 \qquad Q(t)=\frac12\int_I\frac{|\eta_{zz}-r^2\eta|^2}{-\gamma'(z)}\dd z +\frac{|\eta_z(0)|^2}{2(c+\sigma r^2)},
\]
$E_0=E(0)$ and $H_0=Q(0)+(1+r^4)E(0)$.
Then
\begin{equation}\label{eq:direct-all-refined}
 \int_0^\infty Y_r(U(t))\,\dd t
 \le Cr^{-1}\left((1+r)^{1/2}E_0^{1/4}H_0^{1/4} +E_0^{1/8}H_0^{3/8}\right).
\end{equation}
\end{lemma}

\begin{proof}
By Darcy's law and incompressibility of the fluid,
\[
 \partial_zU_{d+1}=-r^2\Pi,\qquad (\partial_z^2-r^2)U_{d+1}=-\gamma'(z)r^2\eta,
 \qquad \partial_t\eta_z(0)=-r^2(c+\sigma r^2)f.
\]
Testing Darcy's law with $U$ and differentiating $Q$ yield the equations
\[
 E'=-\|U\|_{L^2}^2,\qquad Q'=-r^2\|\eta_z\|_{L^2}^2-r^4\|\eta\|_{L^2}^2.
\]
The top term in $Q$ cancels the boundary term from integration by parts.
Thus $Q+(1+r^4)E\le H_0$, which controls $\|\eta\|_{H^2}^2$.
If $E_0=0$, the energy identity gives the zero solution.
We may therefore assume $E_0>0$ and $H_0>0$.
Testing instead with $D_\xi\eta$ gives
\[
 2E=-\operatorname{Re}\langle U,D_\xi\eta\rangle,
 \qquad \|D_\xi\eta\|_{L^2}\le Cr^{-1}E^{1/4}H_0^{1/4}.
\]
Interpolation controls $\eta_z$, and $H_0\ge(1+r^4)E$ absorbs the remaining $\|\eta\|_{L^2}$ term.
Together with $E'=-\|U\|_{L^2}^2$, these bounds give
\[
 E'+\kappa r^2H_0^{-1/2}E^{3/2}\le0.
\]
Integrating this differential inequality yields
\[
 E(t)\le CE_0\left(1+\kappa r^2\sqrt{E_0/H_0}\,t\right)^{-2}.
\]
Differentiating Darcy's law in time and testing with $U$ gives
\[
 \frac12\frac{\dd}{\dd t}\|U\|_{L^2}^2+\int_I(-\gamma'(z))|U_{d+1}|^2\dd z+(c+\sigma r^2)|U_{d+1}(0)|^2=0.
\]
Hence $\|U(t)\|_{L^2}^2\le2t^{-1}E(t/2)$.
Integrating these scalar bounds gives
\[
 \int_0^\infty \|U\|_{L^2}\,\dd t\le Cr^{-1}E_0^{1/4}H_0^{1/4},
 \qquad \int_0^\infty \|U\|_{L^2}^{1/2}E^{1/4}\,\dd t \le Cr^{-3/2}E_0^{1/8}H_0^{3/8}.
\]
Finally, $\partial_zU_x=i\xi(U_{d+1}-\gamma'(z)\eta)$ and $\partial_zU_{d+1}=-i\xi\cdot U_x$ imply
\[
 \|U\|_{H^1}\le C\bigl((1+r)\|U\|_{L^2}+rE^{1/2}\bigr).
\]
Using the interpolation inequality $\|U\|_{B^{1/2}_{2,1}}\le C\|U\|_{L^2}^{1/2}\|U\|_{H^1}^{1/2}$, we prove the estimate \eqref{eq:direct-all-refined}.
\end{proof}

\begin{lemma}[A velocity-kernel bound]\label{lem:direct-velocity-kernel}
For every $\xi\ne0$ and every $G\in\mathcal H_\xi$ with $Y_r(G)<\infty$, the following estimate holds with a constant independent of $r=|\xi|$:
\begin{equation}\label{eq:direct-all-forced-kernel}
 \int_0^\infty Y_r(e^{-tL_r}L_rG)\,\dd t\le CY_r(G).
\end{equation}
\end{lemma}

\begin{proof}
Fix $\xi\ne0$. For the homogeneous solution with
$\eta_0=G_{d+1}$ and $f_0=G_{d+1}(0)$, weak Darcy's law gives
\[
 U(0)=-L_rG,\qquad U_t+L_rU=0.
\]
Thus $U(t)=-e^{-tL_r}L_rG$, and Lemma~\ref{lem:direct-refined-energy} gives an estimate for this kernel.

The bottom condition $G_{d+1}(-1)=0$ and incompressibility $\partial_zG_{d+1}=-i\xi\cdot G_x$ give
\[
 \|G_{d+1}\|_{L^2} \le \min(1,Cr)\|G\|_{L^2}
 \le C\frac r{\langle r\rangle}\|G\|_{L^2}, \qquad G\in\mathcal H_\xi.
\]

We first prove the estimate on the subspace of fields with zero normal trace at the top,
\[
 \mathcal H_\xi^b =\{G\in\mathcal H_\xi:G_{d+1}(0)=0\}.
\]
For $G\in\mathcal H_\xi^b$, we set
\[
 \alpha=\frac{\xi}{r}\cdot G_x,
 \qquad \widetilde G=\left(\frac{\xi}{r}\alpha,G_{d+1}\right).
\]
Since $\xi\cdot\widetilde G_x=\xi\cdot G_x$, the field
$\widetilde G$ remains divergence free and belongs to
$\mathcal H_\xi^b$.
The definition of $L_r$ and boundedness of the horizontal
orthogonal projection give
\[
 L_r\widetilde G=L_rG,\qquad Y_r(\widetilde G)\le CY_r(G).
\]
For this part of the proof, it therefore suffices to consider $G_x=(\xi/r)\alpha$.
Incompressibility gives
\[
 G_{d+1}(z)=-ir\int_{-1}^z\alpha(s)\,\dd s,
 \qquad \int_I\alpha(s)\,\dd s=0.
\]
Since $\alpha$ has zero mean on $I=(-1,0)$, its expansion as the sum of cosine series has no constant term:
\[
 \alpha(z)=\sum_{k=1}^\infty a_k\cos\bigl(k\pi(z+1)\bigr),\qquad
 a_k=2\int_{-1}^0\alpha(z)\cos\bigl(k\pi(z+1)\bigr)\,\dd z,
\]
with convergence in $L^2(I)$.
  For each dyadic $N\geq 1$, we define the dyadic cosine block
\[
 \alpha_N(z)=\sum_{k=N}^{2N-1}a_k\cos\bigl(k\pi(z+1)\bigr),
 \qquad \alpha=\sum_{N\ge1}\alpha_N.
\]
Integrating each block gives the corresponding vector field
\[
 G_N(z)=\left( \frac{\xi}{r}\alpha_N(z), -ir\int_{-1}^z\alpha_N(s)\,\dd s \right).
\]
Direct computation shows that the vertical component is a sum of sine modes with the same indices $N\le k<2N$.
Thus $G_N$ is divergence free, its normal component vanishes at both endpoints, and $G=\sum_{N\ge1}G_N$ in $L^2(I)$.
This decomposition applies to the initial field, and each $G_N$ evolves under the full operator $L_r$.

For $N\ge\langle r\rangle$, the initial energies in
Lemma~\ref{lem:direct-refined-energy} satisfy
\[
 E_0^{1/2}\le CrN^{-1}\|\alpha_N\|_{L^2},
 \qquad
 H_0^{1/2}\le CrN\|\alpha_N\|_{L^2}.
\]
The bulk bounds follow from the normal frequency localization.
As for the boundary term, we use $\eta_{0,z}(0)=-ir\alpha_N(0)$ and $|\alpha_N(0)|\le CN^{1/2}\|\alpha_N\|_{L^2}$.
Lemma~\ref{lem:direct-refined-energy} therefore gives
\begin{equation}\label{eq:direct-high-normal-block}
 \int_0^\infty Y_r(e^{-tL_r}L_rG_N)\,\dd t
 \le CN^{1/2}\|\alpha_N\|_{L^2}.
\end{equation}

We further group all blocks with $N<\langle r\rangle$ into one field $G_{\mathrm{lo}}$, and let $\eta_0=(G_{\mathrm{lo}})_{d+1}$.
It vanishes at $z = 0$, and Bernstein's inequalities give
\[
 E_0^{1/2}\le C\|\eta_0\|_{L^2}, \qquad
 H_0^{1/2}\le C\langle r\rangle^2\|\eta_0\|_{L^2}.
\]
Indeed,
$\|\eta_{0,zz}\|_{L^2}\le C\langle r\rangle^2\|\eta_0\|_{L^2}$
and
$|\eta_{0,z}(0)|\le C\langle r\rangle^{3/2}\|\eta_0\|_{L^2}$.
Applying Lemma~\ref{lem:direct-refined-energy} and the vertical-component bound above, we obtain
\[
 \begin{aligned}
 \int_0^\infty
 Y_r(e^{-tL_r}L_rG_{\mathrm{lo}})\,\dd t
 \le Cr^{-1}\langle r\rangle^{3/2}\|\eta_0\|_{L^2}\le C\langle r\rangle^{1/2}\|G_{\mathrm{lo}}\|_{L^2}.
 \end{aligned}
\]

Note that the right-hand side of \eqref{eq:direct-high-normal-block} can be bounded by
\[
 \sum_{N\ge1}N^{1/2}\|\alpha_N\|_{L^2(I)}
 \le C\|\alpha\|_{B^{1/2}_{2,1}(I)}.
\]
Since $\|G_{\mathrm{lo}}\|_{L^2}\le C\|G\|_{L^2}$,
summing the high- and low-frequency estimates yields,
first for finite cosine sums,
\begin{equation}\label{eq:direct-zero-top-kernel}
 \int_0^\infty Y_r(e^{-tL_r}L_rG)\,\dd t
 \le CY_r(G),\qquad G\in\mathcal H_\xi^b.
\end{equation}

For general $\alpha$, one can choose mean-zero cosine polynomials $\alpha^{(M)}\to\alpha$ in $B^{1/2}_{2,1}(I)$ and reconstruct $G^{(M)}$ by the preceding integration formula.
These fields have zero normal trace at both endpoints.
Since integration maps $L^2(I)$ continuously into $H^1(I)$, we have $Y_r(G^{(M)}-G)\to0$ for fixed $r$.
Boundedness of $e^{-tL_r}L_r$ on $\mathcal H_\xi$ gives
convergence of the kernel outputs in $L^2(I)$ at every time.
Using lower semicontinuity and Fatou's lemma, we can extend \eqref{eq:direct-zero-top-kernel} to every
$G\in\mathcal H_\xi^b$ with finite $Y_r(G)$.

We now return to a general $G\in\mathcal H_\xi$ with $Y_r(G)<\infty$.
To separate its top normal trace, define
\[
 h_r(z)=\frac{\sinh(r(z+1))}{\sinh r},
 \qquad e_r=\frac{D_\xi h_r}{\sqrt{\coth r/r}}.
\]
Integration by parts gives, for every $G\in\mathcal H_\xi$,
\[
 \|e_r\|_{L^2}=1,\qquad \langle G,e_r\rangle(e_r)_{d+1}(0)=G_{d+1}(0).
\]
Thus $\mathcal H_\xi^b$ is the orthogonal complement of $e_r$ in $\mathcal H_\xi$, and
\[
 G=G^b+\langle G,e_r\rangle e_r, \qquad G^b\in\mathcal H_\xi^b.
\]
The explicit formula and interpolation give
\[
 |(e_r)_{d+1}0)|^2=r\tanh r,\qquad \|(e_r)_{d+1}\|_{L^2}\le C\frac r{\langle r\rangle}, \qquad Y_r(e_r)\le C\langle r\rangle^{1/2}.
\]
Thus
\[
 Y_r(G^b)+\langle r\rangle^{1/2}|\langle G,e_r\rangle| \le CY_r(G).
\]
It remains to prove
\[
 \int_0^\infty Y_r(e^{-tL_r}L_re_r)\,\dd t \le C\langle r\rangle^{1/2}.
\]

We decompose $L_re_r$ into its components parallel and orthogonal to $e_r$:
\[
 L_re_r=\omega_re_r+\beta_r,\qquad \omega_r=\langle L_re_r,e_r\rangle.
\]
Since $\langle\beta_r,e_r\rangle=0$, the trace identity above gives $\beta_r\in\mathcal H_\xi^b$.
The definition of $L_r$ implies
\[
 \omega_r\ge(c+\sigma r^2)r\tanh r \ge C^{-1}\frac{r^2}{\langle r\rangle}.
\]
For $\psi\in\mathcal H_\xi^b$, the boundary term vanishes:
\[
 \langle\beta_r,\psi\rangle =\int_I(-\gamma')(e_r)_y\overline{\psi_y}\,\dd z.
\]
Applying the vertical-component bound to $e_r$ and $\psi$ therefore gives
\[
 \|\beta_r\|_{L^2}\le C\frac{r^2}{\langle r\rangle^2}.
\]
The same identity shows that $\beta_r$ is the orthogonal projection of $(-\gamma')(e_r)_{d+1}e_{d+1}$ onto $\mathcal H_\xi^b$.
Its divergence and curl therefore satisfy
\[
 \partial_z(\beta_r)_{d+1}=-i\xi\cdot(\beta_r)_x,\qquad \partial_z(\beta_r)_x =i\xi\bigl((\beta_r)_{d+1}+ \gamma'(e_r)_{d+1}\bigr).
\]
Hence,
\[
 \|\beta_r\|_{H^1} \le C\bigl((1+r)\|\beta_r\|_{L^2}+r\|(e_r)_{d+1}\|_{L^2}\bigr) \le C\frac{r^2}{\langle r\rangle}.
\]
Interpolation yields
\[
 Y_r(\beta_r)\le C\frac{r^2}{\langle r\rangle^{3/2}}, \qquad \omega_r^{-1}Y_r(\beta_r) \le C\langle r\rangle^{-1/2}.
\]

For $k_r(t)=e^{-tL_r}L_re_r$, the decomposition of $L_re_r$ shows that $k_r$ satisfies the equation
\[
 k_r'+\omega_rk_r=-e^{-tL_r}L_r\beta_r, \qquad k_r(0)=L_re_r.
\]
Since $\beta_r\in\mathcal H_\xi^b$, the right-hand side is integrable in $Y_r$ by \eqref{eq:direct-zero-top-kernel}.
Duhamel's formula therefore gives
\[
 \begin{aligned}
 \int_0^\infty Y_r(k_r(t))\,\dd t
 &\le \omega_r^{-1} \left( Y_r(L_re_r) +\int_0^\infty Y_r(e^{-tL_r}L_r\beta_r)\,\dd t \right)\\
 &\le Y_r(e_r)+C\omega_r^{-1}Y_r(\beta_r)\le C\langle r\rangle^{1/2}.
 \end{aligned}
\]
Applying \eqref{eq:direct-zero-top-kernel} to $G^b$ and the last bound to the component along $e_r$, we obtain
\[
 \int_0^\infty Y_r(e^{-tL_r}L_rG)\,\dd t
 \le CY_r(G^b)+C\langle r\rangle^{1/2}|\langle G,e_r\rangle|
 \le CY_r(G).
\]
This proves \eqref{eq:direct-all-forced-kernel} on $\mathcal H_\xi$.
\end{proof}

We now prove Proposition~\ref{prop:global-linear-CL-proposition} using the two lemmas above.

\begin{proof}[Proof of Proposition~\ref{prop:global-linear-CL-proposition}]
Fix $\xi\ne0$, set $r=|\xi|$, and let $\mathbb P_r$ be the $L^2(I)$ orthogonal projection onto $\mathcal H_\xi$.
We first prove that the bulk source satisfies
\begin{equation}\label{eq:direct-application-projection}
 Y_r(\mathbb P_r\widehat G)\le CY_r(\widehat G),
\end{equation}
with $C$ independent of $r$.
For smooth $\widehat G$, let $\widehat\pi$ solve
\[
 (\partial_z^2-r^2)\widehat\pi
 =i\xi\cdot\widehat G_x+\partial_z\widehat G_{d+1},\qquad
 \widehat\pi|_{z=0}=0,\qquad
 \partial_z\widehat\pi|_{-1}=\widehat G_{d+1}|_{-1}.
\]
Then
\[
 \widehat V=\widehat G-(i\xi\widehat\pi,\partial_z\widehat\pi)
\]
belongs to $\mathcal H_\xi$, while the gradient term is orthogonal to $\mathcal H_\xi$ by integration by parts.
Thus $\widehat V=\mathbb P_r\widehat G$, and orthogonality gives
\[
 \|\widehat V\|_{L^2}
 +\|i\xi\widehat\pi\|_{L^2}
 +\|\partial_z\widehat\pi\|_{L^2}
 \le C\|\widehat G\|_{L^2}.
\]
Since $\partial_z\widehat V_{d+1}=-i\xi\cdot\widehat V_x$ and
$\partial_z(i\xi\widehat\pi)=i\xi\partial_z\widehat\pi$, we also have
\[
 \|\widehat V_{d+1}\|_{H^1(I)}
 +\|i\xi\widehat\pi\|_{H^1(I)}
 \le C\langle r\rangle\|\widehat G\|_{L^2}.
\]
The interpolation inequality
$\|w\|_{B^{1/2}_{2,1}(I)}\le C\|w\|_{L^2(I)}^{1/2}\|w\|_{H^1(I)}^{1/2}$ therefore yields
\[
 \|\widehat V_{d+1}\|_{B^{1/2}_{2,1}(I)}
 +\|i\xi\widehat\pi\|_{B^{1/2}_{2,1}(I)}
 \le C\langle r\rangle^{1/2}\|\widehat G\|_{L^2}.
\]
Using $\widehat V_x=\widehat G_x-i\xi\widehat\pi$ and the $L^2$ projection bound, we obtain \eqref{eq:direct-application-projection}.
The estimate can be extended to every $\widehat G\in B^{1/2}_{2,1}(I)$ by density.

We now take the horizontal Fourier transform of \eqref{eq:global-forced-linear-system}.
Testing the forced Darcy equation against
$\psi\in\mathcal H_\xi$, and using the top pressure condition, we get
\[
 \langle \widehat v,\psi\rangle
 =-\langle L_rD_\xi\widehat\eta,\psi\rangle -\langle \widehat G,\psi\rangle +\langle \widehat U_\vartheta,\psi\rangle.
\]
Thus the velocity satisfies
\[
 \widehat v=-\mathbb P_r\widehat G+\widehat U_\vartheta-L_rD_\xi\widehat\eta.
\]

The operator $L_r$ is time independent and depends only on the vertical component of its argument.
Since $\widehat\eta_t=\widehat v_{d+1}$, we have $L_r\partial_t D_\xi\widehat\eta=L_r\widehat v$, and hence
\[
 (\partial_t+L_r)L_rD_\xi\widehat\eta
 =L_r(-\mathbb P_r\widehat G+\widehat U_\vartheta).
\]
Duhamel's formula gives
\[
 L_rD_\xi\widehat\eta(t)=e^{-tL_r}L_rD_\xi\widehat\eta(0)+\int_0^t e^{-(t-\tau)L_r}L_r
 \bigl(-\mathbb P_r\widehat G(\tau)+\widehat U_\vartheta(\tau)\bigr)\,\dd\tau.
\]
Applying \eqref{eq:direct-all-forced-kernel} to the initial term and the two source terms, then using Tonelli's theorem, yields
\[
 \int_0^T Y_r(L_rD_\xi\widehat\eta(t))\,\dd t\le CY_r(D_\xi\widehat\eta(0))+C\int_0^T\bigl(Y_r(\mathbb P_r\widehat G(t))+Y_r(\widehat U_\vartheta(t))\bigr)\,\dd t.
\]
The formula for $\widehat v$ and the projection bound
\eqref{eq:direct-application-projection} therefore imply
\[
 \int_0^T Y_r(\widehat v(t))\,\dd t \le CY_r(D_\xi\widehat\eta(0)) +C\int_0^T\bigl(Y_r(\widehat G(t))+Y_r(\widehat U_\vartheta(t))\bigr)\,\dd t.
\]
Multiplying the above inequality by $1+r$ and integrating in $\xi$ gives
\begin{equation}\label{eq:global-linear-velocity-data}
 \begin{aligned}
 \int_0^T\mathcal Y(v(t))\,\dd t
 \le C\int_{\R^d}(1+r)Y_r(D_\xi\widehat\eta_0)\,\dd\xi +C\int_0^T \bigl(\mathcal Y(G(t))+\mathcal Y(U_\vartheta(t))\bigr)\,\dd t.
 \end{aligned}
\end{equation}

It remains to estimate the initial field
$D_\xi\widehat\eta(0)=D_\xi\widehat\eta_0$.
For $0<r<1$,
\[
 (1+r)Y_r(D_\xi\widehat\eta_0) \le Cr^{-1}\|\widehat\eta_0\|_{H^2(I)}.
\]
For $r\ge1$, the embedding
$H^1(I)\hookrightarrow B^{1/2}_{2,1}(I)$ gives
\[
 (1+r)Y_r(D_\xi\widehat\eta_0) \le C\left(
 \|\widehat\eta_0\|_{H^2(I)} +r\|\widehat\eta_0\|_{H^1(I)} +r^{3/2}\|\widehat\eta_0\|_{L^2(I)}
 \right).
\]
The integral over $r<1$ is controlled by $\mathcal B_{\rm low}(\eta_0)$.
For $r\ge1$, the Cauchy--Schwarz inequality in $\xi$ gives
\begin{equation}\label{eq:global-linear-initial-bound}
 \int_{\R^d}(1+r)Y_r(D_\xi\widehat\eta_0)\,\dd\xi
 \le C\left(
 \mathcal B_{\rm low}(\eta_0)
 +\|\eta_0\|_{H^m(\mathcal S)}
 \right),
 \qquad m>d/2+2.
\end{equation}
Combining this bound with
\eqref{eq:global-linear-velocity-data}
proves \eqref{eq:global-linear-CL-estimate}.
\end{proof}

\section{Global well-posedness}\label{sec:global-estimates}

To finish the proof of Theorem~\ref{thm:global}, we use the linear estimate of Section \ref{sec:global-linear-estimates} to control the time integral of the nonlinear velocity.
We first combine this bound with anisotropic energy and finite-time density estimates, then apply the local continuation criterion.

\subsection{The nonlinear velocity estimate}\label{sec:global-source}
The goal of this subsection is to establish the following estimate.

\begin{proposition}[Uniform velocity integral]\label{prop:global-source-closure}
Let $(f,\eta,U)$ be a local strong solution with a regular density coordinate map, and set
\[
 e_T=\sup_{0\le t\le T} \bigl(\|\eta(t)\|_{X_m}+\|f(t)\|_{H^{m+1}}\bigr).
\]
There exists $\delta>0$, independent of $T$, such that $e_T\le\delta$ implies
\begin{equation}\label{eq:global-source-closed}
 \int_0^T\bigl(\mathcal Y(U(t))+\mathcal Y(\mathsf v(t))\bigr)\,\dd t \le C\bigl(\|\eta(0)\|_{H^m} +\mathcal B_{\rm low}(\eta(0))\bigr).
\end{equation}
For the initial data of Definition~\ref{def:global-data}, the right-hand side is bounded by $C\varepsilon$.
\end{proposition}

The nonlinear equations \eqref{eq:global-density-system} satisfy \eqref{eq:global-forced-linear-system} with $v=\mathsf v$, $G=\mathscr G$, and $\vartheta=\sigma\mathcal N(f)$.
Therefore we can use Proposition~\ref{prop:global-linear-CL-proposition} to bound the time integral of $\mathcal Y(\mathsf v)$ in terms of the initial data and the bulk and capillary sources.
To obtain \eqref{eq:global-source-closed}, we estimate these sources together with the difference $\mathsf v-U$.
We use Lemma~\ref{lem:global-source-product} below to controls both $\mathscr G$ and $\mathsf v-U$ by products of $\nabla_{x,z}\eta$ and $U$, while the tangential pressure relation is used in Lemma~\ref{lem:global-source-capillary} to estimate the capillary contribution.
These bounds give small multiples of $\mathcal Y(U)$, which can be absorbed into the velocity estimate when $e_T$ is sufficiently small.
Local strong regularity ensures that the time integrals are finite.

\begin{lemma}[Product and coefficient estimates]\label{lem:global-source-product}
The velocity norm $\mathcal Y$ satisfies the product estimate
\begin{equation}\label{eq:global-source-algebra}
 \mathcal Y(hF)\le C\mathcal Y(h)\mathcal Y(F).
\end{equation}
For $m\ge2d+5$, one also has
\begin{equation}\label{eq:global-source-coefficient}
 \mathcal Y(\nabla_{x,z}\eta)\le C\|\eta\|_{X_m}.
\end{equation}
\end{lemma}

\begin{proof}
The one-dimensional product and embedding estimates give
\[
 \|ab\|_{B^{1/2}_{2,1}(I)}
 \le C\Big(\|a\|_{L^\infty}\|b\|_{B^{1/2}_{2,1}(I)} +\|b\|_{L^\infty}\|a\|_{B^{1/2}_{2,1}(I)}\Big),
 \qquad \|a\|_{L^\infty}\le C\|a\|_{B^{1/2}_{2,1}(I)}.
\]
Since $1+|\xi+\theta|\le (1+|\xi|)(1+|\theta|)$, Fourier convolution and the preceding estimates give
\[
\begin{aligned}
 &\int_{\R^d}(1+|\kappa|)
   \|\widehat{hF}(\kappa)\|_{B^{1/2}_{2,1}}\,\dd\kappa\le C\iint_{\R^d\times\R^d}(1+|\xi+\theta|)
   \|\widehat h(\xi)\|_{B^{1/2}_{2,1}}
   \|\widehat F(\theta)\|_{B^{1/2}_{2,1}}\,\dd\xi\dd\theta\\
 &\le C\left(\int_{\R^d}(1+|\xi|)\|\widehat h(\xi)\|_{B^{1/2}_{2,1}}\,\dd\xi\right)
   \left(\int_{\R^d}(1+|\theta|)\|\widehat F(\theta)\|_{B^{1/2}_{2,1}}\,\dd\theta\right)\le C\mathcal Y(h)\mathcal Y(F).
\end{aligned}
\]
For the $L_z^2$ part, we use the relation \[(1+|\kappa|)\langle\kappa\rangle^{1/2}\asymp (1+|\kappa|)^{3/2},
\]
and split the convolution according to the larger horizontal frequency.
On $|\xi|\ge|\theta|$, we have $(1+|\xi + \theta|)\le2 (1+|\xi|)$; on the complementary region, $(1+|\xi + \theta|)\le2 (1+|\theta|)$.
Thus
\[
\begin{aligned}
 &\int_{\R^d}(1+|\kappa|)\langle\kappa\rangle^{1/2}
   \|\widehat{hF}(\kappa)\|_{L_z^2}\,\dd\kappa\\
 \le& C\iint_{|\xi|\ge|\theta|}(1+|\xi|)^{3/2}
   \|\widehat h(\xi)\|_{L_z^2}\|\widehat F(\theta)\|_{L_z^\infty}\,\dd\xi\dd\theta+C\iint_{|\xi|<|\theta|}(1+|\theta|)^{3/2}
   \|\widehat h(\xi)\|_{L_z^\infty}\|\widehat F(\theta)\|_{L_z^2}\,\dd\xi\dd\theta\\
 \le& C\left(\int_{\R^d}(1+|\xi|)^{3/2}\|\widehat h(\xi)\|_{L_z^2}\,\dd\xi\right)
   \left(\int_{\R^d}\|\widehat F(\theta)\|_{B^{1/2}_{2,1}}\,\dd\theta\right)\\
   &+C\left(\int_{\R^d}\|\widehat h(\xi)\|_{B^{1/2}_{2,1}}\,\dd\xi\right)
   \left(\int_{\R^d}(1+|\theta|)^{3/2}\|\widehat F(\theta)\|_{L_z^2}\,\dd\theta\right)\le C\mathcal Y(h)\mathcal Y(F).
\end{aligned}
\]
Adding the two bounds proves \eqref{eq:global-source-algebra}.

For \eqref{eq:global-source-coefficient}, we use the vertical extension of Lemma~\ref{lem:global-embeddings} and set
\[
 E(\xi)^2=\int_\R(\langle r\rangle^{2m}+\zeta^4) |\widehat{E_z\eta}(\xi,\zeta)|^2\,\dd\zeta.
\]
Then $\|E\|_{L_\xi^2}\le C\|\eta\|_{X_m}$.
A vertical dyadic decomposition followed by the Cauchy--Schwarz inequality gives
\[
\begin{aligned}
 \|\widehat{\eta_z}\|_{B^{1/2}_{2,1}(I)}
    &\le C\langle r\rangle^{-m/4}E(\xi),&
 \|\widehat{\eta_z}\|_{L_z^2}
    &\le C\langle r\rangle^{-m/2}E(\xi),\\
 \|\widehat{\nabla_x\eta}\|_{B^{1/2}_{2,1}(I)}
    &\le Cr\langle r\rangle^{-3m/4}E(\xi),&
 \|\widehat{\nabla_x\eta}\|_{L_z^2}
    &\le Cr\langle r\rangle^{-m}E(\xi).
\end{aligned}
\]
For the Besov bounds, the relevant dyadic sums are $\sum_{N\ge1}N^3/(\langle r\rangle^{2m}+N^4)\le C\langle r\rangle^{-m/2}$ and $\sum_{N\ge1}N/(\langle r\rangle^{2m}+N^4)\le C\langle r\rangle^{-3m/2}$, where $N$ ranges over dyadic numbers and $N=1$ includes the low block.
For $m\ge2d+5$, inserting the weights in $\mathcal Y$ and applying Cauchy--Schwarz inequality in $\xi$ gives \eqref{eq:global-source-coefficient}.
\end{proof}

In the next lemma, we estimate the harmonic pressure velocity $U_\vartheta$, which is defined before Proposition~\ref{prop:global-linear-CL-proposition}, when the boundary pressure is the nonlinear curvature remainder.

\begin{lemma}[Capillary remainder]\label{lem:global-source-capillary}
Suppose $(f,\eta,U)$ satisfies \eqref{eq:global-coordinate-identities} and $\|f\|_{H^{m+1}}$ is sufficiently small.
Then
\begin{equation}\label{eq:global-source-capbound}
 \mathcal Y(U_{\sigma\mathcal N(f)})
 \le C\|f\|_{H^{m+1}}^2\mathcal Y(U).
\end{equation}
\end{lemma}

\begin{proof}
The harmonic pressure $q_\vartheta$ defined in \eqref{Defqvartheta} satisfies
\[
 \widehat q_\vartheta(\xi,z) =\frac{\cosh(r(z+1))}{\cosh r}\widehat\vartheta(\xi).
\]
Since $U_\vartheta=-\nabla q_\vartheta$, its horizontal Fourier transform is
\[
 \widehat U_\vartheta(\xi,z) =-\bigl(i\xi\widehat q_\vartheta,\partial_z\widehat q_\vartheta\bigr)
 =-\frac{\widehat\vartheta(\xi)}{\cosh r}
   \begin{pmatrix}
    i\xi\cosh(r(z+1))\\
    r\sinh(r(z+1))
   \end{pmatrix}.
\]
For $r\le1$ and $z\in I$, the hyperbolic ratios satisfy
\[
\cosh(r(z+1))/\cosh r\le1, \qquad |\sinh(r(z+1))|/\cosh r\le Cr.
\]
The formula for $\widehat U_\vartheta$ and its $z$-derivative therefore give
\[
 \|\widehat U_\vartheta\|_{L_z^2}\le Cr|\widehat\vartheta|, \qquad
 \|\partial_z\widehat U_\vartheta\|_{L_z^2}\le Cr^2|\widehat\vartheta|, \qquad 
 \|\widehat U_\vartheta\|_{H_z^1}\le Cr|\widehat\vartheta|
\]
For $r\ge1$, the same formula gives 
\[
\|\widehat U_\vartheta\|_{L^2}\le Cr^{1/2}|\widehat\vartheta|, \qquad \|\partial_z\widehat U_\vartheta\|_{L^2}\le Cr^{3/2}|\widehat\vartheta|.
\]
Interpolating between $L^2$ and $H^1$ yields
\begin{equation}\label{eq:global-source-harmonic}
 \mathcal Y(U_\vartheta) \le C\int_{\R^d}r(1+r)|\widehat\vartheta(\xi)|\,\dd\xi.
\end{equation}
We also use the endpoint trace estimate
\begin{equation}\label{eq:global-source-trace}
 |v(0)|\le C\|v\|_{B^{1/2}_{2,1}(I)},
\end{equation}
which follows from Bernstein's inequality and summation of the vertical dyadic blocks.

Set the auxiliary control norm
\[
 A(f)=\int_{\R^d}r(1+r)(c+\sigma r^2)|\widehat f(\xi)|\,\dd\xi.
\]
The curvature remainder satisfies
\begin{equation}\label{eq:global-source-curvature-weight}
 \int_{\R^d}r(1+r) |\widehat{\sigma\mathcal N(f)}(\xi)|\,\dd\xi \le C\|f\|_{H^{m+1}}^2 A(f).
\end{equation}
We write the curvature remainder as $\mathcal N(f)=\diver_xF(p)$, where
\[
 p=\nabla_xf,\qquad
 F(p)=\bigl[1-(1+|p|^2)^{-1/2}\bigr]p =\sum_{k\ge1}c_k|p|^{2k}p,
\]
For $\omega(\xi)=|\xi|^2(1+|\xi|)$, the convolution inequality
\[
 \omega(\xi_1+\cdots+\xi_n) \le n^3\sum_{j=1}^n\omega(\xi_j)
\]
allows us to put the weight on one factor in each term of the series and bound the other factors in Fourier $L^1$.
The series starts at degree three and converges for small $\|\widehat p\|_{L_\xi^1}$.
Summation gives
\[
 \|\omega\widehat{F(p)}\|_{L_\xi^1} \le C\|\widehat p\|_{L_\xi^1}^2 \|\omega\widehat p\|_{L_\xi^1}.
\]
Moreover,
\[
\|\widehat p\|_{L_\xi^1}\le C\|f\|_{H^{m+1}}, \qquad \sigma\|\omega\widehat p\|_{L_\xi^1} =\int_{\R^d}\sigma r^3(1+r)|\widehat f(\xi)|\,\dd\xi \le A(f),
\]
These bounds imply \eqref{eq:global-source-curvature-weight}.
To control $A(f)$ by the velocity, we use the pressure identity on the surface:
\[
 \nabla_x(c-\sigma\Delta_x)f =-U_x|_{z=0}-U_{d+1}|_{z=0}\nabla_xf-\sigma\nabla_x\mathcal N(f).
\]
By \eqref{eq:global-source-trace}, horizontal Fourier convolution, and $\int_{\R^d}(1+|\xi|)|\widehat{\nabla_xf}|\,\dd\xi\le C\|f\|_{H^{m+1}}$,
\[
 A(f)\le C(1+\|f\|_{H^{m+1}})\mathcal Y(U) +C\|f\|_{H^{m+1}}^2A(f).
\]
Since $m>d/2+3$, the control norm $A(f)$ is finite.
For $\|f\|_{H^{m+1}}$ sufficiently small, we absorb the last term and obtain
\begin{equation}\label{eq:global-source-surfacecontrol}
 A(f)\le C\mathcal Y(U).
\end{equation}
Together with \eqref{eq:global-source-curvature-weight} and \eqref{eq:global-source-harmonic}, this proves \eqref{eq:global-source-capbound}.
\end{proof}

\begin{proof}[Proof of Proposition~\ref{prop:global-source-closure}]
From the definitions of $\mathscr G$ and $\mathsf v$,
\[
 \mathscr G=(\nabla_x\eta\,U_{d+1}-\eta_zU_x, \nabla_x\eta\cdot U_x+\eta_zU_{d+1}),
 \qquad \mathsf v-U=(\eta_zU_x,-\nabla_x\eta\cdot U_x).
\]
Taking $\delta\le1$, Lemma~\ref{lem:global-source-product} and Lemma~\ref{lem:global-source-capillary} give
\begin{equation}\label{eq:global-source-Gsmall}
 \mathcal Y(\mathsf v-U)+\mathcal Y(\mathscr G) +\mathcal Y(U_{\sigma\mathcal N(f)}) \le Ce_T\mathcal Y(U).
\end{equation}
The terms in \eqref{eq:global-source-Gsmall} are time-integrable on compact local existence intervals.
Indeed, local regularity gives $U,\mathsf v,\mathscr G\in L_T^2H^{m-1}$, and $\mathcal Y(F)\le C\|F\|_{H^{d+3}}$ with $m-1\ge d+3$.
The integrability of the capillary source follows from \eqref{eq:global-source-Gsmall}.

Using the identities 
\[
\eta_t=\mathsf v_{d+1}, \qquad f_t=\mathsf v_{d+1}|_{z=0},
\] 
the nonlinear solution satisfies \eqref{eq:global-forced-linear-system} with $v=\mathsf v$, $G=\mathscr G$, and $\vartheta=\sigma\mathcal N(f)$ by \eqref{eq:global-Piola-source} and \eqref{eq:global-nonlinear-pressure}.
We apply Proposition~\ref{prop:global-linear-CL-proposition} to these sources and use the triangle inequality and \eqref{eq:global-source-Gsmall} to find
\[
 \begin{aligned}
 \int_0^T\mathcal Y(U)\,\dd t
 &\le\int_0^T\mathcal Y(\mathsf v)\,\dd t
          +Ce_T\int_0^T\mathcal Y(U)\,\dd t\\
 &\le C\bigl(\|\eta(0)\|_{H^m}
                 +\mathcal B_{\rm low}(\eta(0))\bigr)
       +Ce_T\int_0^T\mathcal Y(U)\,\dd t.
 \end{aligned}
\]
Choosing $C\delta\le1/2$ allows us to absorb the last term.
Since $\mathcal Y(\mathsf v)\le(1+Ce_T)\mathcal Y(U)$ by \eqref{eq:global-source-Gsmall}, this proves \eqref{eq:global-source-closed}.
The initial coordinate estimate \eqref{eq:global-initial-coordinate} then bounds its right-hand side by $C\varepsilon$.
\end{proof}

\subsection{Anisotropic energy estimates}\label{sec:global-energy}

We combine the relative energy with normal and tangential energy estimates to prove the following energy inequality.

\begin{proposition}[Anisotropic energy inequality]\label{prop:global-energy}
Let $\mathcal E$ be the energy functional constructed in \eqref{eq:global-total-energy} below.
There exist constants $\delta_0,\kappa,C>0$, depending only on the fixed $d,m,\sigma>0$ and profile bounds, such that every smooth solution with $e(t)\le\delta_0$ satisfies
\begin{equation}\label{eq:global-energy-inequality}
 C^{-1}e\le\mathcal E\le Ce,\qquad
 \mathcal E'+\kappa\mathcal D
       \le Ce^{1/2}\mathcal D+C\mathcal W_d(U)e.
\end{equation}

Suppose a local strong solution in \eqref{eq:global-local-class} admits a regular density coordinate map on $[0,T]$, with $\sup_{t\le T}e(t)$ sufficiently small.
Then $\mathcal E\in C([0,T])$, and for $0\le t_0\le t_1\le T$,
\begin{equation}\label{eq:global-integrated-energy}
 \mathcal E(t_1)+\kappa\int_{t_0}^{t_1}\mathcal D(t)\dd t
 \le\mathcal E(t_0)
 +C\int_{t_0}^{t_1}
       \bigl(e(t)^{1/2}\mathcal D(t)+\mathcal W_d(U(t))e(t)\bigr)\dd t.
\end{equation}
\end{proposition}

\subsubsection{Normal energy estimate}\label{sec:global-normal-section}

To control $\eta_{zz}$, we estimate the physical Laplacian of the function $h=y-\gamma^{-1}(\rho)$.
The comparison below recovers $\eta_{zz}$ from $\Delta h$ and the lower-order and horizontal norms of $\eta$.
The material equation for $h$ also yields a cancellation of the higher derivatives of the background profile.

We first work with smooth, rapidly decaying solutions for which the density coordinate map is a diffeomorphism and the boundary densities satisfy \eqref{eq:global-boundary-density}.
Set
\begin{equation}\label{eq:global-normal-label}
 z_\rho=\gamma^{-1}(\rho),\qquad h=y-\gamma^{-1}(\rho),\qquad
 D_t=\partial_t+u\cdot\nabla_{x,y}.
\end{equation}
Since $D_t\rho=0$ and $D_ty=u_{d+1}$, these quantities satisfy
\begin{equation}\label{eq:global-normal-material}
 D_tz_\rho=0,\qquad D_th=u_{d+1},\qquad D_t\gamma'(z_\rho)=0,
 \qquad h\circ\Phi_\eta=\eta,
\end{equation}
and the boundary conditions give $h=f$ on $\Sig_f$ and $h=0$ on $\{y=-1\}$.
In particular, $\nabla_xh=0$ on the bottom.

\begin{lemma}\label{lem:global-normal-cancellation}
Any smooth solution of \eqref{e:IPM} satisfies the elliptic vorticity-type relations
\begin{equation}\label{eq:global-normal-elliptic-identities}
 \begin{aligned}
  &\Delta u_x =-\partial_y\bigl(\gamma'(z_\rho)\nabla_xh\bigr) =-\gamma'(z_\rho)\nabla_xh_y-\gamma''(z_\rho)(1-h_y)\nabla_xh,\\
  &\Delta u_{d+1} =\diver_x \bigl(\gamma'(z_\rho)\nabla_xh\bigr) =\gamma'(z_\rho)\Delta_xh-\gamma''(z_\rho)|\nabla_xh|^2,\\
  &\partial_y u_x-\nabla_xu_{d+1}=-\gamma'(z_\rho)\nabla_xh,
 \end{aligned}
\end{equation}
and the Laplacian $\Delta h=\Delta_{x,y}h$ obeys the material transport equation
\begin{equation}\label{eq:global-normal-laplacian-transport}
 \begin{aligned}
 D_t\Delta h={}&\gamma'(z_\rho)(1-h_y)\Delta_xh -2D_xu_x:D_x^2h-4\nabla_xu_{d+1}\cdot\nabla_xh_y\\
 &+3\gamma'(z_\rho)\nabla_xh\cdot\nabla_xh_y-2(\partial_y u_{d+1})h_{yy}.
 \end{aligned}
\end{equation}
\end{lemma}

\begin{proof}
Taking the curl and divergence of Darcy's law $u+\nabla p=-\rho e_{d+1}$ with $\diver u=0$ gives $\Delta u_x=\nabla_x\rho_y$, $\Delta u_{d+1}=-\Delta_x\rho$, and $\partial_y u_x-\nabla_xu_{d+1}=\nabla_x\rho$.
Since $\rho=\gamma(y-h)$, the chain rule gives $\nabla_x\rho=-\gamma'(z_\rho)\nabla_xh$ and $\partial_y z_\rho=1-h_y$, which proves \eqref{eq:global-normal-elliptic-identities}.
To obtain the equation for $\Delta h$, we commute $D_t$ with $\Delta$ in $D_th=u_{d+1}$ and use $\diver u=0$:
\[
 D_t\Delta h=\Delta u_{d+1}-(\Delta u)\cdot\nabla h-2\nabla u:\nabla^2h.
\]
The first two equations in \eqref{eq:global-normal-elliptic-identities} show that the terms involving $\gamma''(z_\rho)$ in $\Delta u_{d+1}$ and $(\Delta u)\cdot\nabla h$ both equal $-\gamma''(z_\rho)|\nabla_xh|^2$.
These terms cancel, leaving
\[
 \Delta u_{d+1}-(\Delta u)\cdot\nabla h =\gamma'(z_\rho)(1-h_y)\Delta_xh +\gamma'(z_\rho)\nabla_xh\cdot\nabla_xh_y.
\]
Substituting $\partial_y u_x=\nabla_xu_{d+1}-\gamma'(z_\rho)\nabla_xh$ into the mixed derivative terms of $-2\nabla u:\nabla^2h$ then gives \eqref{eq:global-normal-laplacian-transport}.
\end{proof}

We multiply \eqref{eq:global-normal-laplacian-transport}
by $\Delta h/(-\gamma'(z_\rho))$ and integrate over $\Om_f$.
Since $D_t\gamma'(z_\rho)=0$, the weight is transported;
it also cancels the factor $\gamma'(z_\rho)$ in the principal term.
Integration by parts produces a boundary term on $\Sig_f$.
Set
\begin{equation}\label{eq:global-normal-beta}
 \nu=1-h_y|_{\Sig_f}=\partial_y z_\rho|_{\Sig_f}>0,
 \qquad
 \beta=\frac{1-\nu^3}{3}.
\end{equation}
The leading boundary term is
$\langle\nabla_xf,\nabla_x\beta\rangle$.
By \eqref{eq:global-normal-beta-transport} and
\eqref{eq:global-normal-surface-velocity}, the leading term
in $\beta_t$ is $(c-\sigma\Delta_x)\Delta_xf$.
Pairing it with $(c-\sigma\Delta_x)^{-1}\beta$ gives
$\bigl\langle(c-\sigma\Delta_x)^{-1}\beta, (c-\sigma\Delta_x)\Delta_xf\bigr\rangle =-\langle\nabla_x\beta,\nabla_xf\rangle$.
Thus the time derivative of $\frac12\langle(c-\sigma\Delta_x)^{-1}\beta,\beta\rangle$
cancels the leading boundary term.
We therefore define
\begin{equation}\label{eq:global-normal-energy}
 \mathcal Q_h =\frac12\int_{\Om_f} \frac{(\Delta h)^2}{-\gamma'(z_\rho)}\,\dd x\dd y +\frac12\langle(c-\sigma\Delta_x)^{-1}\beta,\beta \rangle_{L_x^2}.
\end{equation}

\begin{proposition}[Exact normal energy identity]\label{prop:global-normal-exact}
Write 
\begin{equation}\label{eq:global-normal-surface-velocity}
 \begin{aligned}
 &U_x|_{z=0}=-(c-\sigma\Delta_x)\nabla_xf+R_x,\\
 R_x&= \frac{\nabla_xf\,[\nabla_xf\cdot(c-\sigma\Delta_x)\nabla_xf]}{1+|\nabla_xf|^2}
 -\frac{\nabla_xf\,f_t}{1+|\nabla_xf|^2}-\sigma(I_d+\nabla_xf\otimes\nabla_xf)^{-1}\nabla_x\mathcal N(f).
 \end{aligned}
\end{equation}
Then $\mathcal Q_h$ satisfies the exact identity
\begin{equation}\label{eq:global-normal-identity}
 \begin{aligned}
 \mathcal Q_h'
 &+\int_{\Om_f}(1-h_y) \bigl(|D_x^2h|^2+|\nabla_xh_y|^2\bigr)\dd x\dd y\\
 ={}&-\int_{\Om_f}\frac{\Delta h}{-\gamma'(z_\rho)} \bigl(2D_xu_x:D_x^2h+4\nabla_xu_{d+1}\cdot\nabla_xh_y +2(\partial_yu_{d+1})h_{yy}\bigr)\dd x\dd y\\
 &-3\int_{\Om_f}(\Delta h)\nabla_xh\cdot\nabla_xh_y\dd x\dd y+\int_{\Om_f}\nabla_xh\cdot \bigl(D_x^2h-(\Delta_xh)I_d\bigr)\nabla_xh_y\dd x\dd y\\
 &-\int_{\R^d}\nu^3 \bigl((\nabla_xf)^{\mathsf T}(D_x^2f)\nabla_xf -|\nabla_xf|^2\Delta_xf\bigr)\dd x\\
 &+\langle\nabla_x(c-\sigma\Delta_x)^{-1}\beta,R_x\rangle_{L_x^2}+\langle(c-\sigma\Delta_x)^{-1}\beta, 3\beta\diver_x(U_x|_{z=0})-(U_x|_{z=0})\cdot\nabla_x\beta\rangle_{L_x^2}.
 \end{aligned}
\end{equation}
where the third and fourth integrals on the right-hand side vanish identically when $d=1$.
\end{proposition}

\begin{proof}
Differentiating $h(t,x,f(t,x))=f(t,x)$ in $x$ gives $\nabla_xh|_{\Sig_f}=(1-h_y|_{\Sig_f})\nabla_xf=\nu\nabla_xf$.
We differentiate $D_th=u_{d+1}$ with respect to $y$ and restrict to $\Sig_f$.
Using $\diver_{x,y}u=0$ and the chain rule along $\Sig_f$, we obtain
\[
 \bigl(\partial_t+(U_x|_{z=0})\cdot\nabla_x\bigr)(h_y|_{\Sig_f}) =-\nu\diver_x(U_x|_{z=0}).
\]
Multiplying by $\nu^2$ yields the evolution equation for $\beta=(1-\nu^3)/3$:
\begin{equation}\label{eq:global-normal-beta-transport}
 \beta_t=-\diver_x(U_x|_{z=0})+3\beta\diver_x(U_x|_{z=0})  -(U_x|_{z=0})\cdot\nabla_x\beta.
\end{equation}
Taking the tangential derivative of $p|_{\Sig_f}=\sigma H(f)$ and using Darcy's law along $\Sig_f$ gives $U_x|_{z=0}+(U_{d+1}|_{z=0})\nabla_xf=-\nabla_x(cf+\sigma H(f))$.
Eliminating $U_{d+1}|_{z=0}=f_t+(U_x|_{z=0})\cdot\nabla_xf$ and inverting $I_d+\nabla_xf\otimes\nabla_xf$ produces the decomposition \eqref{eq:global-normal-surface-velocity}.

To identify the boundary flux, we use the divergence identities
 \begin{align*}
 (\Delta_xh)^2-|D_x^2h|^2 &=\diver_x\bigl((\nabla_xh)\Delta_xh-(D_x^2h)\nabla_xh\bigr),\\
 h_{yy}\Delta_xh-|\nabla_xh_y|^2 &=\diver_x(h_{yy}\nabla_xh) -\partial_y(\nabla_xh_y\cdot\nabla_xh),
 \end{align*}
Multiplying these identities by $h_y$ and integrating by parts over $\Om_f$ gives
\begin{align*}
 &\int_{\Om_f}h_y(\Delta h)(\Delta_xh)\dd x\dd y
  =\int_{\Om_f}h_y \bigl(|D_x^2h|^2+|\nabla_xh_y|^2\bigr)\dd x\dd y\\
 &+\int_{\Om_f}\nabla_xh\cdot\bigl(D_x^2h-(\Delta_xh)I_d\bigr)\nabla_xh_y\dd x\dd y\\
 &+\int_{\R^d}(h_y|_{\Sig_f})\bigl\{-\nabla_xh\cdot\nabla_x(h_y|_{\Sig_f})+\nabla_xh\cdot(D_x^2h)\nabla_xf -(\nabla_xh\cdot\nabla_xf)\Delta_xh\bigr\}\dd x,
\end{align*}
where all surface traces in the last integral are evaluated on $\Sig_f$ and the bottom boundary contribution vanishes because $\nabla_xh|_{y=-1}=0$.
Applying the same integration by parts to $-\int_{\Om_f}\Delta h\,\Delta_xh\dd x\dd y$ without the $h_y$ weight, and using the tangential derivative of $\nabla_xh|_{\Sig_f}=\nu\nabla_xf$, gives the following upper-boundary contribution from $-\int_{\Om_f}(1-h_y)\Delta h\,\Delta_xh\dd x\dd y$:
\[
 \int_{\R^d}\nabla_xf\cdot\nabla_x\beta\dd x -\int_{\R^d}\nu^3 \bigl((\nabla_xf)^{\mathsf T}(D_x^2f)\nabla_xf -|\nabla_xf|^2\Delta_xf\bigr)\dd x.
\]

Multiplying \eqref{eq:global-normal-laplacian-transport} by $\Delta h/(-\gamma'(z_\rho))$, integrating over $\Om_f$, and applying the Reynolds transport theorem (using $D_t\gamma'(z_\rho)=0$ and $\diver u=0$) gives the time derivative of the bulk integral in $\mathcal Q_h$.
By \eqref{eq:global-normal-beta-transport} and \eqref{eq:global-normal-surface-velocity}, the leading term in the time derivative of the boundary quadratic form in $\mathcal Q_h$ is
\[
 \langle \nabla_x(c-\sigma\Delta_x)^{-1}\beta,U_x|_{z=0}\rangle
 =-\langle \nabla_x\beta,\nabla_xf\rangle +\langle \nabla_x(c-\sigma\Delta_x)^{-1}\beta,R_x\rangle,
\]
and the first term cancels the leading bulk boundary flux $\int_{\R^d}\nabla_xf\cdot\nabla_x\beta\dd x$.
Therefore, we obtain the exact identity \eqref{eq:global-normal-identity}.
\end{proof}

To use  the identity \eqref{eq:global-normal-identity}, we first relate $\mathcal Q_h$ and its dissipation to the second vertical derivatives of $\eta$ and $U$ in $X_m$.

\begin{lemma}[Coercivity and dissipation equivalence]\label{lem:global-normal-comparison}
Let $m>2d$, and suppose that the total energy $e(t)$ defined in \eqref{eq:global-ed} is sufficiently small.
Then
\begin{equation}\label{eq:global-normal-coercivity}
 \|\eta\|_{X_m}^2 \asymp \|\eta\|_{L^2}^2+\||D_x|^m\eta\|_{L^2}^2+\mathcal Q_h, \qquad
 \|\nabla^2h\|_{L^2(\Om_f)}+\|\beta\|_{L_x^2}\le Ce^{1/2},
\end{equation}
and the normal dissipation
\begin{equation}\label{eq:global-normal-dissipation}
 D_h:=\|D_x^2h\|_{L^2(\Om_f)}^2 +\|\nabla_xh_y\|_{L^2(\Om_f)}^2
\end{equation}
satisfies the equivalence
\begin{equation}\label{eq:global-normal-dissipation-comparison}
 \mathcal D\asymp\|U\|_{L^2}^2+\||D_x|^mU\|_{L^2}^2+D_h.
\end{equation}
In addition, the gradients of the velocity and of $h$ and $\eta$ satisfy the pointwise bounds
\begin{equation}\label{eq:global-normal-velocity-bounds}
 \begin{aligned}
 &\|\nabla_xh\|_{L^\infty}+\|\nabla_x\eta\|_{L^\infty}\le C\min(e^{1/2},\mathcal D^{1/2}),\\
& \|D_xu_x\|_{L^\infty}+\|\nabla_xu_{d+1}\|_{L^\infty} +\|\partial_y u_{d+1}\|_{L^\infty}\le C\bigl(\mathcal W_d(U)+\mathcal D\bigr).
 \end{aligned}
\end{equation}
\end{lemma}

\begin{proof}
Using the chain rule, we compute the pullback of $\Delta h$ to the reference strip $\mathcal S$:
\begin{equation}\label{eq:global-normal-pullback-laplacian}
 (\Delta h)\circ\Phi_\eta =\frac{\Delta_x\eta}{J_\eta} -\frac{2\nabla_x\eta\cdot\nabla_x\eta_z}{J_\eta^2} +\frac{1+|\nabla_x\eta|^2}{J_\eta^3}\eta_{zz}.
\end{equation}
$\|\Delta_x\eta\|_{L^2}$ can be controlled by $\|\eta\|_{L^2}+\||D_x|^m\eta\|_{L^2}$ using horizontal interpolation.
Applying Fourier interpolation to $E_z\eta$, with $E_z$ as in Lemma~\ref{lem:global-embeddings}, also gives
\[
 \|\nabla_x\eta_z\|_{L^2} \le C\bigl(\|\eta\|_{L^2} +\||D_x|^m\eta\|_{L^2}+\|\eta_{zz}\|_{L^2}\bigr).
\]
Since $J_\eta\asymp1$ and $\|\nabla_x\eta\|_{L^\infty}\le Ce^{1/2}$, solving \eqref{eq:global-normal-pullback-laplacian} for $\eta_{zz}$ yields
\[
 \|\eta_{zz}\|_{L^2} \le C\bigl(\|\Delta h\|_{L^2(\Om_f)} +\|\eta\|_{L^2}+\||D_x|^m\eta\|_{L^2}\bigr)+Ce^{1/2}\|\eta_{zz}\|_{L^2}.
\]
If $Ce^{1/2}\le\frac12$, the last term can be absorbed into the left-hand side.
The reverse estimate follows from the same pullback formula.
Moreover, $-\gamma'(z_\rho)\asymp1$ and $\nu=J_\eta^{-1}|_{z=0}$.
Hence $\|\beta\|_{L^2}\le C\|\eta_z|_{z=0}\|_{L^2}\le C\|\eta\|_{X_m}$.
These bounds prove \eqref{eq:global-normal-coercivity}.
The chain rule also gives $\|\nabla^2h\|_{L^2(\Om_f)}\le Ce^{1/2}$.

By Lemma~\ref{lem:global-embeddings}, we have $\|\nabla_{x,z}U\|_{L^\infty}\le C\mathcal D^{1/2}$.
Combined with the curl identity $\partial_y u_x-\nabla_xu_{d+1}=-\gamma'(z_\rho)\nabla_xh$ in \eqref{eq:global-normal-elliptic-identities} and the relation $\nabla_xh\circ\Phi_\eta=\nabla_x\eta/J_\eta$, this proves the first line of \eqref{eq:global-normal-velocity-bounds}.
For the horizontal derivatives in physical coordinates, differentiating $U=u\circ\Phi_\eta$ gives
\[
 (\partial_{x_j}u_i)\circ\Phi_\eta
 =\partial_{x_j}U_i-\frac{\partial_zU_i}{J_\eta}\partial_{x_j}\eta,
 \qquad 1\le i\le d+1,\quad 1\le j\le d.
\]
For any extension $\widetilde U$ of $U$ to $\R^d\times\R$, Fourier inversion yields
\[
 \|D_xU\|_{L^\infty(\mathcal S)}
 \le C\int_{\R^d}\int_\R |\xi|\,|\widehat{\widetilde U}(\xi,\zeta)|\dd\zeta\dd\xi
 \le C\int_{\R^d}\int_\R(1+|\xi|)|\widehat{\widetilde U}(\xi,\zeta)|\dd\zeta\dd\xi.
\]
Taking the infimum over these extensions gives $\|D_xU\|_{L^\infty}\le C\mathcal W_d(U)$.
Since $J_\eta^{-1}$ is bounded, the preceding chain rule and the first line of \eqref{eq:global-normal-velocity-bounds} imply
\begin{align*}
 &\|D_xu\|_{L^\infty(\Om_f)}=\|(D_xu)\circ\Phi_\eta\|_{L^\infty(\mathcal S)}\\
 \le&\|D_xU\|_{L^\infty} +\|J_\eta^{-1}\|_{L^\infty}\|\partial_zU\|_{L^\infty}\|\nabla_x\eta\|_{L^\infty} \leq C\bigl(\mathcal W_d(U)+\mathcal D\bigr).
\end{align*}
Finally, incompressibility gives $\partial_yu_{d+1}=-\sum_{j=1}^d\partial_{x_j}u_j$, so that
\[
 \|D_xu_x\|_{L^\infty}+\|\nabla_xu_{d+1}\|_{L^\infty} +\|\partial_yu_{d+1}\|_{L^\infty}\le C\|D_xu\|_{L^\infty} \le C\bigl(\mathcal W_d(U)+\mathcal D\bigr),
\]
which proves the second line of \eqref{eq:global-normal-velocity-bounds}.

For the dissipation comparison, we pull $\Delta u$ back to density coordinates:
\begin{equation}\label{eq:global-normal-velocity-laplacian}
 (\Delta u)\circ\Phi_\eta =\Delta_xU-\frac{2\nabla_x\eta}{J_\eta}\cdot\nabla_x\partial_z U +\frac{1+|\nabla_x\eta|^2}{J_\eta^2}\partial_z^2 U -\bigl((\Delta h)\circ\Phi_\eta\bigr)\partial_z U.
\end{equation}
Interpolation and $\|\Delta h\|_{L^2}\le Ce^{1/2}$ from \eqref{eq:global-normal-coercivity} give $\|\Delta u\|_{L^2(\Om_f)}\le C\mathcal D^{1/2}$.
On each vertical segment $\{x\}\times(-1,f(t,x))$, the quantity $-\gamma'(z_\rho)\nabla_xh$ vanishes at $y=-1$; its vertical derivative is $\Delta u_x$ by \eqref{eq:global-normal-elliptic-identities}.
Applying the one-dimensional Poincar\'e inequality in $y$ and using the boundedness of $\partial_y(-\gamma'(z_\rho))=-\gamma''(z_\rho)(1-h_y)$ therefore yields
\begin{equation}\label{eq:global-normal-bottom-poincare}
 \|\nabla_xh\|_{L^2(\Om_f)}+\|\nabla_xh_y\|_{L^2(\Om_f)}
             \le C\|\Delta u_x\|_{L^2(\Om_f)}.
\end{equation}
From the second identity in \eqref{eq:global-normal-elliptic-identities}, we then deduce $\|\Delta_xh\|_{L^2(\Om_f)}\le C\|\Delta u\|_{L^2(\Om_f)}\le C\mathcal D^{1/2}$.
To pass from $\Delta_xh$ to the full horizontal Hessian $D_x^2h$, we express $(D_x^2h)\circ\Phi_\eta$ as
\[
 (D_x^2h)\circ\Phi_\eta =\frac{D_x^2\eta}{J_\eta} -\frac{\nabla_x\eta\otimes\nabla_x\eta_z+\nabla_x\eta_z\otimes \nabla_x\eta}{J_\eta^2} +\frac{\eta_{zz}}{J_\eta^3}\nabla_x\eta\otimes \nabla_x\eta.
\]
The first bound in \eqref{eq:global-normal-velocity-bounds}, together with $\|\nabla_x\eta_z\|_{L^2}+\|\eta_{zz}\|_{L^2}\le Ce^{1/2}$, gives
\[ \|\nabla_x\eta\otimes\nabla_x\eta_z\|_{L^2}
 +\||\nabla_x\eta|^2\eta_{zz}\|_{L^2}
 \le Ce^{1/2}\mathcal D^{1/2}.
\]
Taking the matrix trace of the Hessian formula and using
$\|D_x^2\eta\|_{L^2}=\|\Delta_x\eta\|_{L^2}$ on $\mathcal S$, we obtain
\[
 \|D_x^2h\|_{L^2(\Om_f)} \le C\|\Delta_xh\|_{L^2(\Om_f)}+Ce^{1/2}\mathcal D^{1/2},
\]
which proves $D_h\le C\mathcal D$.

Conversely, the vertical Poincar\'e inequality and \eqref{eq:global-normal-elliptic-identities} give $\|\Delta u\|_{L^2(\Om_f)}\le CD_h^{1/2}$.
To recover $U_{zz}$ from $\Delta u$, we use horizontal and mixed interpolation:
\[
 \|\Delta_xU\|_{L^2}\le C\bigl(\|U\|_{L^2}+\||D_x|^mU\|_{L^2}\bigr),\qquad
 \|\nabla_xU_z\|_{L^2} \le C\bigl(\|U\|_{L^2}+\||D_x|^mU\|_{L^2}+\|U_{zz}\|_{L^2}\bigr).
\]
The anisotropic embedding gives the same bound for $\|U_z\|_{L^\infty}$.
Thus the last product in \eqref{eq:global-normal-velocity-laplacian} is bounded by $Ce^{1/2}(\|U\|_{L^2}+\||D_x|^mU\|_{L^2}+\|U_{zz}\|_{L^2})$.
Solving that identity for $U_{zz}$ now yields
\[
 \|U_{zz}\|_{L^2} \le C\bigl(D_h^{1/2}+\|U\|_{L^2}+\||D_x|^mU\|_{L^2}\bigr) +Ce^{1/2}\|U_{zz}\|_{L^2}.
\]
Absorbing the last term proves the reverse inequality in \eqref{eq:global-normal-dissipation-comparison}.
\end{proof}

The following curvature and low-order surface estimates then control the boundary terms in \eqref{eq:global-normal-identity}.

\begin{lemma}[Curvature bounds]\label{lem:global-curvature-homogeneous}
For every smooth function $f$ with $\|\nabla_xf\|_{L^\infty}\le 1$, one has
\begin{equation}\label{eq:global-curvature-homogeneous}
 \|\nabla_x\mathcal N(f)\|_{L^2}
 \le C\|\nabla_xf\|_{L^\infty}^2\|\nabla_x^3f\|_{L^2}.
\end{equation}
For each $q\ge0$, the remainder also satisfies the estimate
\begin{equation}\label{eq:global-curvature-tame}
 \|\mathcal N(f)\|_{H^q(\R^d)}\le C_q\|\nabla_xf\|_{L^\infty}^2\|f\|_{H^{q+2}(\R^d)}.
\end{equation}
\end{lemma} 

\begin{proof}
Set $b=1-(1+|\nabla_xf|^2)^{-1/2}$, so that $\mathcal N(f)=\diver_x(b\nabla_xf)$.
The tame product and composition estimates
\cite[Corollary~2.54 and Theorems~2.61,~2.89]{BCD2011}
give, for $X=\dot H^2(\R^d)$ or $X=H^{q+1}(\R^d)$,
\[
 \|b\|_{L^\infty}\le C\|\nabla_xf\|_{L^\infty}^2, \qquad
 \|b\|_X \le C_X\||\nabla_xf|^2\|_X \le C_X\|\nabla_xf\|_{L^\infty}\|\nabla_xf\|_X.
\]
Here the composition constants are uniform because
$\|\nabla_xf\|_{L^\infty}\le1$.
Applying the product estimate once more, we obtain
\[
 \|b\nabla_xf\|_X \le C_X\bigl( \|b\|_{L^\infty}\|\nabla_xf\|_X +\|\nabla_xf\|_{L^\infty}\|b\|_X \bigr) \le C_X\|\nabla_xf\|_{L^\infty}^2\|\nabla_xf\|_X.
\]
Taking $X=\dot H^2$ and using $\mathcal N(f)=\diver_x(b\nabla_xf)$ proves
\eqref{eq:global-curvature-homogeneous}.
Choosing $X=H^{q+1}$, we get \eqref{eq:global-curvature-tame}.
\end{proof}

\begin{lemma}\label{lem:global-surface-low}
Whenever the total energy $e$ is sufficiently small, the free surface satisfies
\begin{equation}\label{eq:global-reconstruction}
 \|\nabla_xf\|_{H^2(\R^d)}+\|f_t\|_{L^2(\R^d)}\le C\mathcal D^{1/2}, \qquad
 \|\nabla_xf\|_{L^\infty(\R^d)}\le C\min(e^{1/2},\mathcal D^{1/2}).
\end{equation}
\end{lemma}

\begin{proof}
Taking the tangential component of Darcy's law along $\Sig_f$, we get
\[
 (c-\sigma\Delta_x)\nabla_xf =-U_x|_{z=0}-U_{d+1}|_{z=0}\nabla_xf-\sigma\nabla_x\mathcal N(f).
\]
The trace inequality on $\mathcal S$ gives $\|U|_{z=0}\|_{L^2}\le C\mathcal D^{1/2}$.
Since $\nabla_x\eta|_{z=0}=\nabla_xf$, Lemma~\ref{lem:global-normal-comparison} also gives $\|\nabla_xf\|_{L^\infty}\le C\min(e^{1/2},\mathcal D^{1/2})$.
The homogeneous curvature bound \eqref{eq:global-curvature-homogeneous} allows us to absorb the last term in the tangential Darcy relation when $\|\nabla_xf\|_{L^\infty}$ is small, giving $\|\nabla_xf\|_{H^2}\le C\mathcal D^{1/2}$.
Using the kinematic condition $f_t=U_{d+1}|_{z=0}-U_x|_{z=0}\cdot\nabla_xf$, we get $\|f_t\|_{L^2}\le C\mathcal D^{1/2}$.
\end{proof}

\begin{lemma}[Normal energy inequality]\label{prop:global-normal-estimate}
Under the smoothness and smallness assumptions of Lemma~\ref{lem:global-normal-comparison}, there exist constants $c_0,C>0$ such that
\begin{equation}\label{eq:global-normal-estimate-equation}
 \mathcal Q_h'+c_0D_h \le Ce^{1/2}\mathcal D+C\mathcal W_d(U)e.
\end{equation}
\end{lemma}

\begin{proof}
Because $1-h_y=\partial_yz_\rho\ge c_0>0$ for small $e$, the dissipative integral on the left-hand side of \eqref{eq:global-normal-identity} is bounded below by $c_0D_h$.
By Lemma~\ref{lem:global-normal-comparison}, the bulk integral containing derivatives of $u$ is bounded by $C(\mathcal W_d(U)+\mathcal D)e$.
The other two bulk integrals satisfy
\[
 C\|\Delta h\|_{L^2}\|\nabla_xh\|_{L^\infty} \|\nabla_xh_y\|_{L^2}\le Ce^{1/2}\mathcal D, \qquad C\|\nabla_xh\|_{L^\infty}\|D_x^2h\|_{L^2} \|\nabla_xh_y\|_{L^2}\le Ce^{1/2}\mathcal D.
\]
Using \eqref{eq:global-reconstruction}, the surface integral containing $D_x^2f$ is bounded by
\[
C\|\nabla_xf\|_{L^\infty}\|\nabla_xf\|_{L^2}\|D_x^2f\|_{L^2} \le Ce^{1/2}\mathcal D.
\]

For the surface transport term, we integrate by parts and use
\[
 \|(c-\sigma\Delta_x)^{-1}\beta\|_{L^2} +\|\nabla_x(c-\sigma\Delta_x)^{-1}\beta\|_{L^2} \le C\|\beta\|_{L^2}.
\]
The resulting term is bounded by
\[
 C\bigl(\|U_x|_{z=0}\|_{L^\infty} +\|\diver_x(U_x|_{z=0})\|_{L^\infty}\bigr)\|\beta\|_{L^2}^2 \le C\mathcal W_d(U)e.
\]
For the term involving $R_x$, \eqref{eq:global-reconstruction} gives 
\begin{equation*}
\|\nabla_x(c-\sigma\Delta_x)f\|_{L^2}+\|f_t\|_{L^2}\le C\mathcal D^{1/2}, \qquad \|\nabla_xf\|_{L^\infty}\le C\min(e^{1/2},\mathcal D^{1/2}).    
\end{equation*}
Together with \eqref{eq:global-curvature-homogeneous}, these bounds control each term in \eqref{eq:global-normal-surface-velocity} and give $\|R_x\|_{L^2}\le C\mathcal D$.
Since
\[
\|\nabla_x(c-\sigma\Delta_x)^{-1}\beta\|_{L^2}\le C\|\beta\|_{L^2}\le Ce^{1/2}, \qquad e\mathcal D\le e^{1/2}\mathcal D
\]
when $e\le 1$, the estimate \eqref{eq:global-normal-estimate-equation} follows.
\end{proof}

\subsubsection{Tangential energy and the corrected pressure}

For $|\alpha|=m$, consider the tangential energy
\begin{equation}\label{eq:global-tangential-energy}
 L_\alpha=\frac12\int_{\mathcal S}(-\gamma'(z))|\partial_x^\alpha\eta|^2\dd x\dd z +\frac c2\|\partial_x^\alpha f\|_{L_x^2}^2 +\frac\sigma2\|\nabla_x\partial_x^\alpha f\|_{L_x^2}^2.
\end{equation}
Differentiating \eqref{eq:global-coordinate-identities} produces terms of the form $(\partial_x^\alpha\nabla_{x,z}\eta)U$.
These involve $m+1$ horizontal derivatives of $\eta$ or $m$ horizontal derivatives of $\eta_z$, neither of which is controlled by $\|\eta\|_{X_m}$.
We use a corrected pressure and a divergence-free corrected flux to derive an energy identity without this loss.

Fix a smooth radial cutoff $\chi\in C_c^\infty(\R)$ satisfying $\chi=1$ on $[0,1/8]$ and $\chi=0$ on $[1/4,\infty)$.
The horizontal Bony paraproduct $T_bh$ has Fourier multiplier $\chi(|\ell|/|k|)$, where $\ell,k$ are the horizontal frequencies of the coefficients $b$ and $h$ respectively.
We set this symbol to zero at $k=0$.
For each multi-index $\alpha\in\mathbb N_0^d$ with $|\alpha|=m$, define the paraproduct sum and commutator remainder by
\begin{equation}\label{eq:global-remainders}
\begin{aligned}
 S_\alpha(\eta,F) &=\sum_{0<\beta\le\alpha}\binom{\alpha}{\beta} T_{\partial_x^{\alpha-\beta}F}\partial_x^\beta\eta,\\
 R_{\alpha,i}(\eta,F) &=\sum_{0<\beta\le\alpha}\binom{\alpha}{\beta}
 (\partial_i\partial_x^\beta\eta)\partial_x^{\alpha-\beta}F -\partial_i S_\alpha(\eta,F), \qquad i=1,\ldots,d,z,
\end{aligned}
\end{equation}
acting componentwise when $F$ is vector-valued.
The good pressure unknown at order $\alpha$ is
\begin{equation}\label{eq:global-good-pressure}
P_\alpha=\partial_x^\alpha\Pi+S_\alpha(\eta,U_{d+1}).
\end{equation}

We introduce vector fields
\[
 F_\alpha=(J_\eta \partial_x^\alpha U_x,\partial_x^\alpha U_{d+1}-\nabla_x\eta\cdot \partial_x^\alpha U_x),\qquad
 T_\alpha=(\partial_x^\alpha U_x+\nabla_x\eta\,\partial_x^\alpha U_{d+1},J_\eta \partial_x^\alpha U_{d+1}),
\]
whose pointwise inner product is $F_\alpha\cdot T_\alpha=J_\eta|\partial_x^\alpha U|^2$, and also the corrected flux
\begin{equation}\label{eq:global-corrected-flux}
G_\alpha=F_\alpha+R_\alpha^F,\qquad R_\alpha^F= \Big((R_{\alpha,z}(\eta,U_i))_{i=1}^d, -\sum_{i=1}^d R_{\alpha,i}(\eta,U_i)\Big).
\end{equation}
Applying $\partial_x^\alpha$ to the coordinate relations \eqref{eq:global-coordinate-identities} and using the definitions \eqref{eq:global-remainders}--\eqref{eq:global-good-pressure} gives
\begin{equation}\label{eq:global-corrected-identities}
 \begin{gathered}
 \diver G_\alpha=0,\qquad
 T_\alpha=-\nabla P_\alpha+\gamma'(z)\partial_x^\alpha\eta e_{d+1} -R_\alpha(\eta,U_{d+1}),\\
G_{\alpha,d+1}=\partial_t\partial_x^\alpha\eta+\diver_xS_\alpha(\eta,U_x), \qquad G_{\alpha,d+1}|_{z=-1}=0,\\
P_\alpha|_{z=0}=(c-\sigma\Delta_x)\partial_x^\alpha f +\sigma\partial_x^\alpha \mathcal N(f)+S_\alpha(f,U_{d+1}|_{z=0}).
 \end{gathered}
\end{equation}
Indeed, the definitions give
\[
 G_\alpha-\partial_x^\alpha\mathsf v =\bigl(-\partial_zS_\alpha(\eta,U_x),\, \diver_xS_\alpha(\eta,U_x)\bigr).
\]
The right-hand side has zero divergence. 
Since $\diver\mathsf v=0$, we get $\diver G_\alpha=0$.
The bottom condition $G_{\alpha,d+1}|_{z=-1}=0$ follows from $\eta=U_{d+1}=0$ on $\{z=-1\}$.

For smooth solutions, taking the $L^2(\mathcal S)$ inner product of $T_\alpha$ with $F_\alpha=G_\alpha-R_\alpha^F$ and integrating by parts using \eqref{eq:global-corrected-identities} gives
\begin{equation}\label{eq:global-tangential-identity}
 \begin{aligned}
 &L_\alpha'+\int_{\mathcal S}J_\eta|\partial_x^\alpha U|^2\dd x\dd z
 =-\int_{\mathcal S}(-\gamma'(z))\partial_x^\alpha\eta
       \diver_xS_\alpha(\eta,U_x)\dd x\dd z\\
 -&\langle(c-\sigma\Delta_x)\partial_x^\alpha f,\diver_xS_\alpha(f,U_x|_{z=0})\rangle-\sigma\langle\partial_x^\alpha\mathcal N(f),G_{\alpha,d+1}|_{z=0}\rangle
  -\langle S_\alpha(f,U_{d+1}|_{z=0}),G_{\alpha,d+1}|_{z=0}\rangle\\
 -&\int_{\mathcal S}R_\alpha(\eta,U_{d+1})\cdot G_\alpha\dd x\dd z -\int_{\mathcal S}T_\alpha\cdot R_\alpha^F\dd x\dd z.
 \end{aligned}
\end{equation}

The following estimates control the transport and commutator terms in \eqref{eq:global-corrected-identities}, as well as the boundary trace of $S_\alpha$.

\begin{lemma}[Commutator and transport bounds]\label{lem:global-tangential-bounds}
For every integer $m>2d$ and $|\alpha|=m$, one has
\begin{equation}\label{eq:global-remainder-bound}
 \|R_\alpha(\eta,F)\|_{L^2(\mathcal S)} \le C\|\eta\|_{X_m}\|F\|_{X_m}.
\end{equation}
Moreover, for real-valued fields, the bulk and surface transport terms satisfy
\begin{align}
 \left|\int_{\mathcal S}(-\gamma'(z))\partial_x^\alpha\eta\diver_x S_\alpha(\eta,U_x)\dd x\dd z\right|
     &\le C\mathcal W_d(U_x)\|\eta\|_{X_m}^2,  \label{eq:global-bulk-transport}\\
 \left|\langle (c-\sigma\Delta_x)\partial_x^\alpha f,\diver_x S_\alpha(f,U_x|_{z=0})\rangle\right|  &\le C\mathcal W_d(U_x)\|f\|_{H^{m+1}}^2, \label{eq:global-surface-transport}
\end{align}
and the surface paraproduct trace obeys
\begin{equation}\label{eq:global-S-trace}
 \|S_\alpha(f,F|_{z=0})\|_{H^{1/2}(\R^d)}  \le C\mathcal W_d(F)\|f\|_{H^{m+1/2}(\R^d)}.
\end{equation}
\end{lemma}

\begin{proof}
We first estimate the commutator remainder using the extension operator $E_z$ of Lemma~\ref{lem:global-embeddings}.
The identity $\partial_iT_bh=T_b\partial_ih+T_{\partial_ib}h$ shows that each summand of $R_{\alpha,i}(\eta,F)$ has the form
\[
(\partial_i\partial_x^\beta\eta)\partial_x^{\alpha-\beta}F
 -T_{\partial_x^{\alpha-\beta}F}(\partial_i\partial_x^\beta\eta)-T_{\partial_i\partial_x^{\alpha-\beta}F}\partial_x^\beta\eta.
\]
In the first difference, the horizontal frequency of $\eta$ is bounded by a constant times that of $F$.
In the last paraproduct, $F$ has lower horizontal frequency.
Bony's decomposition and Young's inequality therefore give
\[
 \|R_{\alpha,i}(\eta,F)\|_{L^2(\mathcal S)} \le C\Bigl(\|\widehat{\partial_iE_z\eta}\|_{L^1}\||D_x|^mE_zF\|_{L^2}+\|\widehat{\partial_iE_zF}\|_{L^1} \||D_x|^mE_z\eta\|_{L^2}\Bigr).
\]
The norms on the right are over $\R^d\times\R$.
The extension bounds and \eqref{eq:global-Fourier-embeddings} prove \eqref{eq:global-remainder-bound}, including the case $i=z$.

For the transport terms, we separate the leading paraproduct by writing
\[
 S_\alpha(h,a)=T_a\partial_x^\alpha h+Q_\alpha(h,a).
\]
Every term of $Q_\alpha$ has a coefficient $\partial_x^{\alpha-\beta}a$ with $r=|\alpha-\beta|\ge1$.
For a horizontal coefficient $a=a(x)$, set
\[
 M(a)=\int_{\R^d}(1+|\xi|)|\widehat a(\xi)|\,\dd\xi.
\]
For $r=1$, the coefficient satisfies
$\|\partial_x^{\alpha-\beta}a\|_{L^\infty}\le CM(a)$.
Whereas for $r>1$, it satisfies
$\|\partial_x^{\alpha-\beta}a\|_{C_*^{1-r}}\le CM(a)$.
The paraproduct estimates in \cite[Theorems~2.6(i) and~2.9(ii)]{AlazardBurqZuily2014}, applied respectively in these two cases, give
\[
 \|\diver_xQ_\alpha(h,a)\|_{H^s} \le CM(a)\|h\|_{H^{m+s}},\qquad s=0,1.
\]
Indeed, the paraproduct with this coefficient has order $r-1$, and its divergence has total order $m$ on $h$.
The cited estimates use a cutoff at low frequencies.
For the homogeneous cutoff used here, the additional terms have bounded input and output frequencies.
They satisfy the same bounds by Young's inequality and the definition of $M(a)$.

The adjoint estimate \cite[Theorem~2.6(iii)]{AlazardBurqZuily2014} shows that $\diver_xT_a+(\diver_xT_a)^*$ has order zero, with $L^2$ operator norm at most $CM(a)$.
Thus
\[
 |\langle w,\diver_xT_aw\rangle|\le CM(a)\|w\|_{L^2}^2.
\]
Moreover,
\[
\partial_j\diver_xT_aw=\diver_xT_a\partial_jw+\diver_xT_{\partial_ja}w, \qquad
 \|\diver_xT_{\partial_ja}w\|_{L^2}
 \le CM(a)\|w\|_{H^1}.
\]
Applying the preceding $L^2$ estimate to $w$ and $\partial_jw$, and integrating by parts in the surface Laplacian, gives
\[
 \bigl|\langle(c-\sigma\Delta_x)w, \diver_xT_aw\rangle\bigr|
 \le CM(a)\|w\|_{H^1}^2.
\]
Combining these estimates with the bound for $\diver_xQ_\alpha$, we obtain the bulk and surface transport bounds.
Fourier inversion in the normal variable gives, for any extension $\widetilde{U_x}$ of $U_x$ to $\R^d\times\R$,
\[
 M(U_x(\cdot,z)) \le C\int_{\R^d}\int_\R(1+|\xi|) |\widehat{\widetilde{U_x}}(\xi,\zeta)|\,\dd\zeta\dd\xi.
\]
Taking the supremum in $z$ and then the infimum over extensions yields $\sup_{z\in\overline I}M(U_x(\cdot,z))\le C\mathcal W_d(U_x)$.
For the bulk term, the $s=0$ bounds applied at each $z$ therefore give
\begin{align*}
 &\left|\int_{\mathcal S}(-\gamma'(z))\partial_x^\alpha\eta \diver_xS_\alpha(\eta,U_x)\,\dd x\dd z\right| \le C\int_I(-\gamma'(z))M(U_x(\cdot,z)) \|\eta(\cdot,z)\|_{H_x^m}^2\,\dd z\\
 & \le C\mathcal W_d(U_x)  \bigl(\|\eta\|_{L^2(\mathcal S)}^2 +\||D_x|^m\eta\|_{L^2(\mathcal S)}^2\bigr)
 \le C\mathcal W_d(U_x)\|\eta\|_{X_m}^2.
\end{align*}
For the surface term, we use the leading-paraproduct estimate together with the $s=1$ bound for $\diver_xQ_\alpha$:
\[
\begin{aligned}
 &\bigl|\langle(c-\sigma\Delta_x)\partial_x^\alpha f, \diver_xS_\alpha(f,U_x|_{z=0})\rangle\bigr|\\
 \le& CM(U_x|_{z=0})\|\partial_x^\alpha f\|_{H^1}^2 +\|(c-\sigma\Delta_x)\partial_x^\alpha f\|_{H^{-1}} \|\diver_xQ_\alpha(f,U_x|_{z=0})\|_{H^1}\\
 \le& CM(U_x|_{z=0})\|f\|_{H^{m+1}}^2 \le C\mathcal W_d(U_x)\|f\|_{H^{m+1}}^2.
\end{aligned}
\]

Finally, paraproduct continuity gives
\[
 \|S_\alpha(f,a)\|_{H^{1/2}}
 \le CM(a)\|f\|_{H^{m+1/2}}.
\]
Taking $a=F|_{z=0}$ and using $M(F|_{z=0})\le C\mathcal W_d(F)$ proves \eqref{eq:global-S-trace}.
\end{proof}

\subsubsection{Completion of the energy estimate}

To estimate the boundary terms in \eqref{eq:global-tangential-identity} and recover the dissipative gain of $3/2$ derivatives on $f$, we use trace inequalities on the strip $\mathcal S$.
The harmonic extension of a trace $\psi$ on $\{z=0\}$, with homogeneous Neumann condition at $\{z=-1\}$, has Dirichlet energy $\int_{\R^d}|\xi|\tanh|\xi|\,|\widehat\psi(\xi)|^2\dd\xi$.
Consequently, for any $P\in H^1(\mathcal S)$ and any divergence-free vector field $G\in L^2(\mathcal S)$ with $G_{d+1}|_{-1}=0$, Green's identity gives the dual trace estimates
\begin{equation}\label{eq:global-strip-traces}
\begin{aligned}
 &\|(|D_x|\tanh|D_x|)^{1/2}P|_{z=0}\|_{L^2(\R^d)}\le\|\nabla P\|_{L^2(\mathcal S)}, \\
& \|(|D_x|\tanh|D_x|)^{-1/2}G_{d+1}|_{z=0}\|_{L^2(\R^d)}\le\|G\|_{L^2(\mathcal S)}.
 \end{aligned}
\end{equation}
Since $r\tanh r\asymp r$ for $r\ge1$ and $r\tanh r\asymp r^2$ for $0<r\le1$, these estimates give $\dot H^{\pm1/2}$ trace bounds at high horizontal frequencies and an $L^2\to \dot H^{-1}$ normal-trace bound at low horizontal frequencies.

\begin{lemma}\label{lem:global-surface-estimates}
Whenever the total energy $e$ is sufficiently small, the free surface satisfies
\begin{equation}\label{eq:global-surface-high}
 \|f\|_{H^{m+5/2}(\R^d)}\le C(e^{1/2}+\mathcal D^{1/2}),
 \qquad \|f_t\|_{H^{m-1/2}(\R^d)}\le C\mathcal D^{1/2}.
\end{equation}
\end{lemma}

\begin{proof}
For any $|\alpha|=m$, the identities \eqref{eq:global-corrected-identities} and the remainder bound \eqref{eq:global-remainder-bound} imply 
\[
\|\nabla P_\alpha\|_{L^2(\mathcal S)}\le C(e^{1/2}+\mathcal D^{1/2}), \qquad \|G_\alpha\|_{L^2(\mathcal S)}\le C\mathcal D^{1/2}.
\]
We apply the first trace bound in \eqref{eq:global-strip-traces} to $P_\alpha|_{z=0}$ on $\{|\xi|\ge1\}$.
In the boundary condition for $P_\alpha$ in \eqref{eq:global-corrected-identities}, one can use \eqref{eq:global-S-trace} to bound the paraproduct by $Ce^{1/2}\mathcal D^{1/2}$, and use \eqref{eq:global-curvature-tame} with $q=m+1/2$ to bound the curvature term.
Adding the low-frequency bound $\|f\|_{L^2}\le e^{1/2}$ yields
\[
 \|f\|_{H^{m+5/2}} \le C(e^{1/2}+\mathcal D^{1/2}) +C\|\nabla_xf\|_{L^\infty}^2\|f\|_{H^{m+5/2}}.
\]
For sufficiently small $e$, we absorb the last term and obtain the first estimate in \eqref{eq:global-surface-high}.

For $f_t$, the second inequality in \eqref{eq:global-strip-traces} and $r\tanh r\le\langle r\rangle$ give
\[
 \|G_{\alpha,d+1}|_{z=0}\|_{H^{-1/2}}
 \le\|(|D_x|\tanh|D_x|)^{-1/2}G_{\alpha,d+1}|_{z=0}\|_{L^2}
 \le\|G_\alpha\|_{L^2(\mathcal S)}.
\]
Using $\partial_x^\alpha f_t=G_{\alpha,d+1}|_{z=0}-\diver_xS_\alpha(f,U_x|_{z=0})$ and \eqref{eq:global-S-trace}, we obtain
\[
\begin{aligned}
 &\|\partial_x^\alpha f_t\|_{H^{-1/2}}
 \le\|G_\alpha\|_{L^2(\mathcal S)} +C\|S_\alpha(f,U_x|_{z=0})\|_{H^{1/2}}\\
 \le& C\mathcal D^{1/2}  +C\mathcal W_d(U_x)\|f\|_{H^{m+1/2}}\le C(1+e^{1/2})\mathcal D^{1/2}
 \le C\mathcal D^{1/2}.
\end{aligned}
\]
Finally, one can use \eqref{eq:global-reconstruction} to control the low frequencies.
Hence,
\[
 \|f_t\|_{H^{m-1/2}} \le C\Big(\|f_t\|_{L^2} +\sum_{|\alpha|=m}\|\partial_x^\alpha f_t\|_{H^{-1/2}}\Big) \le C\mathcal D^{1/2}.
\]
This proves the second estimate in \eqref{eq:global-surface-high}.
\end{proof}

Let $E_{\rm rel}$ denote the relative energy in \eqref{eq:global-relative-energy}.
We define the total energy by
\begin{equation}\label{eq:global-total-energy}
 \mathcal E=E_{\rm rel}+\sum_{|\alpha|=m}L_\alpha+\mathcal Q_h.
\end{equation}
We now show that $\mathcal{E}$ is the energy that we want, and finish the proof of Proposition~\ref{prop:global-energy}.
\begin{proof}[Proof of Proposition~\ref{prop:global-energy}]
Using the results
\[
-\gamma'(z)\ge c_*>0, \qquad c>0, \qquad \int_{\R^d}(\sqrt{1+|\nabla f|^2}-1)\dd x\asymp\|\nabla f\|_{L^2}^2
\]
when $|\nabla f|$ is small, we get $\mathcal E\asymp e$ by Lemma~\ref{lem:global-normal-comparison}.
The same lemma shows that the sum of the dissipative integrals in \eqref{eq:global-relative-energy}, \eqref{eq:global-normal-estimate-equation}, and \eqref{eq:global-tangential-identity} is at least $2\kappa\mathcal D$ for some $\kappa>0$.

For the bulk remainders in \eqref{eq:global-tangential-identity}, one can use \eqref{eq:global-remainder-bound} and the definitions of $T_\alpha$ and $G_\alpha$ to get
\begin{align*}
 &\|R_\alpha(\eta,U_{d+1})\|_{L^2}+\|R_\alpha^F\|_{L^2}\le C\|\eta\|_{X_m}\|U\|_{X_m}
 \le Ce^{1/2}\mathcal D^{1/2},\\
 &\|T_\alpha\|_{L^2}+\|G_\alpha\|_{L^2}\le C\bigl(1+\|\nabla_{x,z}\eta\|_{L^\infty}+\|\eta\|_{X_m}\bigr)\|U\|_{X_m}
 \le C\mathcal D^{1/2}.
\end{align*}

Consequently,
\begin{align*}
 &\left|\int_{\mathcal S}R_\alpha(\eta,U_{d+1})\cdot G_\alpha\,\dd x\dd z\right|+\left|\int_{\mathcal S}T_\alpha\cdot R_\alpha^F\,\dd x\dd z\right|\\
 \le&\|R_\alpha(\eta,U_{d+1})\|_{L^2}\|G_\alpha\|_{L^2} +\|T_\alpha\|_{L^2}\|R_\alpha^F\|_{L^2} \le Ce^{1/2}\mathcal D.
\end{align*}
The transport estimates \eqref{eq:global-bulk-transport}--\eqref{eq:global-surface-transport} yield
\[
\begin{aligned}
 &\left|\int_{\mathcal S}(-\gamma'(z))\partial_x^\alpha\eta \diver_xS_\alpha(\eta,U_x)\,\dd x\dd z\right|+\bigl|\langle(c-\sigma\Delta_x)\partial_x^\alpha f,\diver_xS_\alpha(f,U_x|_{z=0})\rangle\bigr|\\
 &\le C\mathcal W_d(U_x) \bigl(\|\eta\|_{X_m}^2+\|f\|_{H^{m+1}}^2\bigr) \le C\mathcal W_d(U)e.
\end{aligned}
\]
For the remaining noncapillary boundary term, Sobolev duality, \eqref{eq:global-strip-traces}, and \eqref{eq:global-S-trace} give
\[
\begin{aligned}
 &\bigl|\langle S_\alpha(f,U_{d+1}|_{z=0}),G_{\alpha,d+1}|_{z=0}\rangle\bigr|\le\|S_\alpha(f,U_{d+1}|_{z=0})\|_{H^{1/2}}\|G_{\alpha,d+1}|_{z=0}\|_{H^{-1/2}}\\
 \le& C\mathcal W_d(U_{d+1})\|f\|_{H^{m+1/2}}\|G_\alpha\|_{L^2(\mathcal S)}\le Ce^{1/2}\mathcal W_d(U)\mathcal D^{1/2}
 \le Ce^{1/2}\mathcal D.
\end{aligned}
\]
The curvature term in \eqref{eq:global-tangential-identity} is
$-\sigma\langle\partial_x^\alpha\mathcal N(f),G_{\alpha,d+1}|_{z=0}\rangle$.
We estimate it by duality between $H^{1/2}$ and $H^{-1/2}$.
Since $r\tanh r\le\langle r\rangle$, \eqref{eq:global-strip-traces} gives
\[
 \|G_{\alpha,d+1}|_{z=0}\|_{H^{-1/2}}\le C\|G_\alpha\|_{L^2} \le C\mathcal D^{1/2}.
\]
By \eqref{eq:global-reconstruction}, $\|\nabla_xf\|_{L^\infty}^2\le Ce$ and $\|\nabla_xf\|_{L^\infty}^2\le C\sqrt{e\mathcal D}$.
It follows that
\[
\|\nabla_xf\|_{L^\infty}^2(e^{1/2}+\mathcal D^{1/2})\le Ce\,\mathcal D^{1/2}.
\]
Combining the tame curvature estimate \eqref{eq:global-curvature-tame} with the $H^{m+5/2}$ bound for $f$ in \eqref{eq:global-surface-high} gives
\[
 \|\partial_x^\alpha\mathcal N(f)\|_{H^{1/2}} \le C\|\nabla f\|_{L^\infty}^2(e^{1/2}+\mathcal D^{1/2}) \le Ce\,\mathcal D^{1/2}.
\]
Therefore,
\[
\sigma\abs{\langle\partial_x^\alpha\mathcal N(f),G_{\alpha,d+1}|_{z=0}\rangle}
 \le C\sigma e\mathcal D.
\]
Since $\sigma$ is fixed, this is bounded by $Ce\mathcal D$.
Thus \eqref{eq:global-tangential-identity} yields
\begin{equation}\label{eq:global-tangential-estimate}
 \frac{\dd}{\dd t}\sum_{|\alpha|=m}L_\alpha
 +\sum_{|\alpha|=m}\int_{\mathcal S}
     J_\eta|\partial_x^\alpha U|^2\,\dd x\dd z
 \le Ce^{1/2}\mathcal D+C\mathcal W_d(U)e.
\end{equation}

Adding \eqref{eq:global-relative-energy},
\eqref{eq:global-normal-estimate-equation}, and
\eqref{eq:global-tangential-estimate},
and using \eqref{eq:global-normal-dissipation-comparison}
together with $e\mathcal D\le e^{1/2}\mathcal D$ for $e\le1$,
we obtain \eqref{eq:global-energy-inequality}
for smooth solutions.

To extend the estimate to local strong solutions, we use the local regularity and the coordinate transition to obtain $\eta\in C_TH^m$, $U\in L_T^2H^m$, and $\eta_t=U_{d+1}-U_x\cdot\nabla_x\eta\in L_T^2H^{m-1}$.
In particular, $\mathcal E$ is continuous in time.
For $|\alpha|=m$, the pressure relation $\partial_i\Pi=-U_i-(\partial_i\eta)U_{d+1}$, with any $i$ such that $\alpha_i>0$, gives $\partial_x^\alpha\Pi\in L_T^2L^2$.
Thus \eqref{eq:global-good-pressure} and \eqref{eq:global-corrected-identities} imply $P_\alpha\in L_T^2H^1$ and $G_\alpha\in L_T^2L^2$.
The divergence-free condition and \eqref{eq:global-strip-traces} justify Green's identity for these fields, with the boundary terms understood by Sobolev duality.
To justify differentiating the bulk energy term in
\eqref{eq:global-tangential-energy}, we write
\[
 \partial_t\partial_x^\alpha\eta +\diver_xT_{U_x}\partial_x^\alpha\eta =G_{\alpha,d+1}-\diver_xQ_\alpha(\eta,U_x),
\]
Here $Q_\alpha(\eta,U_x)=S_\alpha(\eta,U_x)-T_{U_x}\partial_x^\alpha\eta$, as in the proof of Lemma~\ref{lem:global-tangential-bounds}.
The right-hand side belongs to $L_T^1L^2$ by that lemma.
Horizontal Friedrichs regularization therefore justifies the bulk energy identity.
The bounds $\partial_x^\alpha f\in L_T^2H^{5/2}$ and $\partial_x^\alpha f_t\in L_T^2H^{-1/2}$ justify differentiating the surface $H^1$ energy and yield the integrated tangential identity.

On each compact local existence interval, $u\in L_t^1W^{1,\infty}$ generates a volume-preserving flow carrying the initial fluid domain onto $\Om_{f(t)}$.
For a smooth real-valued scalar field $w$,
\begin{equation}\label{eq:global-material-transport}
 \frac12\frac{\dd}{\dd t}\int_{\Om_{f(t)}}|w|^2\dd x\dd y =\int_{\Om_{f(t)}}wD_tw\dd x\dd y.
\end{equation}
Spatial regularization and time averaging, or pullback along this flow, extend the identity to the regularity used here.

For the normal energy, the right-hand side of \eqref{eq:global-normal-laplacian-transport} belongs to $L_T^1L^2(\Om_f)$, and the weight $(-\gamma'(z_\rho))^{-1}$ is transported.
The weighted form of \eqref{eq:global-material-transport} gives the bulk calculation in \eqref{eq:global-normal-identity}.
Using \eqref{eq:global-normal-beta-transport}, we differentiate the surface energy $\frac12\langle(c-\sigma\Delta_x)^{-1}\beta,\beta\rangle$ and recover the surface terms in the integrated form of \eqref{eq:global-normal-identity}.
To justify the spatial integrations by parts, we insert horizontal cutoffs $\chi_R$ and let $R\to\infty$.
The terms containing derivatives of $\chi_R$ vanish by the Cauchy--Schwarz inequality and the $L^2$ bounds, using $\|\nabla_x\chi_R\|_\infty\le C/R$.
Combining the integrated forms of
\eqref{eq:global-relative-energy},
\eqref{eq:global-normal-estimate-equation}, and \eqref{eq:global-tangential-estimate},
and using \eqref{eq:global-normal-dissipation-comparison},
gives \eqref{eq:global-integrated-energy}.
\end{proof}

\subsection{Finite-time density regularity}\label{sec:global-density-regularity}

To apply the continuation criterion, it remains to control the full Sobolev norm of the density on finite time intervals.

\begin{proposition}[Finite-time density regularity]\label{prop:global-strong-energy}
Suppose a local strong solution in \eqref{eq:global-local-class} admits a regular density coordinate map on $[0,T]$, with $\sup_{t\le T}e(t)$ sufficiently small.
Then, with $C_R$ independent of $T$,
\begin{equation}\label{eq:global-finite-density}
\begin{aligned}
 \sup_{t\le T}\|\widetilde g(t)\|_{H^m(\mathcal S)} \le& C_R\Big(\|\widetilde g_0\|_{H^m} +T\sup_{t\le T}e(t)^{1/2}+\sqrt T \Big(\int_0^T \mathcal D(t)\dd t\Big)^{1/2}\Big)\\
 &\cdot\exp\Big(C_R\big(T+\sqrt T \big(\int_0^T \mathcal D(t)\dd t\big)^{1/2}\big)\Big).
 \end{aligned}
\end{equation}
\end{proposition}

We estimate the physical perturbation $g$ on the moving domain, where the material transport identity \eqref{eq:global-material-transport} eliminates the boundary flux in the commuted $H^m$ estimate.
The tame elliptic bound below controls the velocity by the density norm and the surface regularity given by Proposition~\ref{prop:global-energy} and Lemma~\ref{lem:global-surface-estimates}.
Gronwall's inequality then gives a density bound on each finite interval.

The velocity estimate \eqref{eq:velocity-single}, applied with $s=m$, and the Sobolev equivalence \eqref{eq:sobolev-pullback} give
\begin{equation}\label{eq:global-local-tame}
 \|u\|_{H^m(\Om_f)} \le C_R\bigl(\|g\|_{H^m(\Om_f)}+\|f\|_{H^{m+5/2}}\bigr),
 \qquad \|g\|_{H^m(\Om_f)}\asymp_R\|\widetilde g\|_{H^m(\mathcal S)}.
\end{equation}

\begin{proof}[Proof of Proposition~\ref{prop:global-strong-energy}]
We work in physical coordinates and introduce
\[
 G_m(t)=\Big(\sum_{|\alpha|\le m}         \|\partial_{x,y}^\alpha g(t)\|_{L^2(\Om_{f(t)})}^2\Big)^{1/2}.
\]
We commute physical derivatives with $D_tg=-\gamma'(y)u_{d+1}$ and apply \eqref{eq:global-material-transport} to each derivative.
The tame product estimate then gives, in the integrated sense,
\[
 G_m' \le C_R\|\nabla u\|_{L^\infty} G_m +C_R\|u\|_{H^m(\Om_f)}\|\nabla g\|_{L^\infty} +C_R\|\gamma'(y)u_{d+1}\|_{H^m(\Om_f)}.
\]
The calculation is justified by regularizing the transport equation, with constants uniform for bounded graph geometry.
Using the identity $g\circ\Phi_\eta=\gamma(z)-\gamma(z+\eta)$ and Lemma~\ref{lem:global-embeddings}, we obtain
\[
 \|\nabla g\|_{L^\infty}\le C_R\sup_{t\le T}e(t)^{1/2}, \qquad \|\nabla u\|_{L^\infty}\le C_R\mathcal D(t)^{1/2}.
\]
The bounded profile derivatives, \eqref{eq:global-local-tame}, and \eqref{eq:global-surface-high} also give
\[
 \|\gamma'(y)u_{d+1}\|_{H^m(\Om_f)} \le C_R\|u\|_{H^m(\Om_f)} \le C_R(G_m+\sup_{t\le T}e(t)^{1/2}+\mathcal D(t)^{1/2}).
\]
Since $\sup_{t\le T}e(t)^{1/2}$ is small, we get
\[
 G_m'\le C_R\big(1+\mathcal D(t)^{1/2}\big)G_m+C_R\big(\sup_{t\le T}e(t)^{1/2}+\mathcal D(t)^{1/2}\big).
\]
Using Gronwall's inequality, $\int_0^T \mathcal D(t)^{1/2}\dd t\le\sqrt T (\int_0^T \mathcal D(t)\dd t)^{1/2}$, and the norm equivalence in \eqref{eq:global-local-tame}, we prove the estimate \eqref{eq:global-finite-density}.
\end{proof}

\subsection{Global well-posedness}\label{sec:global-continuation}

We prove Theorem~\ref{thm:global} by continuing the local solution of Corollary~\ref{prop:global-local}, using the velocity, energy, and density bounds above.

\begin{proof}[Proof of Theorem~\ref{thm:global}]
Corollary~\ref{prop:global-local} and Lemma~\ref{lem:global-initial-chart} give a unique local strong solution with a regular initial density coordinate map and
\[
 \mathcal E(0)\le C\varepsilon^2, \qquad \|\eta_0\|_{H^m}+\mathcal B_{\rm low}(\eta_0)\le C\varepsilon.
\]
Let $T_*$ be its maximal existence time.
Choosing $\delta_*>0$ small so that $e\le\delta_*$ ensures the hypotheses of Propositions~\ref{prop:global-source-closure}, \ref{prop:global-energy}, and \ref{prop:global-strong-energy}.
We also have $1/2<J_\eta<3/2$, $\|f\|_{L^\infty}<1/2$, and $Ce^{1/2}\le\kappa/2$.

On any interval $[0,T]\subset[0,T_*)$ where the density coordinate map is regular and $\sup_{t\le T}e(t)\le\delta_*$, Proposition~\ref{prop:global-source-closure} and \eqref{eq:global-linear-Y-Wiener} give
\[
 \int_0^T\bigl(\mathcal Y(\mathsf v(t)) +\mathcal Y(U(t))\bigr)\dd t\le C\varepsilon, \qquad \int_0^T\mathcal W_d(U(t))\dd t\le C\varepsilon.
\]
By Proposition~\ref{prop:global-energy}, absorbing $Ce^{1/2}\mathcal D$ in \eqref{eq:global-integrated-energy} and using $\mathcal E\asymp e$, Gronwall's inequality yields
\[
 \sup_{t\le T}e(t)+\int_0^T \mathcal D(t)\dd t \le C\varepsilon^2\exp(C\varepsilon),
\]
uniformly in $T$.
We choose $\varepsilon_*$ so that the right-hand side is below $\delta_*/2$.
This strictly improves the energy bound and preserves the invertibility of the density coordinate map.
By continuity, these estimates can be extended to every $T<T_*$.

If $T_*<\infty$, Proposition~\ref{prop:global-strong-energy} and \eqref{eq:global-finite-density} give $\sup_{t<T_*}\|\widetilde g(t)\|_{H^m}<\infty$, while $\|f(t)\|_{H^{m+1}}$ remains small and the depth stays positive.
Using Corollary~\ref{prop:global-local}, one can extend the solution beyond $T_*$, which contradicts the assumption $T_*<\infty$.
Thus $T_*=\infty$. 
Letting $T\to\infty$, we obtain the estimate \eqref{eq:global-estimate}.
The identity $\rho_y\circ\Phi_\eta=\gamma'(z)/J_\eta<0$ and the preserved boundary values give $c\le\rho\le\gamma(-1)$ and strict vertical monotonicity.

For each $T<\infty$, \eqref{eq:global-estimate} and
\eqref{eq:global-finite-density} bound
$\|f\|_{C_TH^{m+1}}+\|\widetilde g\|_{C_TH^m}$
uniformly over the compatible small-data class,
with a bound depending on $T$.
The smallness of $\|f\|_{C_TH^{m+1}}$ also ensures
that the fluid depth satisfies $1+f\ge1/2$.
Theorem~\ref{thm:local-target}, specialized in Corollary~\ref{prop:global-local}, gives \eqref{eq:global-local-class} and uniqueness there.
The uniform local lifespan and the continuous-dependence conclusion of Theorem~\ref{thm:local-target}, applied on finitely many subintervals, give continuous dependence in the stated topology.
\end{proof}

\section*{Acknowledgments}
Lizhe Wan is partially supported by the National Natural Science Foundation of China under Grant No.~12671283.
Jiaqi Yang is supported by the National Natural Science Foundation of China under Grant No.~12471225 and by the Natural Science Basic Research Program of Shaanxi under Program No.~2026JC-YXQN-021. 

\section*{Declaration of generative AI and AI-assisted technologies}
 During the preparation of this work the authors used ChatGPT 6.1 for editorial feedback and language editing. 
 After using this tool, the authors reviewed and edited the content as needed.
 They take full responsibility for the content of the published article.

\appendix
\section{Auxiliary estimates and the local theory without surface tension}\label{app:latocca-nguyen}

We collect the elliptic and transport estimates and the local theory without surface tension used in the main text.
The statements below give the relevant parts of the cited results or their consequences on bounded sets, in our notation.
The Dirichlet--Neumann estimates are due to Alazard--Burq--Zuily \cite{AlazardBurqZuily2014} and Nguyen--Pausader \cite{NguyenPausader2020}.
For transport estimates, background multipliers, and the coupled boundary--bulk theory, we use Latocca--Nguyen \cite{LatoccaNguyen2025}.

Several estimates in the cited references are stated on $\R^d$; the periodic case $\cO=\T^d$ follows by the Fourier-localization and extension arguments described in the main text.
The periodic local well-posedness theorem is stated explicitly in \cite[Remark~1.6]{LatoccaNguyen2025}.
Throughout the appendix, finite-depth domains have a common positive depth bound and a fixed bottom.
We use the energy normalization of Section~\ref{sec:elliptic} for harmonic extensions.

\subsection{Flattening and elliptic operators}

The flattening estimate gives the boundary-to-bulk gain used to control the coefficient matrices, while the Dirichlet--Neumann estimates control the boundary operators in the surface equation.
Regularized changes of variables that straighten the free boundary are used in \cite[Section~3.1.1]{AlazardBurqZuily2014} and \cite[Section~3.1]{NguyenPausader2020}.
We recall the graph flattening and Sobolev estimates of \cite[Lemma~3.3]{LatoccaNguyen2025}.

\begin{lemma}[Lemma~3.3 in \cite{LatoccaNguyen2025}]\label{app:ln-flattening}
Let $q>1+d/2$, $f\in H^q(\R^d)$, and, in finite depth, let $b=-\ell_0+b_0$ with $b_0\in H^q$ and $f-b\geq h_*>0$.
For the flattening map \eqref{eq:explicit-flattening}, there exists $K=K(d,q)>0$ such that for any $\delta\in(0,1)$ satisfying
\[
 \delta\leq
 \begin{cases}
 (K\|f\|_{H^q}+1)^{-1},&\text{in infinite depth},\\
 h_*(K\|f\|_{H^q}+K\|b_0\|_{H^q}+1)^{-1},&\text{in finite depth},
 \end{cases}
\]
the map $\Phi_f$ is a Lipschitz diffeomorphism satisfying
\[
 \partial_z\varrho_f\geq\tfrac12\quad\text{in infinite depth},
 \qquad \partial_z\varrho_f\geq \tfrac12h_*\quad\text{in finite depth}.
\]
For every $r\geq0$, whenever the right-hand side is finite,
\begin{align*}
 \|\varrho_f-z\|_{H^r(\R^d\times(-\infty,0))}
 &\leq C_{d,r,\delta}\|f\|_{H^{r-1/2}},\\
 \|\varrho_f-\ell_0z\|_{H^r(\R^d\times(-1,0))}
 &\leq C_{d,r,\delta}
       (\|f\|_{H^{r-1/2}}+\|b_0\|_{H^{r-1/2}}).
\end{align*}
\end{lemma}

The Dirichlet--Neumann bound follows from Alazard--Burq--Zuily \cite[Theorem~3.12]{AlazardBurqZuily2014}; the paralinearization and tame remainder estimate follow from Nguyen--Pausader \cite[Theorem~3.18]{NguyenPausader2020}.

\begin{lemma}[Dirichlet--Neumann estimates]\label{app:ln-dn}
Let $q>1+d/2$, $1/2\leq r\leq q$, $f\in H^q(\R^d)$, and $h\in H^r(\R^d)$.
Then
\begin{equation*}
 \|\DN[f]h\|_{H^{r-1}} \leq C(\|f\|_{H^q})\|h\|_{H^r}.
\end{equation*}
Let $0<\theta\leq1/2$ with $\theta<q-1-d/2$.
If also $f\in H^{q+1/2-\theta}$, then
\begin{equation}\label{eq:app-ln-dn-para}
 \DN[f]h=T_{\lambda_f}(h-T_{B_h}f)-T_{W_h}\cdot\nabla f+\mathcal R[f]h,
\end{equation}
where $\lambda_f$ is defined in \eqref{eq:surface-symbols},
\[
 B_h=\frac{\nabla f\cdot\nabla h+\DN[f]h}{1+|\nabla f|^2}, \qquad W_h=\nabla h-B_h\nabla f,
\]
and the remainder satisfies the estimate
\begin{equation}\label{eq:app-ln-dn-remainder}
 \|\mathcal R[f]h\|_{H^{r-1/2}}
 \leq C(\|f\|_{H^q})
       (1+\|f\|_{H^{q+1/2-\theta}})\|h\|_{H^r}.
\end{equation}
The constants depend on the indices and dimension and, in finite depth, on the fixed bottom geometry and depth bound.
\end{lemma}

\subsection{Transport and background multipliers}

We use the following transport estimate for velocity fields tangent to the boundary of the reference domain.
The whole-space estimate is proved in \cite[Theorem~3.19]{BCD2011}.
The domain version follows from \cite[Theorem~2.8]{LatoccaNguyen2025} using extension and boundary tangency.
The subsequent estimates control the source generated by the variable background.

\begin{lemma}[Theorem~2.8 in \cite{LatoccaNguyen2025}]\label{app:ln-transport}
Let $r\geq1$, $r\ne1+D/2$, $h_0\in H^r(\mathcal S)$, $F\in L_T^1H^r(\mathcal S)$, and $V\in L_T^pL^2(\mathcal S)$ for some $p>1$.
Assume that $V$ is tangent to the boundary and that
\[
 \mathcal V_r(T):=
 \begin{cases}
 \|\nabla V\|_{L_T^1(H^{D/2}\cap L^\infty)},&r<1+D/2,\\
 \|\nabla V\|_{L_T^1H^{r-1}},&r>1+D/2
 \end{cases}
 <\infty.
\]
Then the transport problem
\[
 (\partial_t+V\cdot\nabla) h=F,\qquad h(0)=h_0,
\]
has a solution in $C_TH^r$, which is unique in $C_TH^1$ and satisfies
\[
 \|h\|_{C_TH^r} \leq C_r\left(\|h_0\|_{H^r}+\|F\|_{L_T^1H^r}\right)\exp\bigl(C_r\mathcal V_r(T)\bigr).
\]
\end{lemma}
At the top order, we apply this estimate with $r=s>D/2+1$ and $V\in L_T^2H^s$.

\begin{lemma}[Proposition~E.1 in \cite{LatoccaNguyen2025}]\label{app:ln-multiplier}
Let $J=(-1,0)$ or $(-\infty,0)$, let $a>0$, and $r\geq0$.
Write $\mathcal H^r=H^r(\R^d\times J)\cap L^\infty(\R^d\times J)$, with the sum norm.
If $F\in C_b^{\lceil r\rceil}(\R)$, then
\[
 \|F(h+az)v\|_{H^r} \leq C_{F,a,r}(\|h\|_{\mathcal H^r})\|v\|_{\mathcal H^r}.
\]
If $F'\in C_b^{\lceil r\rceil}(\R)$, then
\[
 \|(F(h_1+az)-F(h_2+az))w\|_{H^r} \leq C_{F,a,r}(\|h_1\|_{\mathcal H^r}+\|h_2\|_{\mathcal H^r}) \|(h_1-h_2)w\|_{\mathcal H^r}.
\]
Both estimates hold whenever all norms on the right-hand side are finite.
\end{lemma}

In our application, $F=\gamma'$, $a=1$ or $\ell_0$, and $h=\varrho_f-az$.
Thus, at indices at most $s$, the estimates require only finitely many bounded derivatives of $\gamma$.

\subsection{The zero-surface-tension local theory and surface comparison}

The background class of \cite[Definition~1.2]{LatoccaNguyen2025} consists of smooth profiles $\gamma$ for which both bounds in Lemma~\ref{app:ln-multiplier}, with $F=\gamma'$, hold at every nonnegative Sobolev index and for every $a>0$.
Section~\ref{sec:vanishing-surface-tension} explains how to use the result at a fixed index under the finite derivative assumptions of this paper.

\begin{lemma}[Theorem~1.4 in \cite{LatoccaNguyen2025}]\label{app:ln-local}
Let $q>(d+3)/2$, let $\gamma$ belong to the background class just defined, and assume
\[
 f_0\in H^q(\R^d),\qquad g_0\in H^q(\Omega_{f_0}),\qquad
 \inf_x\mathcal T(f_0)(x)\geq2a_*>0.
\]
In finite depth, assume $b=-\ell_0+b_0$, $\ell_0>0$, $b_0\in H^{q+1/2}$, and $f_0-b\geq2h_*>0$.
For each $R>0$, there are $\varepsilon(R)>0$ and $T(R)>0$ such that
$\|f_0\|_{H^q}\leq R$ and $\|g_0\|_{H^q(\Omega_{f_0})}\leq\varepsilon(R)$ yield a unique zero-surface-tension solution with
\[
 f\in C_TH^q\cap L_T^2H^{q+1/2},\qquad
 \widetilde g\in C_TH^q,
 \qquad \mathcal T(f)\geq a_*,
\]
and $f-b\geq h_*$ in finite depth.
Smallness may be measured either in the physical initial norm or in the pulled-back norm after adjusting $\varepsilon(R)$ by the bounded geometry equivalence constant.
The solution bounds and time depend only on the data bound, the lower bounds for stability and depth, and the fixed background and bottom.
For every $q'<q$, the pulled-back solution map is continuous into $(C_TH^{q'}\cap L_T^2H^{q'+1/2})\times C_TH^{q'}$.
\end{lemma}

\begin{lemma}[Proposition~5.4 in \cite{LatoccaNguyen2025}]\label{app:ln-surface-difference}
Let $q>(d+3)/2$, $T\le1$, and let two pairs
$(f_j,\widetilde g_j)$ satisfy the first equation
in \eqref{e:IPM} with $\sigma=0$ on $[0,T]$.
Assume
\[
 \sum_{j=1}^2\|f_j\|_{C_TH^q}\leq R,\qquad \sum_{j=1}^2\|f_j\|_{L_T^2H^{q+1/2}}\leq M, \qquad \mathcal T(f_j)\geq a_*>0,
\]
and a common depth bound $f_j-b\geq h_*>0$ in finite depth.
There is a threshold $\varepsilon_*(R,a_*,h_*)>0$ such that, if
$\sum_j\|\widetilde g_j\|_{C_TH^q}\leq\varepsilon_*$, then
\begin{equation*}
 \|f_1-f_2\|_{C_TH^{q-1}}+\|f_1-f_2\|_{L_T^2H^{q-1/2}}
 \leq C_{R,M,a_*,h_*} \left(\|f_1(0)-f_2(0)\|_{H^{q-1}}+\sqrt T\|\widetilde g_1-\widetilde g_2\|_{C_TH^{q-1}}\right).
\end{equation*}
The threshold and constant also depend on the indices and the fixed background and bottom.
\end{lemma}

In Section~\ref{sec:limit-comparison}, we include the capillary residual in the energy argument for this lemma and combine the resulting surface comparison with the elliptic and transport estimates for the density.

\end{document}